\documentclass{amsart}
\usepackage[dvips]{epsfig}
\usepackage{wrapfig}
\input cyrfont.def
\input cyracc.def
\newtheorem{theorem}{Theorem}[section]
\newtheorem{lemma}[theorem]{Lemma}

\newtheorem{proposition}[theorem]{Proposition}
\theoremstyle{definition}
\newtheorem{definition}[theorem]{Definition}
\theoremstyle{remark}
\newtheorem{remark}[theorem]{Remark}
\numberwithin{equation}{section}
\DeclareMathOperator{\spec}{spec}%
\DeclareMathOperator{\cov}{cov}

\begin{document}
\title[Rational Solutions of the
Schlesinger System]{Rational Solutions of the Schlesinger System
and Isoprincipal Deformations of Rational \\ Matrix Functions I.}
\author[V.Katsnelson]{Victor Katsnelson}
\address{Department of Mathematics\\
the Weizmann Institute of Science\\ Rehovot\\ 76100 \\
Israel} \email{victor.katsnelson@weizmann.ac.il}
\thanks{The research of Victor Katsnelson was supported
by the Minerva Foundation.}
\author[D. Volok]{Dan Volok}
\address{Department of Mathematics\\ Ben-Gurion University of the Negev\\
Beer-Sheva\\ 84105 \\ Israel}\email{volok@math.bgu.ac.il} 
\date{April 1, 2003}
\begin{abstract}{\small%
In the present paper we discuss the general facts  concerning the
Schlesinger system: the \(\tau\)-function, the local factorization
of solutions of Fuchsian equations and holomorphic deformations.
We introduce the terminology "isoprincipal" for the deformations
of Fuchsian equations with general (not necessarily non-resonant)
matrix coefficients  corresponding to solutions of the Schlesinger
system. Every isoprincipal deformation is isomonodromic.  The
converse is also true in the non-resonant case, but not in
general.

In the forthcoming sequel we shall give an explicit description of
a class of rational solutions of the Schlesinger system, based on
the techniques developed here, and the realization theory for
rational matrix functions.
}%
\end{abstract}
\maketitle  \noindent
NOTATIONS.\\
\begin{itemize}
\item
$\mathbb C$  stands for the complex plane;
\item
$\overline{\mathbb C}$  stands for the extended complex
 plane (\(=\) the Riemann sphere):\\ \(\overline{\mathbb
C}=\mathbb{C}\cup\infty\); \item ${\mathbb C}^n$  stands for the
\(n\) - dimensional complex space; \item in the coordinate
notation, a point \(\boldsymbol{t}\in{\mathbb C}^n\) is written as
\(\boldsymbol{t}=(t_1, \dots , t_n) ;\) \item
\(\mathbb{C}^n_{\ast}\) is the set of those points
\(\boldsymbol{t}=(t_1, \dots , t_n)\in{\mathbb C}^n,\) whose
coordinates \(t_1, \dots ,t_n\)  are pairwise different:
\({\mathbb{C}}^n_{\ast}={\mathbb C}^n\setminus\big(\bigcup_{1\leq
i,j\leq n,  i\not=j} \{\boldsymbol{t}:t_i=t_j\}\big);\) \item
\({\mathfrak{M}}_k\) stands for the set of all $k\times k$
matrices with complex entries; \item \([ . , . ]\) denotes the
commutator: for
 \(A, B\in\mathfrak{M}_k,\ [A, B]=AB-BA\);
 \item \(I\)  stands for the unity
matrix of the appropriate dimension;
\item
${\mathcal R}
({\mathfrak{M}}_k)$ stands for the set of all rational
${\mathfrak{M}}_k$-valued  functions $R$ such that
 $\mbox{\rm det}R(z)\not\equiv 0$;
\item
${\mathcal P}(R)$ stands for the set of all poles of the
function $R$, ${\mathcal N}(R)$ stands for the set of all poles of
the function $R^{-1}$; ${\mathcal P}(R)$ is said to be {\sf the
pole set of the function $R$}, ${\mathcal N}(R)$ is said to be
{\sf the zero set of the function $R$}.
\end{itemize}

\section{THE SCHLESINGER SYSTEM} \label{SCHLSEC}
 The Schlesinger system  is a
Pfaffian system of differential equations of the form:
\begin{equation}%
\label{Sch} \left\{
\begin{array}{llr}
\displaystyle\frac{\partial Q_{i}}{\partial t_{j}}=
\dfrac{\left[Q_{i},Q_{j}\right]}{t_{i}-t_{j}}, &
 1\leq
i,j\leq n,\ i\not=j,&\qquad (\ref{Sch}.\textup{a})
\\[3.0ex]
\dfrac{\partial Q_{i}}{\partial t_{i}}= -\displaystyle
\sum\limits_{\substack{1\leq j\leq n
\\
j\not=i}} \dfrac{\left[Q_{i},Q_{j}\right]} {t_{i}-t_{j}}, & \quad 1\leq i\leq
n,&\qquad (\ref{Sch}.\textup{b})
\end{array}\right.
\end{equation}
where
 \(Q_1,\ldots,Q_n\) are square \(k\times k\) matrix functions of
 \(\boldsymbol{t}=(t_1,\ldots,t_n),\
 \boldsymbol{t}\in\mathbb{C}^n_{\ast}\).\\
 The system of equations (\ref{Sch}) can be rewritten in a compact way as
 \begin{equation}
 {\mathfrak{d}}Q_{\nu}=\sum\limits_{\substack{\mu=1\\\mu\not=\nu}}^n
 \big[Q_{\mu}, Q_{\nu}\big] \mathfrak{d}\log (t_{\mu}-t_{\nu}), \quad
 \nu=1, \dots , n ,
 \label{SchC}%
 \end{equation}
 where \(\mathfrak{d}\) is the exterior differential.
 The Schlesinger system (\ref{Sch}) is over-determined.
 The compatibility (the integrability) condition for this system is of the form
 \begin{equation}%
 \mathfrak{d} \bigg(\sum\limits_{\substack{\mu=1\\\mu\not=\nu}}^n
 \big[Q_{\mu}, Q_{\nu}\big] \mathfrak{d} \log (t_{\mu}-t_{\nu})\bigg)=0 ,
 \quad \nu=1, \dots , n .
 \label{ClCond}%
 \end{equation}%
The condition (\ref{ClCond}) must be fulfilled  if
 \(Q_1(\boldsymbol{t}), \dots , Q_n(\boldsymbol{t})\)
 satisfy the equations (\ref{SchC}). This condition expresses
 the equality \(\mathfrak{d}\big(\mathfrak{d}Q_v)=0,\) which must be
  satisfied since \(\mathfrak{d}\cdot\mathfrak{d}=0\).

 Let us present the integrability condition (\ref{ClCond}) in a more explicit
 form.
 Since \(\mathfrak{d} \mathfrak{d}\log(t_{\mu}-t_{\nu})=0,\)
 \begin{multline}%
 \label{CE1}%
\mathfrak{ d} \bigg(\sum\limits_{\substack{\mu=1\\\mu\not=\nu}}^n
 \big[Q_{\mu}, Q_{\nu}\big] \mathfrak{d}\log (t_{\mu}-t_{\nu})\bigg)=\\
 =\sum\limits_{\substack{\mu=1\\\mu\not=\nu}}^n
 \big[\mathfrak{d}Q_{\mu}, Q_{\nu}\big] \wedge\mathfrak{d}\log (t_{\mu}-t_{\nu})+
 \sum\limits_{\substack{\mu=1\\\mu\not=\nu}}^n
 \big[Q_{\mu}, \mathfrak{d}Q_{\nu}\big] \wedge\mathfrak{d}\log (t_{\mu}-t_{\nu}),
 \end{multline}%
 where we have to substitute the expression of the form (\ref{SchC})
 for \(\mathfrak{d}Q_{\mu}\) and \(\mathfrak{d}Q_{\nu}\) in the right hand side of (\ref{CE1}).
 Performing the substitutions, we obtain the expressions
\begin{multline}%
 \label{CE2}%
\sum\limits_{\substack{\mu=1\\\mu\not=\nu}}^n
 \big[\sum\limits_{\substack{\lambda=1\\ \lambda\not=\mu}}^n
 \big[Q_{\lambda}, Q_{\mu}\big]
 , Q_{\nu}\big] \mathfrak{d}\log(t_{\lambda}-t_{\mu})\wedge
 \mathfrak{d}\log (t_{\mu}-t_{\nu})+\\
 +\sum\limits_{\substack{\mu=1\\\mu\not=\nu}}^n
 \big[Q_{\mu}, \sum\limits_{\substack{\lambda=1\\ \lambda\not=\nu}}^n
 \big[Q_{\lambda}, Q_{\nu}\big]\big]
\mathfrak{d}\log (t_{\lambda}-t_{\nu})\wedge
\mathfrak{d}\log (t_{\mu}-t_{\nu}) .
\end{multline}%
Since \(\omega\wedge\omega=0\) for every differential form
\(\omega\),
\(\mathfrak{d}\log(t_{\mu}-t_{\nu})\wedge\mathfrak{d}\log(t_{\mu}-t_{\nu})=0\).
Therefore, the integrability condition for the Schlesinger system
(\ref{SchC}) can be presented in the form
\begin{multline}%
\label{CE3}%
\sum\limits_{\substack{\lambda=1, \mu=1;\\
\lambda\not=\mu, \lambda\not=\nu,   \mu\not=\nu}}^n
\big[\big[Q_{\lambda}, Q_{\mu}\big], Q_{\nu}\big]
\mathfrak{d}\log(t_{\lambda}-t_{\mu})\wedge\mathfrak{d}\log (t_{\mu}-t_{\nu})+\\
+\sum\limits_{\substack{\lambda=1, \mu=1;\\
\lambda\not=\mu, \lambda\not=\nu,   \mu\not=\nu}}^n
\big[Q_{\mu},\big[Q_{\lambda},Q_{\nu}\big]\big]
\mathfrak{d}\log (t_{\lambda}-t_{\nu})\wedge
\mathfrak{d}\log (t_{\mu}-t_{\nu})=0,\quad \nu=1, \dots , n .
\end{multline}%

\begin{lemma}%
\label{LOI}%
For every holomorphic square matrix functions
\(Q_1, \dots , Q_n\), the integrability condition
\textup{(\ref{CE3})} holds.
\end{lemma}

\noindent The proof of this lemma is based on two identities. The
first of them is the Jakobi identity
\begin{equation}%
\label{JaId}%
[ [A ,B] ,C]+[ [B ,C] ,A]+[ [C ,A] ,B]=0,
\end{equation}%
which holds for arbitrary square matrices \(A, B, C\) (of the
same dimension). The second identity is an identity related to the
differential holomorphic 1-form \(\omega_{\lambda, \mu}\):
\begin{equation}%
\label{Fms}%
\omega_{\lambda, \mu}=\frac{dt_{\lambda}-dt_{\mu}}{t_{\lambda}-t_{\mu}},\quad
1\leq \lambda, \mu\leq n, \lambda\not=\mu .
\end{equation}%
Since
\begin{equation}%
\label{Ev}%
\omega_{\lambda,\mu}=\omega_{\mu,\lambda}, \quad 1\leq
\lambda, \mu\leq n, \lambda\not=\mu ,
\end{equation}%
there are \(\dfrac{n(n-1)}{2}\) different forms
\(\omega_{\lambda,\mu}\). The second identity is:
\begin{equation}%
\label{SeId}%
\omega_{\lambda,\mu}\wedge \omega_{\mu,\nu}+\omega_{\mu,\nu}\wedge
\omega_{\nu,\lambda}+\omega_{\nu,\lambda}\wedge
\omega_{\lambda,\mu}=0,\quad 1\leq\lambda, \mu, \nu\leq n,\ \ \
\lambda\not=\mu\not=\nu .
\end{equation}%
To prove this identity we observe that the left-hand side of
(\ref{SeId}) is the sum of the expressions that are obtained from
\(\omega_{\lambda,\mu}\wedge \omega_{\mu,\nu}\) by all cyclic
permutations of the indices \(\lambda, \mu, \nu\). Computing the
exterior product \(\omega_{\lambda,\mu}\wedge \omega_{\mu,\nu}\):
\begin{equation}%
\label{ExPr}%
\omega_{\lambda,\mu}\wedge
\omega_{\mu,\nu}=\frac{dt_{\lambda}\wedge dt_{\mu}+dt_{\mu}\wedge
dt_{\nu}+dt_{\nu}\wedge
dt_{\lambda}}{(t_{\lambda}-t_{\mu})(t_{\mu}-t_{\nu})} ,
\end{equation}
we see that it is represented as a fraction whose numerator is
invariant with respect to cyclic permutations of the indices
\(\lambda, \mu, \nu\). Therefore, the right-hand side of
(\ref{SeId}) is equal to
\begin{multline*}%
\left\{\frac{1}{(t_{\lambda}-t_{\mu})(t_{\mu}-t_{\nu})}+\frac{1}{(t_{\mu}-t_{\nu})(t_{\nu}-t_{\lambda})}+
\frac{1}{(t_{\nu}-t_{\lambda})(t_{\lambda}-t_{\mu})}\right\}\\[0.7ex]
\times\left(dt_{\lambda}\wedge dt_{\mu}+dt_{\mu}\wedge
dt_{\nu}+dt_{\nu}\wedge dt_{\lambda}\right)
\end{multline*}%
Evidently, the expression in the curly brackets vanishes. The
identity (\ref{SeId}) is proved.

\noindent
\begin{remark}
The forms \(\dfrac{1}{2\pi i} \omega_{\lambda,\mu},\ 1\leq
\lambda, \mu\leq n,\ \lambda\not=\mu,\) generate the cohomology
ring \(H^{\ast}(\mathbb{C}^n_{\ast}, \mathbb{Z})\) of the space
\(\mathbb{C}^n_{\ast}\). This result was established by
V. Arnold, \cite{Arn}.
\end{remark}

\noindent%
\textsf{PROOF of LEMMA \ref{LOI}}. %
Let us fix \(\nu\in 1, ,\dots , n .\) First, we transform the
first summand in (\ref{CE3}). Applying the Jakobi identity to the
factor \(\big[\big[Q_{\lambda}, Q_{\mu}\big], Q_{\nu}\big]\) in
the first sum, we rewrite (\ref{CE3}) as
\begin{multline}%
\label{CE4}%
\sum\limits_{\substack{\lambda=1, \mu=1;\\
\lambda\not=\mu, \lambda\not=\nu,   \mu\not=\nu}}^n
\big[\big[Q_{\mu}, Q_{\nu}\big], Q_{\lambda}\big]
\omega_{\lambda,\mu}\wedge\omega_{\mu,\nu}\\
+\sum\limits_{\substack{\lambda=1, \mu=1;\\
\lambda\not=\mu, \lambda\not=\nu,   \mu\not=\nu}}^n
\big[\big[Q_{\nu}, Q_{\lambda}\big], Q_{\mu}\big]
\omega_{\lambda,\mu}\wedge\omega_{\mu,\nu}\\
+\sum\limits_{\substack{\lambda=1, \mu=1;\\
\lambda\not=\mu, \lambda\not=\nu,   \mu\not=\nu}}^n
\big[\big[Q_{\lambda}, Q_{\nu}\big], Q_{\mu}\big]
\omega_{\lambda,\nu}\wedge\omega_{\mu,\nu}=0
\end{multline}
Since the summation runs over all pairs \(\lambda, \mu ,\) we can
interchange \(\lambda\) and \(\mu\) in the second and third sums
of (\ref{CE4}). Using also (\ref{Ev}) and the anti-commutativity
of the exterior product (\(
\omega\wedge\sigma=-\sigma\wedge\omega\) for every differential
1-forms  \(\omega \) and \(\sigma\)), we rewrite (\ref{CE4}) as
\begin{multline}%
\label{CE41}%
\sum\limits_{\substack{\lambda=1, \mu=1;\\
\lambda\not=\mu, \lambda\not=\nu,   \mu\not=\nu}}^n
\big[\big[Q_{\mu}, Q_{\nu}\big], Q_{\lambda}\big]
\omega_{\lambda,\mu}\wedge\omega_{\mu,\nu}\\
+\sum\limits_{\substack{\lambda=1, \mu=1;\\
\lambda\not=\mu, \lambda\not=\nu,   \mu\not=\nu}}^n
\big[\big[Q_{\mu}, Q_{\nu}\big], Q_{\lambda}\big]
\omega_{\nu,\lambda}\wedge\omega_{\lambda,\mu}\\
+\sum\limits_{\substack{\lambda=1, \mu=1;\\
\lambda\not=\mu, \lambda\not=\nu,   \mu\not=\nu}}^n
\big[\big[Q_{\mu}, Q_{\nu}\big], Q_{\lambda}\big]
\omega_{\mu,\nu}\wedge\omega_{\nu,\lambda}=0
\end{multline}
The last equality holds because its left-hand side is of the form
\begin{equation*}
\sum\limits_{\substack{\lambda=1, \mu=1;\\
\lambda\not=\mu, \lambda\not=\nu,   \mu\not=\nu}}^n
\big[\big[Q_{\mu}, Q_{\nu}\big], Q_{\lambda}\big]
\left\{\omega_{\lambda,\mu}\wedge\omega_{\mu,\nu}+\omega_{\nu,\lambda}\wedge
\omega_{\lambda,\mu}
+\omega_{\mu,\nu}\wedge\omega_{\nu,\lambda}\right\} ,
\end{equation*}
where the sums in the curly brackets vanish according to the
identity (\ref{SeId}). \\
\mbox{}\hfill  Q.E.D.\\[1ex]

\noindent From Lemma \ref{LOI} and from Frobenius theorem it
follows that the Cauchy problem for the Schlesinger system is
locally solvable for every initial data. More precisely, the
following result holds:

\begin{theorem}%
\label{LocExistTh}%
 Let
\(\boldsymbol{t}^0\in\boldsymbol{\mathbb{C}^n_*}\),
\(Q_1^0, \dots , Q_n^0\in{\mathfrak{M}}_k\).
Then there exist a neighborhood\ \ %
\(\mathcal{V}(\boldsymbol{t}^0, \delta)\) of the point
\(\boldsymbol{t}^0\),
\[\mathcal{V}(\boldsymbol{t}^0, \delta)=\{\boldsymbol{t}\in\mathbb{C}_n:
|t_1-t_1^0|<\delta, \dots , |t_n-t_n^0|<\delta\},
\mathcal{V}(\boldsymbol{t}^0, \delta)\subset\mathbb{C}^n_{\ast} ,\]
 and matrix functions \(Q_1(\boldsymbol{t}), \dots , Q_n(\boldsymbol{t})\),
\[Q_{\nu}(\boldsymbol{t}):
\mathcal{V}(\boldsymbol{t}^0, \delta)\mapsto \mathfrak{M}_k ,
\nu=1, \dots , n ,\] which are holomorphic in
\(\mathcal{V}(\boldsymbol{t}^0, \delta)\), satisfy the
Schlesinger system \textup{(\ref{Sch})} for
\(\boldsymbol{t}\in\mathcal{V}(\boldsymbol{t}^0, \delta)\), and
meet the initial condition
\begin{equation}%
\label{IniCond}%
Q_{\nu}(\boldsymbol{t}^0)=Q_{\nu}^0,  \nu=1, \dots , n
\end{equation}%
at the point \(\boldsymbol{t}^0\). Here \(\delta\) is a positive
number which depends on \(\boldsymbol{t}^0\) and \(Q_{\nu}^0,
\nu=1, \dots , n .\)
\end{theorem}%

\noindent Thus Theorem \ref{LocExistTh} is a local existence result. Actually, a
much more stronger global existence result holds for the Schlesinger system. It
turns out that the solution \(Q_1(\boldsymbol{t}), \dots ,  Q_n(\boldsymbol{t})\)
of the above introduced Cauchy problem  for this system exists \textsl{globally} as
meromorphic (multi-valued) functions in \(\mathbb{C}^n_{\ast}\). We do not discuss
this more general result in detail -- this will lead us too far. However, it is
worth mentioning  that every solution of the Schlesinger system can be represented
as a fraction of holomorphic functions in \(\mathbb{C}^n_{\ast}\), where some
complex valued holomorphic function \(\tau( . )\) on \(\mathbb{C}^n_{\ast}\)
appears in the denominator. This is the so-called \(\tau\)-function for the
considered solution. The singularities of this solution are precisely those
\(\boldsymbol{t}\in\mathbb{C}^n_{\ast}\) where \(\tau(\boldsymbol{t})=0.\)
\textsl{The "zeros" of the function \(\tau( . )\) are responsible for the
singularities of the solution \(Q_1( . ), \dots , Q_n( . )\)}.

\section{THE TAU-FUNCTION Of THE SCHLESINGER SYSTEM} \label{TauSec}
\begin{lemma} \label{CTF} Let matrix-functions
\(Q_1(\boldsymbol{t}), \dots , Q_n(\boldsymbol{t})\) satisfy
the Schlesinger system \textup{(\ref{Sch})}. Then the differential
1-form \(\omega\):
\begin{equation}%
\label{EF}%
\omega\stackrel{\textup{def}}{=}%
\frac{1}{2}\sum\limits_{\substack{\lambda, \mu=1\\
\lambda\not=\mu}}^n \textup{\large
tr} (Q_{\lambda}(\boldsymbol{t})Q_{\mu}(\boldsymbol{t}))
\mathfrak{d}\log (t_{\lambda}-t_{\mu}),
\end{equation}%
is closed:
\begin{equation}%
\label{ClF}%
\mathfrak{d}\omega=0 .
\end{equation}%
\end{lemma}

\noindent \textsf{PROOF.} It is clear that
\begin{eqnarray*}
\mathfrak{d}\omega&=&\frac{1}{2} \sum\limits_{\substack{\lambda, \mu=1\\
\lambda\not=\mu}}^n
\textup{tr} \Big(\big(\mathfrak{d}Q_{\lambda}(t)\big)Q_{\mu}(t)\Big)
 \mathfrak{d}\log(t_{\lambda}-t_{\mu})\\
&&+\frac{1}{2} \sum\limits_{\substack{\lambda, \mu=1\\
\lambda\not=\mu}}^n
\textup{tr} \Big(Q_{\lambda}(t)\big(\mathfrak{d}Q_{\mu}(t)\big)\Big)
 \mathfrak{d}\log(t_{\lambda}-t_{\mu}) .
\end{eqnarray*}
Interchanging the summation indices \(\lambda\) and \(\mu\) in the first sum
(and using the rule \(\textup{tr} AB=\textup{tr} BA\)), we obtain (after
renaming \(\lambda\to\nu, \nu\to\lambda)\) that
\begin{equation}%
\label{CC1}%
\mathfrak{d}\omega=\sum\limits_{\substack{\lambda, \nu=1\\
\lambda\not=\nu}}^n
\textup{tr} \Big(Q_{\lambda}(t)\big(\mathfrak{d}Q_{\nu}(t)\big)\Big)
 \mathfrak{d}\log(t_{\nu}-t_{\lambda}) .
\end{equation}%
Substituting the expression \textup(\ref{SchC}) for
\(\mathfrak{d}Q_{\nu}\) into \textup(\ref{CC1}), we obtain the
following expression for \(\mathfrak{d}\omega\):
\begin{equation}%
\label{CC2}%
\mathfrak{d}\omega=\sum\limits_{\substack{\lambda=1, \mu=1, \nu=1\\
\lambda\not=\mu, \lambda\not=\nu,   \mu\not=\nu}}^n
\textup{tr} \Big(
Q_{\lambda} \big[Q_{\mu},Q_{\nu}\big]\Big)\omega_{\mu,\nu}\wedge
\omega_{\nu,\lambda} ,
\end{equation}%
where \(\omega_{\lambda,\mu}\) is defined by (\ref{Fms}),  or
\begin{equation}%
\label{CC3}
\mathfrak{d}\omega=\sum\limits_{\substack{\lambda=1, \mu=1, \nu=1\\
\lambda\not=\mu, \lambda\not=\nu,   \mu\not=\nu}}^n
\textup{tr} \Big(
Q_{\lambda}Q_{\mu}Q_{\nu}-Q_{\lambda}Q_{\nu}Q_{\mu}\Big)
  \omega_{\mu,\nu}\wedge\omega_{\nu,\lambda} .
\end{equation}
The traces in the last formula are invariant under a circular
permutation of the three indices \(\lambda, \mu, \nu\). So we
can replace \(\omega_{\mu,\nu}\wedge\omega_{\nu,\lambda}\) by
\(\frac{1}{3}\big(\omega_{\mu,\nu}\wedge\omega_{\nu,\lambda}+
\omega_{\nu,\lambda}\wedge\omega_{\lambda,\mu}+\omega_{\lambda,\mu}\wedge\omega_{\mu,\nu
}\big).\) But this quantity vanishes identically according to the
identity (\ref{SeId}). Thus, \(\mathfrak{d}\omega=0\). \hfill
Q.E.D.\\[1ex]

\noindent Since a closed form \(\mathfrak{d}\omega\) is locally
exact, the following theorem is an immediate consequence of Lemma
\ref{CTF}:

\begin{theorem}%
\label{ETF}%
Let holomorphic matrix-functions
\(Q_1(\boldsymbol{t}), \dots , Q_n(\boldsymbol{t})\) satisfy
the \\ Schlesinger system \textup{(\ref{Sch})} in some simply
connected domain \(\mathcal{D},
\mathcal{D}\subset\mathbb{C}^n_{\ast}\) . Then there exists a
holomorphic non-vanishing function \(\tau(\boldsymbol{t})\) in
\(\mathcal{D}\) such that
\begin{equation}%
\label{TFR}%
\mathfrak{d}\log \tau(\boldsymbol{t})=\frac{1}{2}\sum\limits_{\substack{\lambda, \mu=1\\
\lambda\not=\mu}}^n \textup{\large
tr} (Q_{\lambda}(\boldsymbol{t})Q_{\mu}(\boldsymbol{t})) \mathfrak{d}\log
(t_{\lambda}-t_{\mu})
\end{equation}%
Such a function \(\tau\) is unique up to a constant factor: if
\(\tau_{1}(\boldsymbol{t})\) and \(\tau_{2}(\boldsymbol{t})\) are two functions
satisfying the relation \textup{(\ref{TFR})}, then
\(\tau_{2}(\boldsymbol{t})=c \tau_{1}(\boldsymbol{t})\), where \(c\not=0\) is a
constant.
\end{theorem}%

\begin{definition}%
\label{DTF}%
A function \(\tau\) related to a solution \(Q_1, \dots , Q_n\) of the
Schlesinger system \textup{(\ref{Sch})} by the equality \textup{(\ref{TFR})} is
said to be \textsf{the \(\tau\)-function for this solution.}
\end{definition}%

\noindent The exactness of the differential form
\textup{(\ref{EF})}, stated in our Lemma \ref{CTF}, was
established by M. Sato, T. Miwa and M. Jimbo in \cite{SMJ2}
(see \S 2.4 there). In \cite{SMJ2}, the potential
\textup{(\ref{TFR})} was introduced as well as the term
"\(\tau\)-function" for this potential.

\begin{remark}%
\label{DHTF}%
>From Lemma \ref{CTF} it follows that the \(\tau\)-function for a
solution of the Schlesinger equation is a holomorphic
non-vanishing (multi-valued) function for those \(\boldsymbol{t}\)
where the solution is holomorphic. Actually, the \(\tau\)-function
is holomorphic on the \textsf{whole}   \(\mathbb{C}^n_{\ast}\). It
vanishes for those \(\boldsymbol{t}\) where the solution is
singular. (See the last paragraph of the section \ref{SCHLSEC}).
In a very general setting, the holomorphy of \(\tau\)-functions
was proved in \textup{\cite{Miw}}.
\end{remark}%

\begin{remark}%
\label{RTF}%
The \(\tau\)-function is a very general concept related to
non-linear differential equations. This function was first
discovered by R. Hirota as a tool for generating many-soliton
solutions. (References to the works of H. Hirota can be found in
bibliography of the book \textup{\cite{New}}). Let us give
examples of \(\tau\)-functions for the Korteveg - deVries
equation
\begin{equation}%
\label{KdV}%
u_t-6uu_x+u_{xxx}=0 .
\end{equation}%
\(N\) - soliton solution
\begin{equation}%
\label{NSol}%
u(x, t)=-2\frac{\partial^2}{\partial x^2}\log \Delta(x, t) ,
\end{equation}%
\begin{equation}%
\label{NDet}%
\Delta(x, t)=\det
\bigg(\delta_{jk}+\frac{c_jc_k}{\eta_j+\eta_k}\exp[-(\eta_j^3+\eta_k^3)t
-(\eta_j+\eta_k)x]\bigg)_{j,k=1, \dots , N} ;
\end{equation}%
\(N\) - phase solution
\begin{equation}%
\label{NPhS}%
u(x, t)=-2\frac{\partial^2}{\partial x^2}\log
\vartheta(\boldsymbol{a}x+\boldsymbol{b}t+\boldsymbol{c}) ,
\end{equation}%
where \(\vartheta\) is the Riemann theta function;\\
rapidly decreasing solution
\begin{equation}%
\label{RDeS}%
u(x, t)=-2\frac{\partial^2}{\partial x^2}\log \Delta_{[x, \infty)}(t) ,
\end{equation}%
with \(\Delta_{[x, \infty)}(t)\) the Fredholm determinant of the
Gelfand - Levitan equation. In these examples the
\(\tau\)-function is
\(\Delta(x, t), \vartheta(\boldsymbol{a}x+\boldsymbol{b}t+\boldsymbol{c}),
\Delta_{[x, \infty)}(t)\), respectively.

In the paper \textup{\cite{WNTB}}, \(\tau\)-function appeared as
the spin - spin correlation function for the two - dimensional
Ising model. In \textup{\cite{JMMS}}, the Fredholm determinant
\(\det(\boldsymbol{I}-\lambda \boldsymbol{K})\) of the integral
operator \(\boldsymbol {K}\) with the kernel \(K(x, x^{\prime})=
\frac{\sin (x-x^{\prime})}{x-x^{\prime}}\) on the union of
intervals \(\boldsymbol{\cup}_{1\leq j\leq n} [a_{2j-1}, a_{2j}],
\ a_1<\dots<a_{2n},\) is expressed as the \(\tau\)-function
\(\tau(a_1, \dots , a_{2n})\) of the Schlesinger system
constructed from the operator. (See in particular section 7 of
\textup{\cite{JMMS}}, where the deformation theory of the Fredholm
integral operator with the kernel \(\frac{\sin
(x-x^{\prime})}{x-x^{\prime}}\) was related to the deformation
theory of linear differential equations). The interrelation
between the deformation theory of linear differential operators,
the Schlesinger system and its \(\tau\)-function, correlation
function in statistical physics theory, quantum fields and
Fredholm determinants are covered in the expository survey
\textup{\cite{SMJ1}}. The \(\tau\)-function for general matrix
linear differential equation with rational coefficients was
considered in \textup{\cite{JMU}}, and was studied in great detail
in \textup{\cite{JM1}, \cite{JM2}}.

However, the importance of the \(\tau\)-function and its central
role in the soliton theory was not understood before the works of
the research group from Kioto (M. Sato, T. Miwa, M. Jimbo, E.
Date, M. Kashiwara, Y. Sato). The work of this group, surrounding
Mikio Sato, reached its peak in the years around 1981. In the
paper \textup{\cite{SaSa}}, the time evolution of a solution of a
non-linear equation is interpreted as the dynamical motion of a
point of the infinite dimensional Grassman manifold. In the papers
of E. Date, M. Jimbo, M. Kashiwara and T. Miwa, Lie algebras and
vertex operators were incorporated into the soliton theory. (See
the series of seven notes \textup{\cite{DJKM1}}, and the paper
\textup{\cite{DJKM2}}). In this setting, the \(\tau\)-function
appears in a very general way. The paper \textup{\cite{Sa}}
explains some ideas from \textup{\cite{SaSa}} in more detail. The
further development of the approach related to the consideration
of soliton equations as dynamic systems on an infinite dimensional
Grassman manifold is done in \textup{\cite{SeWi}.} In particular,
in section 3 of \textup{\cite{SeWi}} the geometric meaning of the
\(\tau\)-function is explained: this function appears as an
infinite determinant related to the determinant bundle over the
Grassman manifold. In \textup{\cite{MSW}}, the \(\tau\)-function
is defined as the generalized cross-ratio of four points of the
Grassman manifold. Tau function is discussed in Chapters 7 and 8
 of the book \cite{Dick}.
The book \textup{\cite{MJD}} is a presentation of the algebraic
analysis of the nonlinear differential equations of the KdV type
written for students. In particular, it contains a very clear
presentation of the circle of problem related to
\(\tau\)-functions. Chapter 4 of the book \textup{\cite{New}} is
dedicated to the \(\tau\)-function and its relation to the
Painlev{\'e} property and to the B\"acklund transformations. This
non-formally written book covers a very widespread circle of
problems.
\end{remark}%

\noindent In the present work we  construct a class of rational
solutions of the Schlesinger system (\ref{Sch}). For such a
solution, we calculate its \(\tau\)-function explicitly. This
\(\tau\)-function is equal to the determinant of some (finite)
matrix which is intimately related to this solution. In the
considered special case, \(\tau\)-function is a rational function
of the variables \(t_1 , \dots , t_n\).

\section{THE POTENTIAL FUNCTION FOR A SOLUTION OF THE SCHLESINGER SYSTEM.\label{PotFunc}}
From (\ref{Sch}.b) it follows that if the matrix functions \(Q_1(\boldsymbol{t}),
\dots , Q_n(\boldsymbol{t})\) satisfy the Schlesinger system in a domain
\(\mathcal{G}\in\mathbb{C}^n_{\ast}\), then
\begin{equation}%
\label{PolnDif}%
\frac{\partial Q_i(\boldsymbol{t})}{\partial t_j}\equiv
\frac{\partial Q_j(\boldsymbol{t})}{\partial t_i} .
\end{equation}%

Hence, there exists a matrix function \(V(\boldsymbol{t})\) in
\(\mathcal{G}\) such that
\begin{equation}%
\label{Grad}%
Q_j(\boldsymbol{t})= \frac{\partial V(\boldsymbol{t})}{\partial
t_j},\ \ \boldsymbol{t}\in\mathcal{G},\quad j=1, \dots , n .
\end{equation}%

\begin{remark}%
\label{MultPot} If the domain \(\mathcal{G}\) is multiconnected,
the function \(V(\boldsymbol{t})\) can be multivalued. In this
case we can consider \(V(\boldsymbol{t})\) as single-valued
function on the universal covering \(\cov (\mathcal{G},
\boldsymbol{t}_0)\) of the domain \(\mathcal{G}\), where
\(\boldsymbol{t}_0\) is a chosen distinguished point of
\(\mathcal{G}\).
\end{remark}%
\begin{definition}
The matrix function \(V(\boldsymbol{t})\) which satisfies the
condition \eqref{Grad} is said to be the \textsf{potential matrix
function} for the solution
\(Q_1(\boldsymbol{t}), \dots , Q_n(\boldsymbol{t})\) of the
Schlesinger system.
\end{definition}

We have proved that the potential matrix function exists for any
solution of the Schlesinger system. It is clear that the potential
function is unique up to an arbitrary constant summand.

\begin{lemma}%
\label{TPOT1}%
Let matrices
\(Q_1(\boldsymbol{t}), \dots , Q_n(\boldsymbol{t})\) satisfy
the Schlesinger system \eqref{Sch}, and moreover the condition
\begin{equation}%
\label{ZCSS}%
\sum\limits_{1\leq i\leq n}Q_i(\boldsymbol{t})=0
\end{equation}%
for every \(\boldsymbol{t}\in\mathcal{G}\).

Then the function
\begin{equation}%
\label{EFPF}%
V(\boldsymbol{t})\stackrel{\textup{\tiny
def}}{=}\sum\limits_{1\leq i\leq n}t_i\cdot Q_i(\boldsymbol{t}),\
\ \boldsymbol{t}\in\mathcal{G} ,
\end{equation}%
is a potential matrix function for the solution
\(Q_1(\boldsymbol{t}), \dots , Q_n(\boldsymbol{t})\) .
\end{lemma}%

\noindent%
\textsf{PROOF}. Differentiating the expression on the left hand
side of \eqref{EFPF}, we obtain
\begin{equation}%
\label{Dif1}%
\frac{\partial V(\boldsymbol{t})}{\partial
t_j}=Q_j(\boldsymbol{t})+\sum\limits_{1\leq i\leq n}t_i\cdot
\frac{\partial Q_i(\boldsymbol{t})}{\partial t_j} .
\end{equation}%
Thus, we have to prove that
\begin{equation}%
\label{Dif2}%
\sum\limits_{1\leq i\leq n}t_i\cdot \frac{\partial
Q_i(\boldsymbol{t})}{\partial t_j}=0 .
\end{equation}%
 Substituting the expression
\eqref{Sch} for \(\frac{\partial Q_i(\boldsymbol{t})}{\partial
t_j}\) into the left hand side of \eqref{Dif2}, we obtain
\begin{multline}%
\label{Dif3}%
\sum\limits_{1\leq i\leq n}t_i\cdot \frac{\partial
Q_i(\boldsymbol{t})}{\partial t_j}=\sum\limits_{\substack{1\leq
i\leq
n,\\i\not=j}}t_i\cdot\frac{[Q_i(\boldsymbol{t}), Q_j(\boldsymbol{t})]}{t_i-t_j}
-t_j\cdot \sum\limits_{\substack{1\leq i\leq
n,\\i\not=j}}\frac{[Q_i(\boldsymbol{t}), Q_j(\boldsymbol{t})]}{t_i-t_j}
\\[1.0ex]
=\sum\limits_{\substack{1\leq i\leq
n,\\i\not=j}}[Q_i(\boldsymbol{t}), Q_j(\boldsymbol{t})]=
\sum\limits_{\substack{1\leq
i\leq n}}[Q_i(\boldsymbol{t}), Q_j(\boldsymbol{t}]=0 .
\end{multline}%
\large
Thus, \eqref{Dif2} holds.%
\hfill Q.E.D.

\begin{lemma}
\label{TPOT2}%
Let \(\mathcal{G}\) be a domain in \(\mathbb{C}^n_{\ast}\), and
let matrix functions
\(Q_1(\boldsymbol{t}), \dots ,\) \(Q_n(\boldsymbol{t})\) satisfy
the Schlesinger system \eqref{Sch} for
\(\boldsymbol{t}\in\mathcal{G}\). Assume that the equality
\eqref{ZCSS} holds at least at one point
\(\boldsymbol{t}=\boldsymbol{t}_0\in\mathcal{G}\). Then this
equality  holds at every point \(\boldsymbol{t}\in\mathcal{G}\).
\end{lemma}

\noindent%
\textsf{PROOF}. The equality \eqref{ZCSS} will be proved if we
prove that
\begin{equation}%
\label{ZerDer}%
\frac{\partial {\ }}{\partial t_j} \bigg(\sum\limits_{1\leq i\leq
n}Q_i(\boldsymbol{t})\bigg)=0\ \ \text{for}\ \
\boldsymbol{t}\in\mathcal{G} .
\end{equation}%
>From \eqref{Sch} it follows that
\begin{multline*}
\sum\limits_{1\leq i\leq n}\frac{\partial {Q_i(\boldsymbol{t})
}}{\partial t_j}=\frac{\partial {Q_j(\boldsymbol{t}) }}{\partial
t_j} + \sum\limits_{\substack{1\leq i\leq
n\\i\not=j}}\frac{\partial {Q_i(\boldsymbol{t}) }}{\partial
t_j}\\
=-\sum\limits_{\substack{1\leq i\leq
n\\i\not=j}}\frac{[Q_i(\boldsymbol{t}),Q_j(\boldsymbol{t})]}{t_i-t_j}+
\sum\limits_{\substack{1\leq i\leq
n\\i\not=j}}\frac{[Q_i(\boldsymbol{t}),Q_j(\boldsymbol{t})]}{t_i-t_j}=0
\ \ \textup{for} \ \ \boldsymbol{t}\in\mathcal{G} .
\end{multline*}
Thus, \eqref{ZerDer} holds. \hfill Q.E.D.

The following lemma is an immediate consequence of two previous
ones.

\begin{lemma}%
\label{TPOT3}%
Let the matrix functions
\(Q_1(\boldsymbol{t}), \dots , Q_n(\boldsymbol{t})\) satisfy
the Schle\-singer system \eqref{Sch} for every
\(\boldsymbol{t}\in\mathcal{G}\), where \(\mathcal{G}\) is a
domain in \(\mathbb{C}^n_{\ast}\). Assume that the condition
\eqref{ZCSS} holds at least for one point
\(\boldsymbol{t}=\boldsymbol{t}_0\in\mathcal{G}\).

Then the function \(V(\boldsymbol{t})\) which is defined by
\eqref{EFPF} is a potential matrix function for the solution
\(Q_1(\boldsymbol{t}), \dots , Q_n(\boldsymbol{t})\) .
\end{lemma}

The expression \(\sum t_iQ_i(\boldsymbol{t})\),
which appears on the right hand side
of the equality \eqref{EFPF}, has a useful interpretation in
terms of an auxiliary ordinary differential equation. Let
\(\boldsymbol{t}=(t_1, \dots , t_n)\in\mathbb{C}^n_{\ast}\),
and \(Q_1, \dots , Q_n\) be \(k\times k\) matrices. (For the
moment, we consider the point \(\boldsymbol{t}\) and the matrices
\(Q_1, \dots , Q_n\) as fixed rather than variable.) Let us
introduce an auxiliary complex variable \(x\) belonging to the
domain %
\footnote{ We recall that
\(\overline{\mathbb{C}}=\mathbb{C}\cup\infty\) is the extended
complex plane.}
\(\overline{\mathbb{C}}\setminus\{t_1, \dots , t_n\}\).
Consider the linear differential equation
\begin{equation}%
\label{AuxDE}%
\frac{dY(x)}{dx}=\bigg(\sum\limits_{1\leq i\leq
n}\frac{Q_i}{x-t_i}\bigg)\cdot Y(x)%
\end{equation}%
with respect to a \(k\times k\) matrix function \(Y(x)\) of the
variable \(x\). The coefficient matrix function
\(\sum\frac{Q_i}{x-t_i}\) of
the equation \eqref{AuxDE} is holomorphic in the domain
\(\overline{\mathbb{C}}\setminus\{t_1, \dots , t_n\}\), and if
the condition
\begin{equation}%
\label{ZAI}%
\sum\limits_{1\leq i\leq n}Q_i=0
\end{equation}%
holds, then the residue of the coefficient matrix function at the
point \(\infty\) vanishes. Thus, under the condition \eqref{ZAI},
the Cauchy problem  with the initial condition
\begin{equation}%
\label{IniCondInf}%
 Y(\infty)=I
\end{equation}%
for the differential equation \eqref{AuxDE} in the domain
\(\overline{\mathbb{C}}\setminus\{t_1, \dots , t_n\}\) is
solvable. The solution \(Y(x)\) of this Cauchy problem is a
function holomorphic in the disk \(R<|x|\leq\infty\), where
\(R=\max_{j}|t_j|\). The function \(Y(x)\) admits a
Laurent expansion in this disk:
\begin{equation}%
\label{LaurExpa}%
Y(x)=I+\sum\limits_{1\leq_l\leq \infty}\frac{Y_{-l}}{x^l},\quad
R<|x|\leq\infty .
\end{equation}%

\begin{lemma}%
\label{IFLCS}%
Under the condition \eqref{ZAI}, the Laurent coefficient
\(Y_{-1}\)  of the expansion \eqref{LaurExpa} for the solution
\(Y(x)\) of the Cauchy problem
\eqref{AuxDE} - \eqref{IniCondInf} can be expressed as
\begin{equation}%
\label{EFRI} Y_{-1}=- \sum\limits_{1\leq i\leq n}t_iQ_i .
\end{equation}%
\end{lemma}%

\noindent%
\textsf{PROOF}. Under the condition \eqref{ZAI}, the Laurent
expansion at \(\infty\) of the coefficient matrix function
\(\sum\frac{Q_i}{x-t_i}\)   is
of the form
 \begin{equation}%
 \label{LECMF}%
 \sum\limits_{1\leq i\leq
n}\frac{Q_i}{x-t_i}=\bigg(\sum\limits_{1\leq i\leq
n}t_iQ_i\bigg)\cdot x^{-2}+o(|x|^{-2})\ \ \textup{as}\ \
x\to\infty .
\end{equation}%
Substituting the expansions \eqref{LaurExpa} and \eqref{LECMF}
into the differential equation \eqref{AuxDE} and comparing the
Laurent coefficients, we obtain \eqref{EFRI}. \\ %
\hspace*{0.0ex} \hfill Q.E.D.

\begin{theorem}%
Let \(k\times k\) matrix functions \
\(Q_1(\boldsymbol{t}), \dots , Q_n(\boldsymbol{t})\) \ satisfy
the  Schlesinger system for \(\boldsymbol{t}\in\mathcal{G}\),
where \(\mathcal{G}\) is a domain in \(\mathbb{C}^n_{\ast}\). Let
the equality \eqref{ZCSS} hold at least at one point
\(\boldsymbol{t}=\boldsymbol{t}_0\in\mathcal{G}\), or, what is the
same \footnote{Lemma \ref{TPOT2}}, at every point
\(\boldsymbol{t}\in\mathcal{G}\). Assume that for each fixed
\(\boldsymbol{t}\in\mathcal{G}\), \(k\times k\) matrix function
\(Y(x,\boldsymbol{t})\) satisfies the linear differential equation
\begin{equation}%
\label{AuxLDEq}%
\frac{dY(x,\boldsymbol{t})}{dx}=\bigg(\sum\limits_{1\leq i\leq
n}\frac{Q_i(\boldsymbol{t})}{x-t_i}\bigg)\cdot
Y(x,\boldsymbol{t}),\quad
x\in\mathbb{C}\setminus\{t_1, \dots , t_n\},
\end{equation}%
with respect to \(x\)  and the initial condition \footnote{Due to
\eqref{ZCSS}, the value \(Y(x,\boldsymbol{t})_{|x=\infty}\)
exists.}
\begin{equation}%
\label{IniCondi}%
Y(x,\boldsymbol{t})_{|x=\infty}=I.
\end{equation}%
Then the coefficient \footnote{Here we redenote the Laurent
coefficient \(Y_{-1}(\boldsymbol{t})\) by \(V(\boldsymbol{t})\).}
\(V(\boldsymbol{t})\) of the Laurent expansion
\begin{equation}%
\label{LEDPT}%
Y(x,\boldsymbol{t})=I-\frac{V(\boldsymbol{t})}{x}+\frac{Y_{-2}(\boldsymbol{t})}{x^2}+
\frac{Y_{-3}(\boldsymbol{t})}{x^2}+\cdots ,\quad
(\max\limits_i|t_i|)<|x|\leq\infty ,
\end{equation}%
is the potential matrix function for the solution
\(Q_1(\boldsymbol{t}), \dots , Q_n(\boldsymbol{t})\), i.e. the
equalities \eqref{Grad} hold.
\end{theorem}%

\noindent%
\textsf{PROOF}. This theorem is an immediate consequence of Lemmas
\ref{TPOT1}  and \ref{IFLCS}. \hspace*{0.0ex}\hfill Q.E.D.

The differential equation \eqref{AuxLDEq} plays an important role
in what follows.

\section{THE SCHLESINGER SYSTEM AS THE COMPATIBILITY
CONDITION FOR SOME OVER-DETERMINED SYSTEM OF LINEAR DIFFERENTIAL
EQUATIONS.\label{CompCond}}
 The result stated in this section is an important step
in the study of the nonlinear Schlesinger system (\ref{Sch}) by
linear methods.

\begin{theorem}
\label{SEC} Let
\(Q_1(\boldsymbol{t}), \dots, Q_n(\boldsymbol{t})\) be \(k\times
k\) matrix functions,  holomorphic with respect to
\(\boldsymbol{t}= (t_1, \dots , t_n)\) in some open domain
\(\mathcal{G}\), \(\mathcal{G}\subseteq\mathbb{C}^n_{\ast}\).

Assume that the following system of linear differential equations
\begin{subequations}%
\label{ODS}%
\begin{align}
\frac{\partial Y}{\partial x}=\bigg(\sum\limits_{1\leq j\leq
n}\frac{Q_j(\boldsymbol{t})}{x-t_j}\bigg) Y , \hspace{8.5ex}
\label{ODS1}\\[1.0ex]
 \frac{\partial Y}{\partial t_k}=-\frac{Q_k(\boldsymbol{t})}%
{x-t_k} Y\ \ \ (k=1, \dots , n)
 \label{ODS2}
\end{align}%
\end{subequations}%
for a square matrix function \(Y\) is compatible, that is, there
exists a \(k\times k\) matrix function \(Y= Y(x;\boldsymbol{t}),\)
which is holomorphic, satisfies the equations \textup{(\ref{ODS})}
for \(x\in\mathbb{C}\), \(\boldsymbol{t}\in\mathcal{G}\),
\(x\not=t_1, \dots , t_n\), and, moreover, is not degenerate:
\begin{equation}%
\label{NDG}%
\det Y(x;\boldsymbol{t})\not=0 , \ \ \text{for}\ \ \
x\in\mathbb{C}, \ \boldsymbol{t}\in\mathcal{G} , \
x\not=t_1, \dots , x\not=t_n .
\end{equation}%
Then the matrix functions
\(Q_1(\boldsymbol{t}), \dots , Q_n(\boldsymbol{t})\) satisfy
the Schlesinger system \textup{(\ref{Sch})} for
\(\boldsymbol{t}\in\mathcal{G}\).
\end{theorem}

\noindent \textsf{PROOF.} Let us exploit the equalities
\begin{equation}
\label{IO1}%
\frac{\partial }{\partial t_k} \frac{\partial Y}{\partial
x}=\frac{\partial }{\partial x} \frac{\partial Y}{\partial
t_k},\quad 1\leq k\leq n .
\end{equation}
Differentiating  equation(\ref{ODS1}) with respect to \(t_k\), we
obtain
\begin{equation*}
\label{PD}%
\frac{\partial }{\partial t_k} \frac{\partial Y}{\partial x}=
\left(\sum\limits_{1\leq j\leq n}
\frac{Q_j}{x-t_j}\right) \frac{\partial Y}{\partial t_k}+
\left(\sum\limits_{1\leq j\leq n}\frac{\frac{\partial
Q_j}{\partial t_k}}{x-t_j}\right) Y+\frac{Q_k}{(x-t_k)^2} Y .
\end{equation*}
Substituting the expression for \(\frac{\partial Y}{\partial
t_k}\) from (\ref{ODS2}) into the last equality, we get
\begin{equation}
\label{ESD1}%
\frac{\partial{\ }}{\partial t_k}\frac{\partial Y}{\partial x}=
-\sum\limits_{\substack{1\leq j\leq
n}}\frac{Q_j Q_k}{(x-t_j)(x-t_k)} Y +\frac{Q_k}{(x-t_k)^2} Y + \sum\limits_{1\leq
j\leq n}\frac{\dfrac{\partial Q_j}{\partial t_k}}{x-t_j} Y .
\end{equation}
Differentiating  equation(\ref{ODS2}) with respect to \(x\), we
obtain
\begin{equation*}%
\frac{\partial }{\partial x} \frac{\partial Y}{\partial t_k}=
\frac{Q_k}{(x-t_k)^2} Y - \frac{Q_k}{x-t_k} \frac{\partial
Y}{\partial x} .
\end{equation*}%
Substituting the expression for \(\frac{\partial Y}{\partial x}\)
from (\ref{ODS1}) into the last equality, we get
\begin{equation}%
\label{ESD2}%
\frac{\partial }{\partial x} \frac{\partial Y}{\partial t_k}=
\frac{Q_k}{(x-t_k)^2} Y -\sum\limits_{\substack{1\leq j\leq
n}}\frac{Q_k Q_j}{(x-t_j)(x-t_k)} Y
\end{equation}%
Thus,  equation (\ref{IO1}) takes the form
\begin{equation}%
\label{ESD3}%
\sum\limits_{1\leq l\leq n}\frac{\dfrac{\partial Q_l}{\partial
t_k}}{x-t_l}=\sum\limits_{\substack{1\leq l\leq
n\\[0.2ex]l\not=k}}\frac{\big[ Q_l , Q_k \big]}{(x-t_l)(x-t_k)}
\end{equation}%
(Since matrix \(Y\) is non-degenerate, we can cancel out the
factor \(Y\).)

The numerators in (\ref{ESD3}) do not depend on \(x\). Thus,
(\ref{ESD3}) expresses the equality of two rational functions.
Comparing the residues of both sides of (\ref{ESD3}) at \(x=t_j,
j\not=k\), we obtain (\ref{Sch}.a). Comparing the residues at
\(x=t_k\), we obtain (\ref{Sch}.b).
\hfill  Q.E.D.\\[1ex]

\begin{remark}%
\label{RI}%
In the proof of Theorem \textup{\ref{SEC}} we have not exploited
the equalities
\[\frac{\partial }{\partial t_k} \frac{\partial Y}{\partial
t_j}=\frac{\partial }{\partial t_j} \frac{\partial Y}{\partial
t_k},\quad 1\leq k,j\leq n .\] They lead to (\ref{Sch}.a) and
 provide no additional information on \(Q_1,\ldots,Q_n.\)
Thus,  reversing the reasoning in the proof, we can demonstrate
that the Schlesinger system \textup{(\ref{Sch})} is not only a
necessary, but also a sufficient condition for the solvability of
the linear system \textup{(\ref{ODS1})} - \textup{(\ref{ODS2})}
for a non-degenerate matrix function \(Y\).
\end{remark}%

\begin{remark}%
\label{CW}%
The result, which we formulated as Theorem \ref{SEC}, is contained
in the paper \textup{\cite{Sch1}} of L. Schlesinger.
\textup{(}See Section II of this paper, especially p. 294. See
also p. 324 of the book \textup{\cite{Sch2}}\textup{)}. This was
the paper \textup{\cite{Sch1}} where the Schlesinger equations
\textup{(\ref{Sch})} originally appeared. In \textup{\cite{Sch1}},
the equations \textup{(\ref{Sch})} were derived as the
compatibility condition of some over-determined system of linear
differential equations. In \textup{\cite{Gar}}, R. Garnier has
derived the Painlev\'{e} equation \(P_{VI}\) from the
compatibility condition of the appropriate over-determined system
of linear differential equations. \textup{(}This over-determined
system itself appeared in the paper \textup{\cite{FuR}} of
R. Fuchs, which, in its turn, was based on works of his father
L.\.Fuchs published in Sitzungsberichte der Berliner Akademie der
Wissenschaften in 1888 - 1898\textup{)}. Since the publication
of the works \textup{\cite{ZaSh}} and \textup{\cite{AKNS}},
representing a solution of a nonlinear differential equation as
the compatibility condition for some over-determined system of
linear differential equations has become a common wisdom. In
\textup{\cite{AKNS}}, under the strong influence of the paper
\textup{\cite{ZaSh}}, M. Ablovitz, D. Kaup, A.Newell and
H. Segur, have represented the Korteveg de Vries equation, the
nonlinear Schr\"{o}dinger equation, the Sine-Gordon equation and a
number of other non-linear equations of physical significance as
the compatibility conditions for the appropriate systems of linear
differential equations.
\end{remark}%

\section{FUCHSIAN EQUATIONS AND BEHAVIOR OF THEIR
SOLUTIONS -- GLOBAL AND LOCAL}
 \label{FEBTS}
 A Fuchsian equation is a linear differential equation of the
form
\begin{equation}%
\label{FDE}%
\frac{d Y}{d x}=\bigg(\sum\limits_{1\leq j\leq n}\frac{Q_j}{x-t_j}\bigg)Y ,
\end{equation}%
where \(t_1, \dots , t_n\) are pairwise distinct points of the
complex plane (\(t_1\not=t_2\not= \dots \not= t_n\)) ,
\(Q_j, j=1, \dots , n\) are square matrices of the same
dimension (say, \(k\times k\)) , \(x\) is the complex variable
which belongs to the punctured Riemann sphere
\(\overline{\mathbb{C}}\setminus\{t_1, \dots ,t_n\}\). The
points \(t_1, \dots ,t_n\) are regular singular points for the
equation (\ref{FDE}). The residue of the coefficient matrix
\(\sum\limits_{1\leq k\leq n}\frac{Q_j}{x-t_j}\) of the equation
(\ref{FDE}) at the point \(\infty\) is equal to \
\(-\!\!\!\sum\limits_{1\leq k\leq n}Q_k\). So, if the condition
\begin{equation}%
\label{IRP}%
\sum\limits_{1\leq k\leq n}Q_k=0
\end{equation}%
is satisfied, the point \(\infty\) is a regular point for the
equation (\ref{FDE}). In what follows we assume that the condition
(\ref{IRP}) is satisfied.

We consider the differential equation (\ref{FDE}) as a \textsl{matricial} one.
This means that we consider solutions \(Y(x)\) of this equation which are
\(k\times k\) matrix functions. If \(x_0\) is a regular point for the equation
(\ref{FDE}) (that is \(x_0\not=t_1, \dots , t_n\)), then for every
\(k\times k\) matrix \(Y_0\) there exists the solution \(Y(x; x_0, Y_0)\)
which satisfies the initial condition
\begin{equation}
\label{InCond}%
Y(x; x_0, Y_0)_{|x=x_0}=Y_0 .
\end{equation}%
This solution is a holomorphic function in a disk neighborhood of
the initial point \(x_0\). It can be continued analytically along
any path which does not pass through the singular points
\(t_1, \dots , t_n\). Thus the solution of the \textsl{linear}
differential equation (\ref{FDE}) exists \textsl{globally} and is
a (multi-valued) holomorphic function on the punctured Riemann
sphere \(\overline{\mathbb{C}}\setminus \{t_1, \dots , t_n\}\).
For a point \(x\in\overline{\mathbb{C}}\setminus
\{t_1, \dots , t_n\}\), the value of the analytic continuation
of the solution \(Y( . ; x_0, Y_0)\) along a path
\(\gamma(x_0, x)\subset\overline{\mathbb{C}}\setminus
\{t_1, \dots , t_n\}\), which starts at the point \(x_0\) and
ends at the point \(x\),
 can depend on this path \(\gamma\). If
\(\gamma^{\prime}(x_0, x)\) and \(\gamma^{\prime\prime}(x_0,x)\)
are two such paths,  homotopic in \(\overline{\mathbb{C}}\setminus
\{t_1, \dots , t_n\}\), then the values of analytic
continuations of the solution \(Y( . ; x_0, Y_0)\) along
\(\gamma^{\prime}\) and \(\gamma^{\prime\prime}\) coincide. Thus,
the value of the analytic continuation of the solution
\(Y( . ; x_0, Y_0)\) to the point \(x\) depends on \textsl{the
homotopy class of the path} leading from \(x_0\) to \(x\), along
which the continuation is made, but \textsl{not on the path
itself}. In other words, the solution \(Y( . , x_0, Y_0)\) of
the differential equation (\ref{FDE}) with the initial condition
(\ref{InCond}) is
 holomorphic function  on the universal covering space
\(\text{cov}(\overline{\mathbb{C}}\setminus \{t_1, \dots , t_n\}; x_0 )\)
of the space \(\overline{\mathbb{C}}\setminus \{t_1, \dots , t_n\}\) with
the distinguished point \(x_0\).

\begin{definition}
\label{Fund}%
A solution \(Y(x)\) of a linear differential equation in a connected domain is
said to be \textsf{fundamental} if \(\det Y(x)\not\equiv 0\) (or, what is
equivalent, if \(\det Y(x_0)\not= 0\) for some point \(x_0\) from this
domain).
\end{definition}

\noindent If \(Y(x)\) is a  fundamental solution of a linear
differential equation, and \(Z(x)\) is another solution of this
equation then \(Z(x)=Y(x)C\), where \(C\) is a constant (with
respect to \(x\)) matrix. In what follows we distinguish the
solution \(Y(x)\) of the differential equation (\ref{FDE}), which
satisfies the normalizing (initial) condition
\begin{equation}%
\label{ICond}%
Y(x_0)=I
\end{equation}%
at a distinguished point \(x_0\) which is regular for the equation
(\ref{FDE}).

Let \(\gamma(x_0, x_0)\) be a \textsl{closed} path, which starts at the point
\(x_0\) and ends \textsl{at the same point} \(x_0\). We call such a path
\textsf{\textsl{a loop with  the distinguished point}} \(x_0\). Let
\(x_0^{\ast}(\gamma)\) be the point of the universal covering space
\(\text{cov}(\overline{\mathbb{C}}\setminus \{t_1, \dots , t_n\}; x_0 )\),
which corresponds to the homotopic class \(\gamma\) of the loop
\(\gamma(x_0, x_0)\). We denote by \(M_{\gamma}\) the value
\(Y(x_0^{\ast}(\gamma); x_0, I)\) of the solution of the equation
(\ref{FDE}), normalized by the initial condition (\ref{ICond}), at the point
\(x_0^{\ast}(\gamma)\). The correspondence \(\gamma\mapsto M_{\gamma}\) is a
representation %
\footnote{ As usual, the product
\(\gamma^{\prime}\gamma^{\prime\prime}\) is  defined by the
composed loop, obtained by traversing first the loop, generating
\(\gamma^{\prime\prime}\), and then the loop, generating
\(\gamma{\prime}\).
}%
of the fundamental group \(\pi(\overline{\mathbb{C}}\setminus
\{t_1, \dots , t_n\}; x_0 )\) into the general linear group
\(GL(k; \mathbb{C})\):
\begin{equation*}%
\label{MR}%
M_{\gamma^{\prime}\gamma^{\prime\prime}}
=M_{\gamma^{\prime}}M_{\gamma^{\prime\prime}}\text{ \ for every \ }
\gamma^{\prime}, \gamma^{\prime\prime}\in\pi(\overline{\mathbb{C}}\setminus
\{t_1, \dots , t_n\}; x_0 )
\end{equation*}%
and
\begin{equation*}%
\label{MRU}%
M_{\varepsilon}=I,\text{ \ where \ }\varepsilon \text{ \ is the unity of the
group \ }\pi(\overline{\mathbb{C}}\setminus \{t_1, \dots , t_n\}; x_0 )
\end{equation*}%

\begin{definition}
\label{DMM} The matrix \(M_{\gamma}\) is said to be \textsf{ the
monodromy matrix} for the equation \textup{(\ref{FDE})},
corresponding to the homotopic class \(\gamma\). The
correspondence \(\gamma\mapsto M_{\gamma}\) is said to be
\textsf{the monodromy representation} of the equation
\textup{(\ref{FDE})}.
\end{definition}

\noindent The group \(\pi(\overline{\mathbb{C}}\setminus
\{t_1, \dots , t_n\}; x_0 )\) is finitely generated. Let us
choose a special
\begin{wrapfigure}[15]{r}{0.5\linewidth}
\begin{minipage}[H]{1.0\linewidth}
\epsfig{file=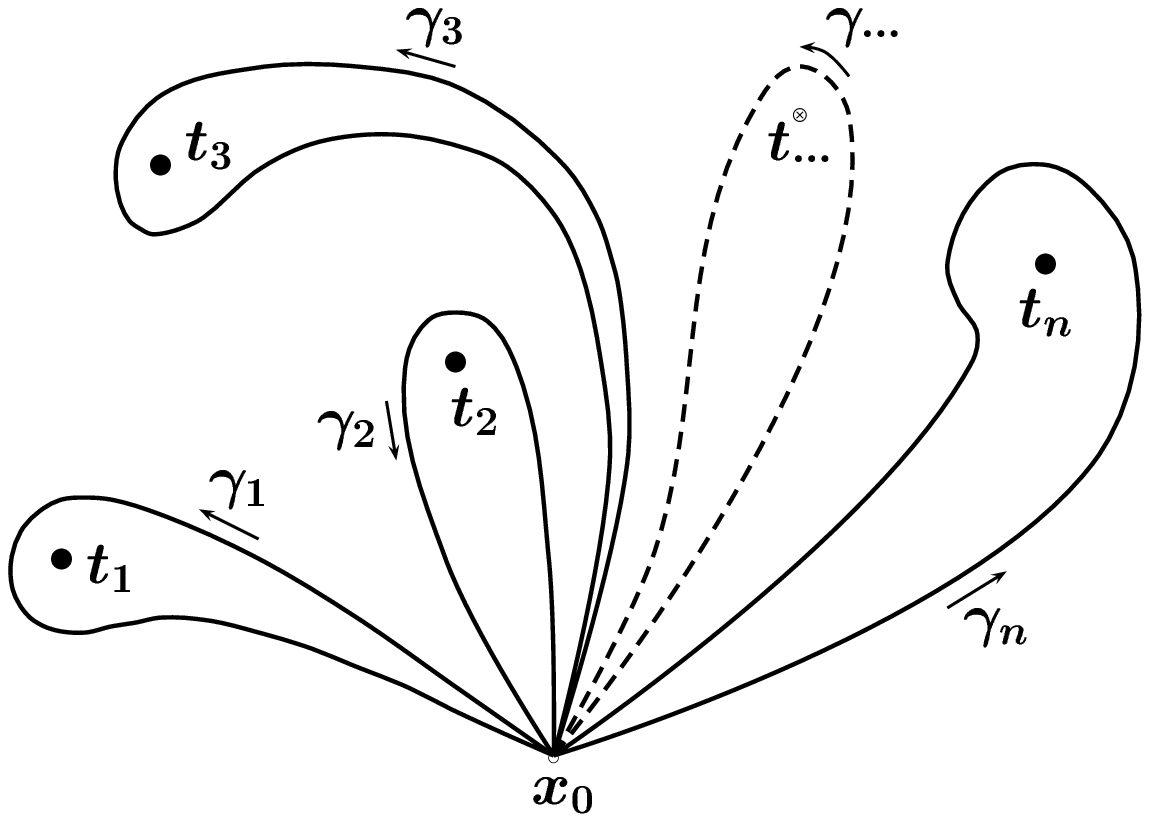,width=1.0\linewidth,clip=} \caption{}
\label{BGH}%
\end{minipage}
\end{wrapfigure}
system of generators of this group. Every simple loop (i.e. the
loop which is homeomorphic to a circle) divides the Riemann sphere
\(\overline{\mathbb{C}}\) into two components. Let \(\gamma_j\) be
a simple loop in the space \(\overline{\mathbb{C}}\setminus
\{t_1, \dots , t_n\}\), with the distinguished point \(x_0\),
 such that the point \(t_j\) and the set of all other points
\(\{t_p\}_{p\not=j}\) belong to the different components and,
moreover, the loop \(\gamma_j\)   goes around the point \(t_j\) in
the counterclockwise direction. (See Figure \ref{BGH}). We denote
the \mbox{homotopic} class of this loop by \(\gamma_j\) as well.
It is known that the elements \(\gamma_1, \dots , \gamma_n\)
generate the fundamental group
\(\pi(\overline{\mathbb{C}}\setminus \{t_1, \dots , t_n\}; x_0
)\) and that the equality
\(\gamma_1\cdot \cdots \cdot\gamma_n=\varepsilon\) is the
only generating relation
\footnote{ %
The generating relation in this form holds for the loops, arranged
as shown in \mbox{Figure \ref{BGH}}. In the general case, this
relation is of the form
\(\gamma_{p_1}\cdot \dots \cdot \gamma_{p_n}=\varepsilon\)
where \(p_1, \dots , p_n\) is a permutation of the indices
\(1, \dots , n\) which is determined by the arrangement of the
loops.
}.

\noindent For what follows, it is crucial to know how  a solution
of a Fuchsian differential  equation behaves in a neighborhood of
its singular point.

More generally, let us discuss the behavior of a solution of the
differential equation of the form

\begin{equation}%
\label{NSP}%
\frac{dY}{dx}=\bigg(\frac{Q}{x-t}+\Phi(x)\bigg)Y(x) ,
\end{equation}%
where \(Q\) is a constant \(k\times k\) matrix, and \(\Phi\) is a \(k\times
k\) matrix function, which is holomorphic in the disk \(D_{t, \rho}\),
centered at the point \(x=t\), of radius \(\rho\). Solutions of the equation
(\ref{NSP}) are holomorphic (in general, multi-valued) functions in the
punctured disk \(D_{t,\rho}\setminus \{t\}\).

The equation (\ref{NSP}) can be considered as a perturbation of the equation
\begin{equation}%
\label{UE}%
\frac{dY}{dx}=\frac{Q}{x-t} Y(x) .
\end{equation}%
The "unperturbed" equation (\ref{UE}) can be solved explicitly. It has a
fundamental solution  of the form
\begin{equation}%
\label{SUE}%
Y(x)=(x-t)^Q .
\end{equation}%
(By definition, \((x-t)^Q\stackrel{\textup{\tiny
def}}{=}e^{Q\ln(x-t)})\). In general, the function \((x-t)^Q\)
considered as a function in the complex plane is multivalued. It
can be considered as a univalued holomorphic matrix function on
the universal covering \(\cov (\mathbb{C}\setminus \{t\})\) of the
complex plane punctured at the point \(\{t\}\) (or, in other
terminilogy, on the Riemann surface of \(\log (x-t)\)). The
differentiation rule
\begin{equation}%
\label{DifRu}%
\frac{d }{d x} (x-t)^{Q}=\frac{Q}{x-t} (x-t)^{Q} ,\quad x\in \cov
\big(\mathbb{C}\setminus\{t\}\big) ,
\end{equation}%
holds. (The values of \((x-t)^Q\) in both sides of (\ref{DifRu}) are the
same).

The structure of the solution of the original "perturbed" equation (\ref{NSP})
depends on the nature of the spectrum of the matrix \(Q\).

\begin{definition}%
\label{DNRM}%
A square matrix \(Q\) is said to be \textsf{non-resonant} if no
two eigenvalues of \(Q\) differ by a nonzero integer, or in other
words, if the spectra of the matrices \(Q+nI\) and \(Q\) are
disjoint for every \(n\in\mathbb{Z}\setminus 0\).
\end{definition}%

\begin{remark}%
\label{MEV}%
It may happen that some eigenvalues of a non-resonant matrix \(Q\)
coincide (i.e., the eigenvalues of a non-resonant matrix \(Q\) may
be multiple).
\end{remark}

\begin{proposition} %
\label{SNRQ}%
If \(Q\) is a non-resonant (constant with respect to \(x\)) matrix, and
\(\Phi(x)\) is a matrix-function holomorphic in the disk \(D_{t,\rho}\), then
the differential equation \textup{(\ref{NSP})} has a fundamental solution
\(Y(x)\) of the form
\begin{equation}%
\label{FNRQ}%
Y(x)=H(x)(x-t)^{Q} ,\quad x\in\cov (D_{t,\rho}\setminus\{t\}) ,
\end{equation}
defined on the universal covering \(\cov
(D_{t,\rho}\setminus\{t\})\) of the punctured disc
\(D_{t,\rho}\setminus\{t\}\),
where \(H(x)\) is a matrix function, holomorphic and invertible %
\footnote{The inverse matrix function \(H^{-1}(x)\) exists and is holomorphic
in \(D_{t,\rho}\)}%
\ in the disk \(D_{t,\rho}\).
\end{proposition}

\noindent The proof of Proposition \ref{SNRQ}, and even of a more
general statement, is given in Appendix.

In the general case, when the  eigenvalues of the matrix \(Q\) can
differ by non-zero integers, the situation is more complicated.

\begin{proposition}%
\label{SRQ}%
The differential equation  \textup{(\ref{NSP})}, with a constant
matrix \(Q\) and a matrix function \(\Phi(x)\), holomorphic in the
disk \(D_{t, \rho}\), has a fundamental solution \(Y(x)\) of the
form
\begin{equation}%
\label{FRQ}%
Y(x)=H(x)(x-t)^L(x-t)^{\hat{Q}} , \quad x\in\cov
(D_{t,\rho}\setminus\{t\}) ,
\end{equation}
where \(H(x)\) is a matrix function holomorphic and invertible in
the disk \(D_{t,\rho}\), \(L\) is a diagonalizable matrix with
integer eigenvalues \(l_1,\dots,l_k\), and \(\hat{Q}\) is a
non-resonant matrix whose eigenvalues
\(\hat{\lambda}_1,\dots,\hat{\lambda_k}\) are related to the
eigenvalues \(\lambda_1,\dots,\lambda_k\) of the matrix \(Q\)
according to
\begin{equation}%
\label{CN}%
 \hat{\lambda}_p=\lambda_p-l_p , \qquad 1\leq p\leq k .
\end{equation}%
For a certain choice of the integers \(l_1\geq \dots \geq l_n\), the
matrices \(L\) and \(Q\) can be simultaneously reduced to the forms
\begin{equation*}%
L=T^{-1} \textup{Diag}\big(l_1, \dots , l_n\big) T,\ %
\hat{Q}=T^{-1} \Big(\textup{Diag}\big(\hat{\lambda}_1, \dots ,
 \hat{\lambda}_n\big)+U\Big) T ,
\end{equation*}%
where \(U\) is an upper-triangular matrix with zero diagonal.
\end{proposition}%

\noindent Proposition \textup{\ref{SRQ}} is applicable to equation
(\ref{NSP}) with the general matrix \(Q\). In particular, it is
applicable to the equation with the non-resonant matrix \(Q\).
However, if \(Q\) is non-resonant, one can choose \(L=0\) and
\(\hat{Q}=Q\).

\begin{remark}%
\label{NoUn}%
It should be mentioned that the matrix function \(H(x)\) in
(\ref{LocS}) (or in (\ref{LocSR})) is not quite unique. The matrix
\(H(t)=H(x)_{x=t}\) can be an arbitrary invertible matrix
commuting with the matrix \(Q\). Of course, the solution
(\ref{FNRQ}) (or (\ref{FRQ})) depends on this choice. In the
\textsf{\textsl{non-resonant case}}, as soon as this choice is
made, the local solution of the equation \eqref{NSP} at the point
\(x=t\), i.e. the solution of the form \eqref{FNRQ}, is determined
uniquely (up to the choice of \(\arg (x-t)\) which leads to the
non-uniqueness of \((x-t)^Q=e^{Q\ln (x-t)}\)). In what follows, we
choose
\begin{equation*}%
H(x)_{|x=t}=I.
\end{equation*}%
 \textsl{In the non-resonant case},
the latter equality can be considered as a \textit{normalizing
condition} for the local solution (\ref{FNRQ}) of \eqref{NSP} at
the point \(x=t\).

However, in the resonant case, the condition \(H(x)_{|x=t}=I\) do
\textsf{\textsl{not}} determine the solution of the form
\eqref{FRQ} of the equation \eqref{NSP} uniquely. The following
example illustrates this non-uniqueness. Consider the differential
equation
\begin{equation*}%
\frac{\partial Y(x)}{\partial x}=\frac{Q}{x}\cdot Y(x),
\end{equation*}%
where
\begin{equation*}%
Q=
\begin{bmatrix}
1&0\\
0&0
\end{bmatrix}\,\cdot
\end{equation*}%
This equation, which is of the equation \eqref{NSP} with \(t=0,\Phi(x)\equiv 0\),
has the solution %
\begin{equation*}%
Y_1(x)=%
\begin{bmatrix}%
x&0\\
0&1
\end{bmatrix},%
\end{equation*}%
which is of the form \eqref{FRQ} with \(H(x)=\begin{bmatrix}
1&0\\
0&1
\end{bmatrix}\),
\(L=\begin{bmatrix}
1&0\\
0&0
\end{bmatrix}\),
\(\hat{Q}=\begin{bmatrix}
0&0\\
0&0
\end{bmatrix},\)\\
\(t=0\).
The matrix function
\begin{equation*}
Y_2(x)=Y_1(x)\cdot%
\begin{bmatrix}%
1&1\\
0&1
\end{bmatrix},%
\quad\textup{or, more explicitely,}\quad Y_2(x)=
\begin{bmatrix}%
x&x\\
0&1
\end{bmatrix},
\end{equation*}
is another solution of the same differential equation. However,
the function \(Y_2(x)\) is also representable in the form
\begin{equation*}
Y_2(x)=
\begin{bmatrix}%
1&x\\
0&1
\end{bmatrix}\cdot
\begin{bmatrix}%
x&0\\
0&1
\end{bmatrix},%
\end{equation*}
i.e. in the form \eqref{FRQ}, with \(H(x)=\begin{bmatrix}%
1&x\\
0&1
\end{bmatrix}\),
\(L=\begin{bmatrix}
1&0\\
0&0
\end{bmatrix}\),
\(\hat{Q}=\begin{bmatrix}
0&0\\
0&0
\end{bmatrix},t=0.\)
 Thus, the considered differential equation possesses two
different solutions \(Y_1(x)\) and \(Y_2(x)\), both these
solutions are of the form \eqref{FRQ}, with the same \(L\) and
\(\hat{Q}\), but with different \(H(x)\), and both these factors
\(H(x)\) satisfy the condition \(H(x)_{|x=0}=I\).

Thus, \textit{in the resonant case, the condition
\(H(x)_{|x=t}=I\) does not uniquely determine the local solution
\eqref{FRQ} of the equation \eqref{NSP} corresponding to the
singular point \(x=t\), even sough the matrices \(L\) and
\(\hat{Q}\) are specified.}
\end{remark}%

\begin{remark}%
\label{ED}%
There is another essential difference between the non-resonant and
the "resonant" cases. If the matrix \(Q\) is not non-resonant,
then the matrix \(\hat{Q}\), which appears in
\textup{(\ref{FRQ})}, depends on the value
\(\Phi(t)=\Phi(x)_{x=t}\) of the ``perturbing'' (see (\ref{NSP}))
function \(\Phi\) at the point \(t\). In particular, the Jordan
structure of the matrix \(\hat{Q}\) may depend on the matrix
\(\Phi(t)\). This Jordan structure, together with the eigenvalues
of \(\hat{Q}\), determine the multi-valued behavior of the matrix
function \((x-t)^{\hat{Q}}\), and finally, the multi-valued
behavior of \(Y(x)\) in \(D_{t, \rho}\). Thus, in the "resonant"
case the multi-valued behavior of the solution \(Y(x)\) depends
not only on the matrix \(Q\) (the coefficient of the ``leading
term'')  but also on the ``perturbing'' matrix \(\Phi\). As  can
be seen from \textup{(\ref{FNRQ})}, in the non-resonant case the
multi-valued behavior of the solution \(Y(x)\) depends only on
\(Q\), but not on \(\Phi\).
\end{remark}%

\begin{remark}%
Of course, for every matrices \(A\) and \(B\), the matrix \(e^A\cdot e^B\) can
be presented as \(e^C\):
\[e^A\cdot e^B=e^C,\qquad \text{where} \qquad C=\log (e^A\cdot e^B) ,\]
where an arbitrary value for the matrix \(\log (e^A\cdot e^B)\)
can be chosen. However, the case of exponential \textup{functions}
is different. In general, a matrix \textup{function}
\(e^{A\zeta}\cdot e^{B\zeta}\) of \(\zeta\) (with constant
matrices \(A\) and \(B\)) is not representable as \(e^{C\zeta}\)
with a constant matrix \(C\). In particular, in the resonant case
the solution \textup{(\ref{FRQ})} of the differential equation
\textup{(\ref{NSP})} is not, in general, representable in the form
\(Y(x)=H(x)\cdot (x-t)^{R}\) with a holomorphic invertible matrix
function \(H\) and a constant matrix \(R\).
\end{remark}

\noindent Propositions \ref{SNRQ} and \ref{SRQ} go back to
L. Fuchs and G. Frobenius. The systematic study of linear
differential equations with analytic coefficients in the
neighborhood of a singular point was initiated by Lazarus Fuchs in
1866 and followed by G. Frobenius in 1873. They discussed the
"scalar" differential equation of the \(n^{\text{th}}\) order.
Matrix differential equations in the complex domain were
considered later. The statement which we formulated as Proposition
\ref{SNRQ} can be found in a lot of excellent textbooks now. See,
for example, \cite{CoLe}, Theorem 4.1; \cite{Gant}, Chapter XV, \S
10; \cite{Hart}, Chapter IV, Corollary 11.2; \cite{IKSY}, Chapter
I, Theorem 2.3. However, it is difficult to point out the source,
where this statement first appeared precisely in the same form as
in Proposition \ref{SNRQ}. The earliest reference  we know is
\cite{Ra} (see Section 3 there). The statement which we formulated
as Proposition \ref{SRQ} can be found in
[Gant], Chapter XV, \S 10, in
\cite{HsSi}, Chapter V, Section 5;
\cite{Sib}, Chapter III, Theorem 3.9.4,page 89. The earliest
of this references is \cite{Gant} (the first edition of
Gantmacher's book appeared in 1954). However, Y. Sibuya,
\cite{Sib}, page 89, refers to Hukuhara's book \cite{Huk} (Theorem
41.1, page 143) published in 1950. Unfortunately, the book
\cite{Huk} is unavailable to us.

Let us return to the Fuchsian equation (\ref{FDE}). We consider
this equation in a disk neighborhood \(D_{t_j, \rho}\) of the
singular point \(t_j\), which does not contain other singular
points \(t_p, p\not=j\). Let us choose and fix a point \(x_j\) of
the punctured disk \(D_{t_j, \rho}\setminus\{t_j\}\). We also
choose and fix a path \(\alpha_j\) in
\(\overline{\mathbb{C}}\setminus \{t_1, \dots , t_n\}\), which
starts at the distinguished point \(x_0\) and ends at the

\begin{figure}[h]%
\begin{minipage}[t]{0.31\linewidth}%
\centering\epsfig{figure=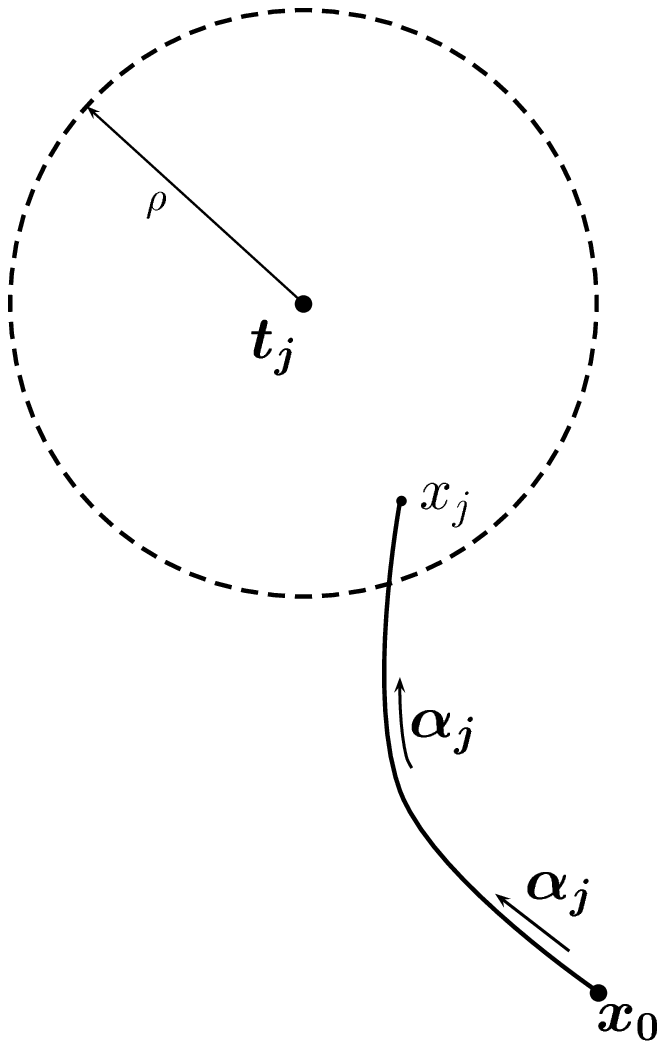,width=1.0\linewidth}%
\caption{a: \mbox{path \(\alpha_j\)}.}
\end{minipage}%
\addtocounter{figure}{-1}%
\hfill
\begin{minipage}[t]{0.31\linewidth}%
 \centering\epsfig{figure=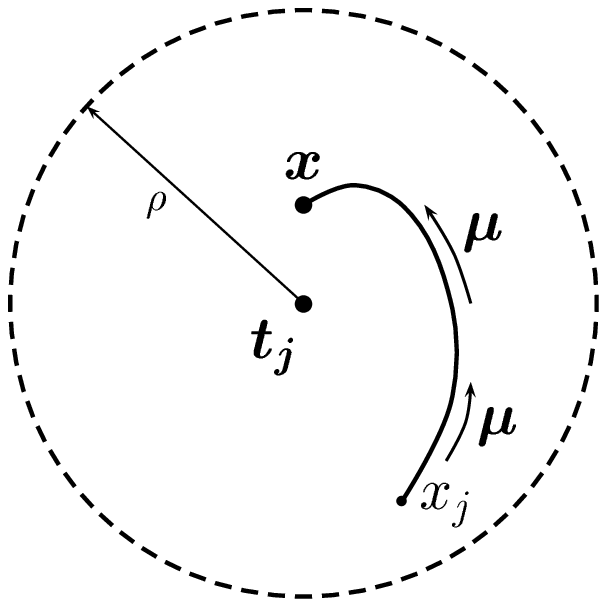,width=1.0\linewidth}%
\caption{b: \mbox{path \(\mu\)}.}
\end{minipage}%
\addtocounter{figure}{-1}%
\hfill
\begin{minipage}[t]{0.31\linewidth}%
\centering\epsfig{figure=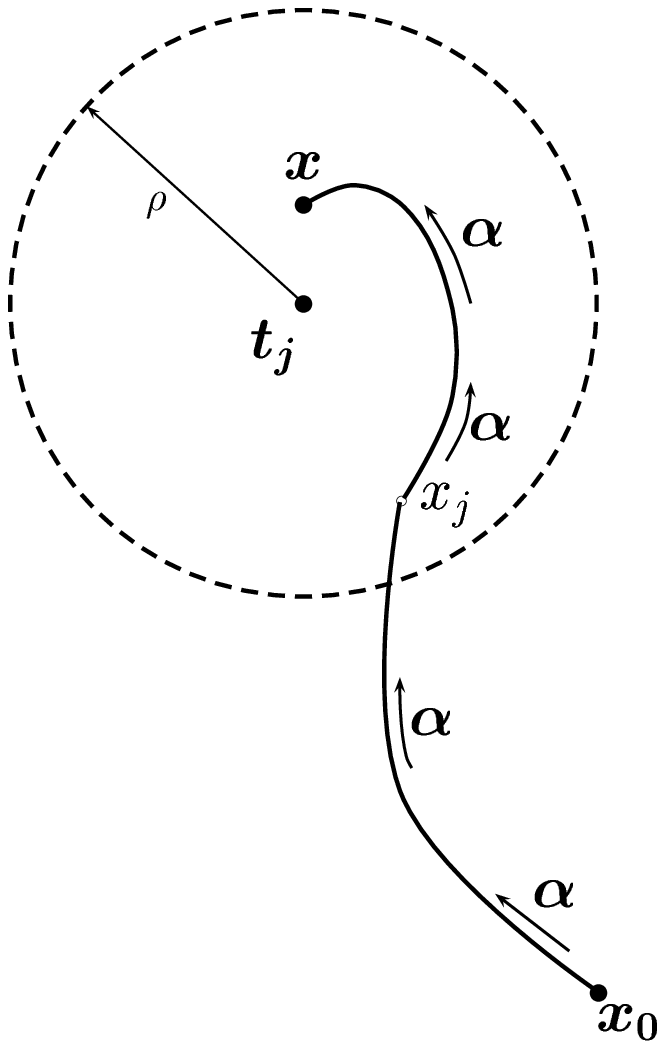,width=1.0\linewidth}%
\caption{c: \mbox{path \(\alpha\)}.}\label{FundDj}%
\end{minipage}%
\end{figure}%

\noindent%
 point \(x_j\).
(See Figure \ref{FundDj}.a.)
 We need these
\(x_j\) and \(\alpha_j\) to distinguish a certain analytic
continuation \(Y_{\alpha_j}\) of the solution of (\ref{FDE})
normalized by (\ref{ICond}) into the punctured disk \(D_{t_j,
\rho}\setminus\{t_j\}\). Let \(x\) be an arbitrary point of the
punctured disk \(D_{t_j, \rho}\setminus\{t_j\}\), and \(\mu\) be
an arbitrary path in \(D_{t_j, \rho}\setminus\{t_j\}\), starting
at \(x_j\) and ending at \(x\). (See Figure \ref{FundDj}.b.). So,
the pair \((x,\mu)\) can be considered as a point of the universal
covering \(\cov(D_{t_j, \rho}\setminus\{t_j\};x_j)\). We construct
the composed path \(\alpha = \mu\cdot \alpha_j\). The path
\(\alpha\) starts at the point \(x_0\) and ends at the point
\(x\). (See Figure \ref{FundDj}.c.) We  continue analytically the
solution \(Y\) along the path \(\alpha\) from the point \(x_0\) to
the point \(x\). In this manner we distinguish the continuation
\(Y_{\alpha_j}(x)\) of the solution \(Y( . )\), normalized by the
condition (\ref{ICond}), from a neighborhood of \(x_0\) to the
punctured disk \(D_{t_j, \rho}\setminus\{t_j\}\). The function
\(Y_{\alpha_j}( . )\) can be considered as a multi-valued function
in \(D_{t_j, \rho}\setminus\{t_j\}\), or, better to say, as a
single-valued function  in the universal covering \(\cov(D_{t_j,
\rho}\setminus\{t_j\};x_j)\) with the distinguished point \(x_j\).
The notation \(Y_{\alpha_j}\) reflects how the path \(\alpha_j\)
is involved in the process of the analytic continuation. In fact,
this process involves two steps. On the first step, the solution
is continued along the path \(\alpha_j\) from a neighborhood of
the initial point \(x_0\) to a neighborhood of the point \(x_j\),
where the point \(x_j\) is considered as a distinguished point of
the punctured disk \(D_{t_j,\rho}\setminus\{t_j\}\). On the second
step, the solution is continued along paths which are contained in
\(D_{t_j, \rho}\setminus\{t_j\}\) from the neighborhood of the
point \(x_j\) to points of \(\cov (D_{t_j,
\rho}\setminus\{t_j\};x_j)\).

The matrix function \(Y_{\alpha_j}( . )\) is a fundamental
solution of (\ref{FDE}). Propositions \ref{SNRQ}-\ref{SRQ} provide
another fundamental solution of (\ref{FDE}). As a preliminary we
assume that the matrix \(Q_j\) is non-resonant.  According to
Proposition \ref{SNRQ}, there exists the fundamental solution
\(Y_j(x)\) of (\ref{FDE}), which is representable as
\begin{equation}%
\label{LocS}%
Y_j(x)=H_j(x)(x-t_j)^{Q_j} ,\qquad x\in
\cov\big(D_{t_j,\rho}\setminus\{t_j\},x_j\big) ,
\end{equation}%
where \(H_j(x)\) is a matrix function, holomorphic and invertible
in the whole (non-punctured) disc \(D_{t_j, \rho}\). Under the
normalizing condition
\begin{equation}%
\label{NCSiP}%
H_j(x)_{|x=t_j}=I,
\end{equation}%
the matrix function \(H_j(x)\) is unique.
 The matrix function \((x-t_j)^{Q_j}\) is defined up to a
right factor of the form \(\exp[2\pi iQ_j]\). To avoid this non-uniqueness, we
choose a certain value
\begin{equation}
\label{ChArg}%
\theta_j=\arg (x_j-t_j) ,
\end{equation}
and fix this choice. Under the agreement (\ref{ChArg}), the matrix
function \((x-t)^{Q_j}\) is defined uniquely on
\(\cov\big(D_{t_j,\rho}\setminus\{t_j\},x_j\big)\). Thus,
\textit{if the matrix \(Q_j\) is non-resonant, the local solution
\(Y_j(x)\) of the differential equation \eqref{FDE} corresponding
to the singular point \(t_j\) can be distinguished uniquely by
means of the normalizing conditions \eqref{NCSiP} and
\eqref{ChArg}.}

If the matrix \(Q_j\) is resonant, the situation is different.
According to Proposition \ref{SRQ}, there exists the fundamental
solution \(Y_j(x)\) of (\ref{FDE}), which is representable as
\begin{equation}%
\label{LocSR}%
Y_j(x)=H_j(x)(x-t_j)^{L_j}(x-t_j)^{\hat{Q}_j} ,\qquad x\in
\cov\big(D_{t_j,\rho}\setminus\{t_j\}\big) ,
\end{equation}%
where the matrix \(H_j(x)\) has the same properties as before, \(L_j\) is a
diagonalizable matrix with integer eigenvalues (so, the matrix function
\((x-t_j)^{L_j}\) is single-valued in \(D_{t,\rho}\setminus\{t_j\}\)), and
\(\hat{Q}_j\) is a non-resonant matrix, whose eigenvalues differ from eigenvalues of
the matrix \(Q_j\) by some integer numbers. However, even though the normalizing
conditions \eqref{NCSiP} and \eqref{ChArg} are assumed and the matrices
\(L_j,\hat{Q}_j\) are chosen, the matrix function \(H_j(x)\) in (\ref{LocSR})) is
determined {\em non-uniquely}\footnote{See Remark \ref{NoUn}. In what follows, such
an arbitrariness in the choice of the local solution does not play any essential
role. We need the notions of the local solution and of the connection coefficient
only to \textsl{illustrate} and \textsl{motivate} the notion of the
\textsf{isoprincipal deformation} which will be introduced further, in section
\ref{IsoDef}.}. Thus, \textit{if the matrix \(Q_j\) is resonant, we may deal with
some local solution \(Y_j(x)\) of the form \eqref{LocSR} corresponding to the
singular point \(t_j\) of the differential equation \eqref{FDE}, but in general, we
can not distinguish naturally such a solution.}

\begin{definition}
\label{DefLoSo}%
Assume that the matrix \(Q_j\) is non-resonant. The solution \(Y_j\) of the Fuchsian
equation \textup{(\ref{FDE})} representable in the form \textup{(\ref{LocS})}   near
\(\{t_j\}\), with \(H_j\) satisfying the normalizing condition (\ref{NCSiP}), is
said to be \textsf{the normalized local solution of (\ref{FDE}) corresponding to the
singular point \(\{t_j\}\)}.
\end{definition}

Let us choose  some local solution \(Y_j\) of the equation (\ref{FDE}) corresponding
to the singular point \(t_j\). If the matrix \(Q_j\) is non-resonant, it is natural
to choose the normalized local solution. If the matrix  \(Q_j\) is not non-resonant,
we choose some local solution of the form \eqref{FRQ} and fix such a choice.
 Since both solutions,
\(Y_{x_j,\alpha_j}\) and \(Y_j\), are fundamental, the equality
\begin{equation}%
Y_{x_j,\alpha_j}(x)=Y_j(x)C,
\label{ConRel}%
\end{equation}%
holds,  where \(C\) is a constant (with respect to \(x\))
non-degenerate matrix.

\begin{definition}%
\label{Connec}%
The equality \textup(\ref{ConRel}) is said to be \textsf{the connection relation.}
The matrix \(C\) from \textup(\ref{ConRel}) is said to be \textsf{the connection
coefficient}.
\end{definition}%

 According to Definition \ref{Connec}, the connection
coefficient \(C\) depends on the choice of the point \(x_j\in
D_{t_j,\rho}\setminus\{t_j\}\) and on the choice of the path \(\alpha_j\),
connecting the points \(x_0\) and \(x_j\): \(C=C_{x_j,\alpha_j}\). Now we will
discuss this dependence.

Let \(\alpha_j^{\prime}\) and  \(\alpha_j^{\prime\prime}\) be two
paths, both of which start at the point \(x_0\) and end at the
point \(x_j\), \(\alpha_j^{\prime},
\alpha_j^{\prime\prime}\subset\overline{\mathbb{C}}\setminus
\{t_1, \dots , t_n\}\). It is clear that
\(\alpha_{j}^{\prime\prime}=\alpha_{j}^{\prime}\cdot\gamma\),
where \(\gamma\) is the loop with the distinguished point \(x_0\),
which is composed of the path \(\alpha_{j}^{\prime\prime}\) and of
the path \((\alpha_{j}^{\prime})^{-1}\) (the latter is the path
\(\alpha_{j}^{\prime}\) which is traversed in the opposite
direction: from the point \(x_j\) to the point \(x_0\)):
\(\gamma=(\alpha_{j}^{\prime})^{-1}\cdot\alpha_{j}^{\prime\prime}\).
The paths \(\alpha^{\prime}\) and \(\alpha^{\prime\prime}\),
connecting the point \(x_0\) with the point \(x\), are defined as
the compositions \(\alpha^{\prime}=\mu\cdot\alpha_{j}^{\prime}\)
and \(\alpha^{\prime\prime}=\mu\cdot\alpha_{j}^{\prime\prime}\),
where, as before, \(\mu\) ,  \(\mu\subset
D_{t_j,\rho}\setminus\{t_j\}\),  is a path connecting the points
\(x_j\) and \(x\). Thus,
\(\alpha^{\prime\prime}=\alpha^{\prime}\cdot\gamma\), and
\begin{equation*}
Y_{\alpha_{j}^{\prime\prime}}(x)=%
Y_{\alpha_{j}^{\prime}}(x)M_{\gamma},
\end{equation*}
where \(M_{\gamma}\) is the monodromy matrix corresponding to the loop
\(\gamma\). On the other hand,
\begin{equation*}%
Y_{\alpha_j^{\prime}}(x)=Y_j(x) C_{x_j,\alpha_j^{\prime}},\qquad
Y_{\alpha_j^{\prime\prime}}(x)=Y_j(x) C_{x_j,\alpha_j^{\prime\prime}} .
\end{equation*}%
Therefore,
\begin{equation}%
\label{CConC}%
C_{x_j,\alpha_j^{\prime\prime}}=C_{x_j,\alpha_j^{\prime}} M_{\gamma} , \quad
\text{where }
\gamma=(\alpha_{j}^{\prime})^{-1}\cdot\alpha_{j}^{\prime\prime} .
\end{equation}%

Let us discuss, how the connection coefficient \(C\) depends  on the point
\(x_j\). Given two points \(x_j^{\prime}, x_j^{\prime\prime}\in
D_{t_j, \rho}\setminus\{t_j\}\) and a path
\(\alpha_j^{\prime}\subset\overline{\mathbb{C}}\setminus
\{t_1, \dots , t_n\}\), connecting \(x_0\) and \(x_j^{\prime}\), there
exists a path \(\alpha_j^{\prime\prime}\subset\overline{\mathbb{C}}\setminus
\{t_1, \dots , t_n\}\), connecting \(x_0\) and \(x_j^{\prime\prime}\), such
that
\begin{equation}%
\label{ICP}%
C_{x_j^{\prime},\alpha_j^{\prime}}=C_{x_j^{\prime\prime},\alpha_j^{\prime\prime}},
\end{equation}%
Indeed, let \(\nu\) be a path, connecting \(x_j^{\prime}\) and
\(x_j^{\prime\prime}\), \(\nu\subset D_{t_j, \rho}\setminus\{t_j\}\). We define the
path \(\alpha_j^{\prime\prime}\) connecting \(x_0\) with \(x_j^{\prime\prime}\) as
the composition \(\alpha_j^{\prime\prime}=\nu\cdot\alpha_j^{\prime}\). The pairs
\((x_j^{\prime},\alpha_j^{\prime})\) and
\((x_j^{\prime\prime},\alpha_j^{\prime\prime})\) allow us to distinguish the
solutions \(Y_{\alpha_j^{\prime}}\) and \(Y_{\alpha_j^{\prime\prime}}\) of the
equation (\ref{FDE}) in the punctured neighborhood \(D_{t_j,\rho}\setminus\{t_j\}\).
To obtain values of these solutions at a point \mbox{\(x\in
D_{t_j,\rho}\setminus\{t_j\}\)}, we have to choose paths \(\mu_j^{\prime}\) and
\(\mu_j^{\prime\prime} \), connecting the points \(x_j^{\prime}\) and
\(x_j^{\prime\prime}\), respectively, with the point \(x\). If we make this choice
coherently, we obtain the equality \(Y_{\alpha_j^{\prime}}(x)=
Y_{\alpha_j^{\prime\prime}}(x)\) which implies (\ref{ICP}). Such a coherent choice
can be done in the following way.
 If \(\mu^{\prime}\) is a
path in \(D_{t_j,\rho}\setminus\{t_{j}\}\), connecting \(x_j^{\prime}\) with \(x\),
we define the path \(\mu^{\prime\prime}\), connecting \(x_j^{\prime\prime}\) with
\(x\), as the composition \(\mu^{\prime\prime}=\mu^{\prime}\cdot\nu^{-1}\). To
obtain the values of the solutions \(Y_{x_j^{\prime},\alpha_j^{\prime}}\) and
\(Y_{x_j^{\prime\prime},\alpha_j^{\prime\prime}}\) at the point \(x\), we have to
continue analytically the solution \(Y\), which is normalized by (\ref{ICond}), from
the point \(x_0\) to the point \(x\) along the paths
\(\alpha^{\prime}=\mu^{\prime}\cdot\alpha_j^{\prime}\) and
\(\alpha^{\prime\prime}=\mu^{\prime\prime}\cdot\alpha_j^{\prime\prime}\),
respectively. However, these two paths are homotopic. To see this, we remark that
\(\alpha^{\prime\prime} = (\mu^{\prime}\cdot\nu^{-1})\cdot(\nu\cdot\alpha^{\prime})
= \mu^{\prime}\cdot(\nu^{-1}\cdot\nu)\cdot\alpha^{\prime}\) and that the path
\(\nu^{-1}\cdot\nu\) is contractible to the point \(x_j^{\prime}\). Hence,
\(Y_{\alpha_j^{\prime}}(x)= Y_{\alpha_j^{\prime\prime}}(x)\).
 Thus, in some sense the connection coefficient  does not depend on the
choice of the point \(x_j\in D_{t_j, \rho}\setminus\{t_j\}\).

\begin{wrapfigure}[18]{r}{0.35\linewidth}
\begin{minipage}[H]{1.0\linewidth}
 \epsfig{file=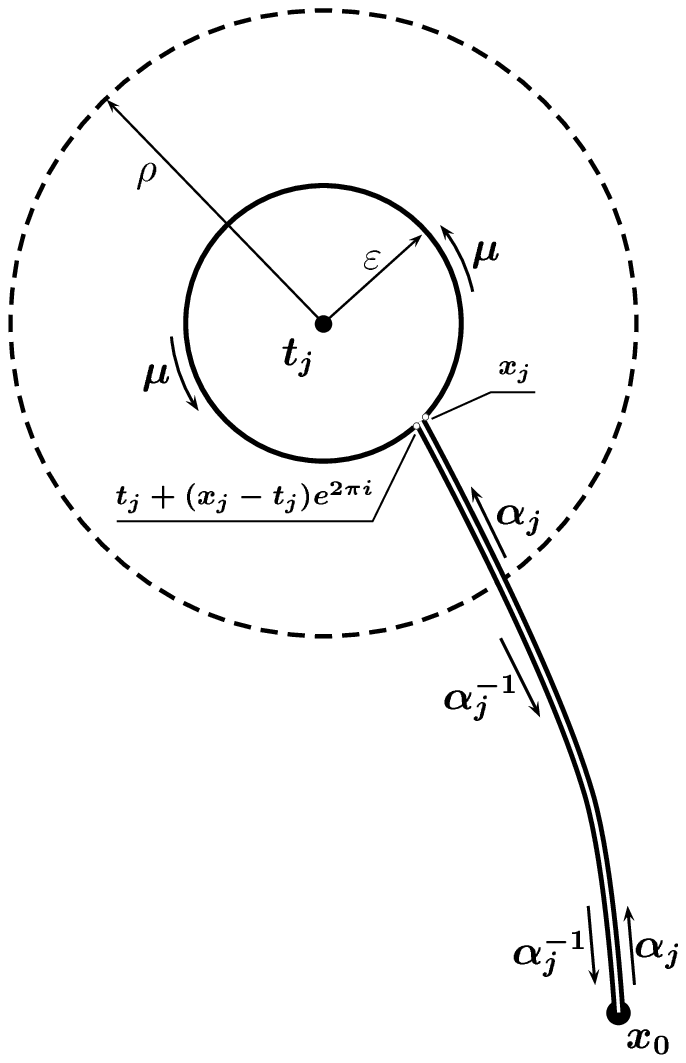,width=\linewidth,clip=}
\caption{\!\!\hbox{Loop \(\gamma_{\alpha_j}.\)}} %
\label{BGHj}%
\end{minipage}
\end{wrapfigure}
 Let us calculate the monodromy matrix corresponding to the loop
``going around'' the point \(t_j\). We relate such a loop \(\gamma_{\alpha_j}\) to
every path \(\alpha_j\), \(\alpha_j\subset\overline{\mathbb{C}}\setminus\{t_1, \dots
, t_n\}\), which starts at the point \(x_0\) and ends at a point \(x_j\) belonging
to a punctured disk \(D_{t_j,\rho}\setminus\{t_j\}\). The radius \(\rho\)  must be
chosen so small that the disk does not contain any other point \(t_p, p\not=j\). Let
\(\mu\) be the circular arc, centered at the point \(t_j\), of the radius
\(\varepsilon=|x_j-t_j|\), which starts at the point \(x_j\), ends at the point
\(t_j+(x_j-t_j)e^{2\pi i}\) and is oriented counterclockwise, making precisely one
turn around its center \(t_j\).
 We construct the loop \(\gamma_{\alpha_j}\) as the
composition of the path \(\alpha_j\), of the arc \(\mu\) and of the path
\(\alpha_j^{-1}\): \(\gamma=(\alpha_j)^{-1}\cdot\mu\cdot\alpha_j\) (see Figure
\ref{BGHj}).

\begin{definition}%
\label{LoGen}%
The loop \(\gamma_{\alpha_j}\) is said to be \textsf{the loop generated by the
path \(\alpha_j\)}.
\end{definition}%

\begin{lemma}%
\label{MALGBP}%
The monodromy matrix \(M_{\gamma_{\alpha_j}}\) along the loop
\(\gamma_{\alpha_j}\) generated by a path \(\alpha_j\) can be represented as
\begin{equation}%
\label{EFM}%
M_{\gamma_{\alpha_j}}=(Y_{\alpha_j}(x_j))^{-1} Y_{\alpha_j}(t_j+(x_j-t_j)e^{2\pi
i}) ,
\end{equation}%
where \(x_j\) is the endpoint of the path \(\alpha_j\).
\end{lemma}%

\noindent%
\textsf{PROOF.} %
For linear differential equations, the following \textsl{semigroup
property} holds: Let \(\gamma\) be a path starting at a point
\(a\) and ending at the point \(c\), and let \(b\) is an
``intermediate''  point of the path \(\alpha\), that is the point
which is located ``between'' the points \(a\) and \(c\). Then
\(Y(c, a)=Y(c,b) Y(b,a)\), where for \(\xi,\eta\in\alpha\),
\(Y(\eta,\xi))\) is the value \(Y(x,\xi)_{|x=\eta}\) of the
solution of the equation satisfying the initial condition
\(Y(x,\xi)_{|x=\xi}=I\). Let us take the loop
\(\gamma_{\alpha_j}\) as such a path \(\gamma\). The point \(a\)
is the point \(x_0\) considered as the starting point of the loop,
the point \(c\) is the same point \(x_0\), but considered as the
end point of the loop. The point \(b\) is the point
\(t_j+(x_j-t_j)e^{2\pi i}\). From the semigroup property it
follows that \(Y(b,c)Y(c,b)=I\). It is clear that
\(Y(b,a)=Y_{\alpha_j}(t_j+(x_j-t_j)e^{2\pi i})\),
\(Y(b,c)=Y_{\alpha_j}(x_j)\), and
\(Y(c,a)=M_{\gamma_{\alpha_j}}\). \hfill Q.E.D.\\[1ex]

\noindent From (\ref{ConRel}) it follows that in the non-resonant case,%
\begin{equation*}%
Y_{x_j,\alpha_j}(x_j)=H_j(x_j)(x_j-t_j)^{Q_j}C_{\alpha_j}
\end{equation*}%
and
\begin{equation*}%
Y_{x_j,\alpha_j}(t_j+(x_j-t_j)e^{2\pi i})=H_j(x_j)e^{2\pi
iQ_j}(x_j-t_j)^{Q_j}C_{\alpha_j} .%
\end{equation*}%
(The function \(H_j(x)\) is single-valued in \(D_{t_j, \rho}\),
and the function \((x-t_j)^{Q_j}\) receives the factor \(e^{2\pi
iQ_j}\), when \(x\) traverses the arc \(\mu\) once, in the
counterclockwise direction.) Therefore, according to Lemma
\ref{MALGBP},
\begin{equation}%
\label{expMon}%
 M_{\gamma_j}=
 C_{\alpha_j}^{-1}(x_j-t_j)^{-Q_j}e^{2\pi iQ_j}(x_j-t_j)^{Q_j}C_{\alpha_j}.
\end{equation}%
Finally, we obtain

\begin{lemma}%
\label{EMM}%
The monodromy matrix \(M_{\gamma_{\alpha_j}},\) corresponding to the loop
\(\gamma_{\alpha_j}\) generated by a path \(\alpha_j\) can be represented as
 \begin{equation}%
 \label{MME}%
 M_{\gamma_j}=C_{\alpha_j}^{-1}\exp[2\pi iQ_j]C_{\alpha_j}
 \end{equation}%
 in the non-resonant case, and as
\begin{equation}%
 \label{MMER}%
 M_{\gamma_j}=C^{-1}_{\alpha_j}\exp[2\pi i\hat{Q}_j]C_{\alpha_j}
 \end{equation}%
 in the resonant case, where \(C_{\alpha_j}\) is the connection
 coefficient,
 corresponding to the path \(\alpha_j\).
 \end{lemma}%

 \noindent(The reasoning in the resonant case is analogous. One needs to use
 the equality (\ref{LocSR}), taking into account that the matrix-function
 \((x-t_j)^{L_j}\) is single-valued in the disk \(D_{t_j,\rho}\)).

For what follows, it is very fruitful to rewrite the representations
(\ref{LocS}) - (\ref{ConRel}), (\ref{LocSR}) - (\ref{ConRel}) for a
``local'' solution %
\footnote{ I.e., for a solution in a neighborhood of the singular point
\(t_j\).%
}%
of (\ref{FDE}), as well as the expressions (\ref{MME}) and
(\ref{MMER}) for the monodromy matrices, in a different form. The
respective transformations are trivial, but nevertheless very
useful. The form in which we present the local solutions hints at
how to introduce the notion of the isoprincipal deformation of the
Fuchsian differential equation. This notion will be of crucial
importance to us.

Let  \(D_{t_j,\rho}\) be a disk neighborhood of the singular point \(t_j\) of
the differential equation (\ref{FDE}), which is small enough, so that the
point \(t_j\) is the only singular point which is contained in this
neighborhood. In this neighborhood lives a ``local'' solution \(Y_j\) of the
equation (\ref{FDE}), which in the non-resonant case is of the form
(\ref{LocS}), and in the resonant case is of the form (\ref{LocSR}). However,
this ``local'' solution \(Y_j\)  does not, in general, satisfy the initial
condition (\ref{ICond}). To be precise, we should  speak about the analytic
continuation of this solution from the neighborhood \(D_{t_j,\rho}\) to the
point \(x_0\). To discuss such an analytic continuation, we have to choose a
point \(x_j\) of the punctured neighborhood \(D_{t_j,\rho}\setminus\{t_j\}\),
and a path \(\alpha_j\),
 \(\alpha_j\subset\overline{\mathbb{C}}\setminus\{t_1, \dots , t_n\}\),
 leading from the point \(x_0\) to the point \(x_j\). The solution
 \(Y_{\alpha_j}\) which satisfies the initial condition (\ref{ICond}) at the
 distinguished point \(x_0\) and which is continued analytically to the
 neighborhood \(D_{t_j,\rho}\) along the path \(\alpha_j\) is of the
 form
 \begin{equation}
 \label{LSRefN}
 Y_{\alpha_j}(x)=H_j(x)(x-t_j)^{Q_j} C_{\alpha_j} ,\qquad x\in
\cov\big(D_{t_j,\rho}\setminus\{t_j\},x_j\big) ,
\end{equation}
in the non-resonant case, and of the form
\begin{equation}
 \label{LSRefR}
 Y_{\alpha_j}(x)=H_j(x)(x-t_j)^{L_j}(x-t_j)^{\hat{Q_j}} C_{\alpha_j} ,\qquad x\in
\cov\big(D_{t_j,\rho}\setminus\{t_j\},x_j\big) ,
\end{equation}
in the resonant case, where \(C_{\alpha_j}\) is the connection coefficient
corresponding to the choice of the homotopic class of the path \(\alpha_j\)
and the choice of the arguments \(\arg (x_j-t_j)\), see (\ref{ChArg}),
(\ref{LocS}) and (\ref{LocSR}).

We introduce the matrices
\begin{equation}%
\label{RewrQ}%
A_{\alpha_j}=C_{\alpha_j}^{-1}Q_jC_{\alpha_j} ,\quad
H_{\alpha_j}=H_j C_{\alpha_j}
\end{equation}%
in the non-resonant case, and the matrices
\begin{equation}%
\label{RewrQRes}%
A_{\alpha_j}=C_{\alpha_j}^{-1}\hat{Q}_jC_{\alpha_j} ,\quad
B_{\alpha_j}=C_{\alpha_j}^{-1}L_jC_{\alpha_j} ,\quad
H_{\alpha_j}=H_j C_{\alpha_j}
\end{equation}%
in the resonant case. Propositions \ref{SNRQ} and \ref{SRQ} can be
``developed'' in the following way.

\begin{theorem}
\label{SNRQI}%
If the matrix \(Q_j\) is non-resonant, then the solution \(Y_{\alpha_j}\) of
the differential equation \textup{(\ref{FDE})} that satisfies the initial
condition \textup{(\ref{ICond})} at the distinguished point \(x_0\) and is
continued analytically
 along a path \(\alpha_j\) into a small punctured neighborhood
 \(D_{t_j,\rho}\setminus\{t_j\}\)  of the singular point
 \(t_j\), is representable in the form
\begin{equation}%
\label{RSNRQI}%
Y_{\alpha_j}(x)=H_{\alpha_j}(x)(x-t_j)^{A_{\alpha_j}} , \qquad
x\in \cov\big(D_{t_j,\rho}\setminus\{t_j\},x_j\big) ,
\end{equation}
 where \(A_{\alpha_j}\) is a matrix,  similar to
 the matrix \(Q_j\), and \(H_{\alpha_j}\) is a matrix function,
 holomorphic and invertible in the entire neighborhood
\(D_{t_j,\rho}\) (including the point \(t_j\)).
\end{theorem}

\begin{theorem}
\label{SRQI} If the matrix \(Q_j\) is resonant, then the solution
\(Y_{\alpha_j}\) of the differential equation
\textup{(\ref{FDE})}, that satisfies the initial condition
\textup{(\ref{ICond})} at the distinguished point \(x_0\) and is
continued analytically
 along a path \(\alpha_j\) into a small punctured neighborhood
 \(D_{t_j,\rho}\setminus\{t_j\}\)  of the singular point
 \(t_j\), is representable in the form
\begin{equation}%
\label{RSNRQ}%
Y_{\alpha_j}(x)=H_{\alpha_j}(x)(x-t_j)^{Z_{\alpha_j}}(x-t_j)^{A_{\alpha_j}}
, \qquad x\in \cov\big(D_{t_j,\rho}\setminus\{t_j\},x_j\big) ,
\end{equation}
 where \(Z_{\alpha_j}\) is a diagonalizable matrix with integer eigenvalues
 \(l_1, \dots , l_k\),
 \(A_{\alpha_j}\) is a non-resonant matrix, whose
  eigenvalues \(\hat{\lambda}_p\) are related
 to the eigenvalues
\(\lambda_p\) of the
matrix \(Q_j\) by the relation %
\begin{equation}%
\label{EVD}%
\hat{\lambda}_p=\lambda_p-l_j,\quad 1\leq p\leq k ,
\end{equation}%
and \(H_{\alpha_j}(x)\) is a matrix function,  holomorphic and
invertible in the entire neighborhood (including the point
\(t_j\)). For a certain choice of the integers
\(l_1\geq \dots \dots \geq l_n\), the matrices \(Z_{\alpha_j}\)
and \(A_{\alpha_j}\) can be simultaneously reduced to the forms
\begin{equation*}%
Z_{\alpha_j}=T_{\alpha_j}^{-1}
 \textup{Diag}\big(l_1, \dots , l_n\big) T_{\alpha_j},\ %
A_{\alpha_j}=
T_{\alpha_j}^{-1} \Big(\textup{Diag}\big(\hat{\lambda}_1, \dots ,
 \hat{\lambda}_n\big)+U\Big) T_{\alpha_j} ,
\end{equation*}%
where \(U\) is an upper-triangular matrix with zero diagonal.
\end{theorem}

\begin{remark}%
\label{ChP}%
If \(\alpha_{j}^{\prime}\) is another path with the same starting point
\(x_0\) and endpoint \(x_j\),
 then, in the non-resonant case,
 the matrices \(A_{\alpha_j^{}}\) and \(H_{\alpha_j^{}}\) from
(\ref{RSNRQI}) are
 transformed according to the rule
\begin{equation}%
\label{TrNR}%
A_{\alpha_j^{\prime}}=M_{\gamma}^{-1}A_{\alpha_j^{}}M_{\gamma},\qquad
H_{\alpha_j^{\prime}}(x)=H_{\alpha_j}(x)M_{\gamma} ,
\end{equation}%
and in the resonant case, the matrices \(Z_{\alpha_j}\), \(A_{\alpha_j}\) and
\(H_{\alpha_j}\) from (\ref{RSNRQ}) are
 transformed according to the rule
\begin{equation}%
\label{TrR}%
Z_{\alpha_j^{\prime}}=M_{\gamma}^{-1}Z_{\alpha_j^{}}M_{\gamma},\quad
A_{\alpha_j^{\prime}}=M_{\gamma}^{-1}A_{\alpha_j^{}}M_{\gamma},\quad
H_{\alpha_j^{\prime}}(x)=H_{\alpha_j}(x)M_{\gamma} ,
\end{equation}%
 where \(M_{\gamma}\) is the monodromy matrix
corresponding to the loop \(\gamma=(\alpha_j)^{-1}\alpha_j^{\prime}\). In
particular, if the paths
\(\alpha_j\) and \(\alpha_j^{\prime}\) are homotopic then
\[Z_{\alpha_j^{\prime}}=Z_{\alpha_j^{}}, \quad
A_{\alpha_j^{\prime}}=A_{\alpha_j^{}}, \quad
H_{\alpha_j^{\prime}}(x)=H_{\alpha_j}(x) .\]
\end{remark}%

\begin{theorem}%
\label{EMMT}%
Let \(M_{\gamma_{\alpha_j}}\) be the monodromy matrix for the differential
equation (\ref{FDE}) corresponding to the loop \(\gamma_{\alpha_j}\), which
goes around the singular point \(t_j\) and is constructed from the path
\(\alpha_j\), connecting the initial point \(x_0\) to a small neighborhood of
\(t_j\): \(\gamma_j=(\alpha_j)^{-1}\cdot\mu\cdot\alpha_j\) (see \textup{Figure
\ref{BGHj}}). Then this monodromy matrix can be expressed as
\begin{equation}%
\label{FEMM}%
M_{\gamma_{\alpha_j}}=\exp (2\pi i A_{\alpha_j}) ,
\end{equation}%
where the matrix \(A_{\alpha_j}\) is the same as in the exponent
in (\ref{RSNRQI}), if \(Q_j\) is non-resonant, and the exponent in
(\ref{RSNRQ}), if \(Q_j\) is resonant.
\end{theorem}%

\begin{theorem}%
\label{SpMaTe} %
Let \(M_{\gamma_j}\) be the monodromy matrix for the differential
equation (\ref{FDE}) corresponding to a loop \(\gamma_{j}\) with
the distinguished point \(x_0\) which makes one turn
counterclockwise around the singular point \(t_j\) and makes no
turns around the other singular points \(t_p, p=1, \dots , t_n,
p\not=j\). Then, whatever the matrix \(Q_j\) is, non-resonant or
not, the spectra \(\spec M_{\gamma_j}\) and \(\spec Q_j\) are
related by
\begin{equation}%
\label{SpMT}%
\spec M_{\gamma_j}=\exp{\{2\pi i \spec Q_j }\} ,
\end{equation}%
where the latter equality holds ``with multiplicities'', that is the
multiplicity \(\kappa_{\mu}\) of the eigenvalue \(\mu\) of the matrix
\(M_{\gamma_j}\) and multiplicities \(\kappa_{\lambda}\) of the eigenvalues
\(\lambda\) of the matrix \(Q_j\) are related by
\(\kappa_{\mu}=\sum\kappa_{\lambda}\), were the sum is taken over all
\(\lambda\) such that \(\exp\{\lambda\}=\mu\).
\end{theorem}%

\noindent%
\textsf{PROOF}. %
This theorem is the immediate consequence of Theorem \ref{EMMT},
of the spectral mapping theorem for matrices (from
\textup(\ref{FEMM}) it follows that \(\spec
M_{\gamma_j}=\exp{\{2\pi i \spec A_{\alpha_j}}\}\)), and of
Theorems \ref{SNRQI} or \ref{SRQI}, which relate the spectra of
the matrices \(A_{\alpha_j}\) and \(Q_j\). However, here we can do
without Theorem \ref{SRQI}, restricting ourselves to non-resonant
matrices \(Q_j\) only. Every matrix can be approximated by a
non-resonant one. The spectrum \(\spec Q_j\) depends on \(Q_j\)
continuously. The monodromy matrix \(M_{\gamma_j}\) depends on
\(Q_j\) continuously as well. \hfill Q.E.D.

\section{DEFORMATIONS OF FUCHSIAN DIFFERENTIAL EQUATIONS \label{DefFuchsEq}}
 A \textsl{deformation} of a Fuchsian equation is \textsl{a family} of Fuchsian
 differential equations depending on parameters:
\begin{equation*}%
\label{FDEP}%
\frac{d Y}{d x}=\bigg(\sum\limits_{1\leq j\leq
n}\frac{Q_p ({\boldsymbol\alpha})}{x-t_p}\bigg)Y .
\end{equation*}%
 If the coefficients
\(Q_p({\boldsymbol\alpha})\) depend on the parameters
\(\boldsymbol{\alpha}=(\alpha_1, \dots , \alpha_m)\) in
some open set of the space
\(\mathbb{C}^m\) holomorphically, such a deformation is said
to be a \textsl{holomorphic deformation of a Fuchsian differential equation}.
The study of linear differential equations that depend on parameters
holomorphically was started by Lazarus Fuchs, see
\cite{FuL1} - \cite{FuL4}. This study was continued by L. Schlesinger, see
\cite{Sch1} - \cite{Sch3}, and by Richard Fuchs, see \cite{FuR}. In
particular, L. Schlesinger has considered linear equation, where loci
\(t_1, \dots , t_n\) of singular points serve as the parameters
\(\boldsymbol{\alpha}\). In other words, L. Schlesinger has considered
deformations of Fuchsian differential equations of the form

\begin{equation}%
\label{FDES}%
\frac{d Y}{d x}=\bigg(\sum\limits_{1\leq p\leq
n}\frac{Q_p ({\boldsymbol{t}})}{x-t_p}\bigg)Y ,
\end{equation}%
where the the coefficients \(Q_p(t_1, \dots , t_n), p=1, \dots ,
n,\) are \(k\times k\) matrix functions which are defined and
holomorphic for \(\boldsymbol{t}\) from a domain \(\mathcal{D}\)
of the space \(\mathbb{C}^n_{\ast}\). In the further we assume
that the condition
\begin{equation*}%
\sum\limits_{1\leq p\leq n}Q_p ({\boldsymbol{t}})\equiv 0, \quad
\boldsymbol{t}\in\mathcal{D},
\end{equation*}%
is satisfied.

In what follows we study the normalized by (\ref{ICond}) solution
\(Y(x,\boldsymbol{t}))\) of the family (\ref{FDES}) of Fuchsian equations as a
function of \(x\) and \(\boldsymbol{t}\). In particular, we consider the
monodromy matrices \(M_{\gamma_j}\) for the equation (\ref{FDES}), as well as
the matrices \(A_{\alpha_j}\) and functions \(H_{\alpha_j}(x)\) from the
representations (\ref{RSNRQI}) and (\ref{RSNRQ}) for various
\(t_1, \dots , t_n\): \(M_{\gamma_j}=M_{\gamma_j}(\boldsymbol{t})\),
\(A_{\alpha_j}=A_{\alpha_j}(\boldsymbol{t})\) and
\(H_{\alpha_j}=H_{\alpha_j}(\boldsymbol{t})\), intending to study how these
matrices depend on \(\boldsymbol{t}=(t_1, \dots , t_n)\). Here a problem
arises.  The matrices \(M_{\gamma_j}\), \(A_{\alpha_j}\), \(H_{\alpha_j}\)
depend on the homotopy classes of \(\gamma_j\) and \(\alpha_j\),
which are considered on the punctured Riemann sphere
\mbox{\(\overline{\mathbb{C}}\setminus\{t_1, \dots , t_n\}\)}. However, for
different \(\boldsymbol{t}=(t_1, \dots , t_n)\), the sets
\mbox{\(\overline{\mathbb{C}}\setminus\{t_1, \dots , t_n\}\)} are
different.
 Thus we  need  to explain what does it mean that,
for the deformation (\ref{FDES}) of the Fuchsian equation, the monodromy
matrices \(M_{\gamma_j}\), or the exponents \(A_{\alpha_j}\) from
(\ref{RSNRQI}) (or (\ref{RSNRQ})) are holomorphic with respect to
\(\boldsymbol{t}\), or what does it mean that these matrices do not depend on
\(\boldsymbol{t}\). The latter problems are of the local nature. %

Let a point
\(\boldsymbol{t}^0=(t_1^0, \dots , t_n^0)\in\mathbb{C}^n_{\ast}\)
be fixed, \(t_j^0\not=x_0,\ j=1, \dots , n\), where
\begin{wrapfigure}[18]{r}{0.60\linewidth}
\begin{minipage}[H]{1.0\linewidth}
\epsfig{file=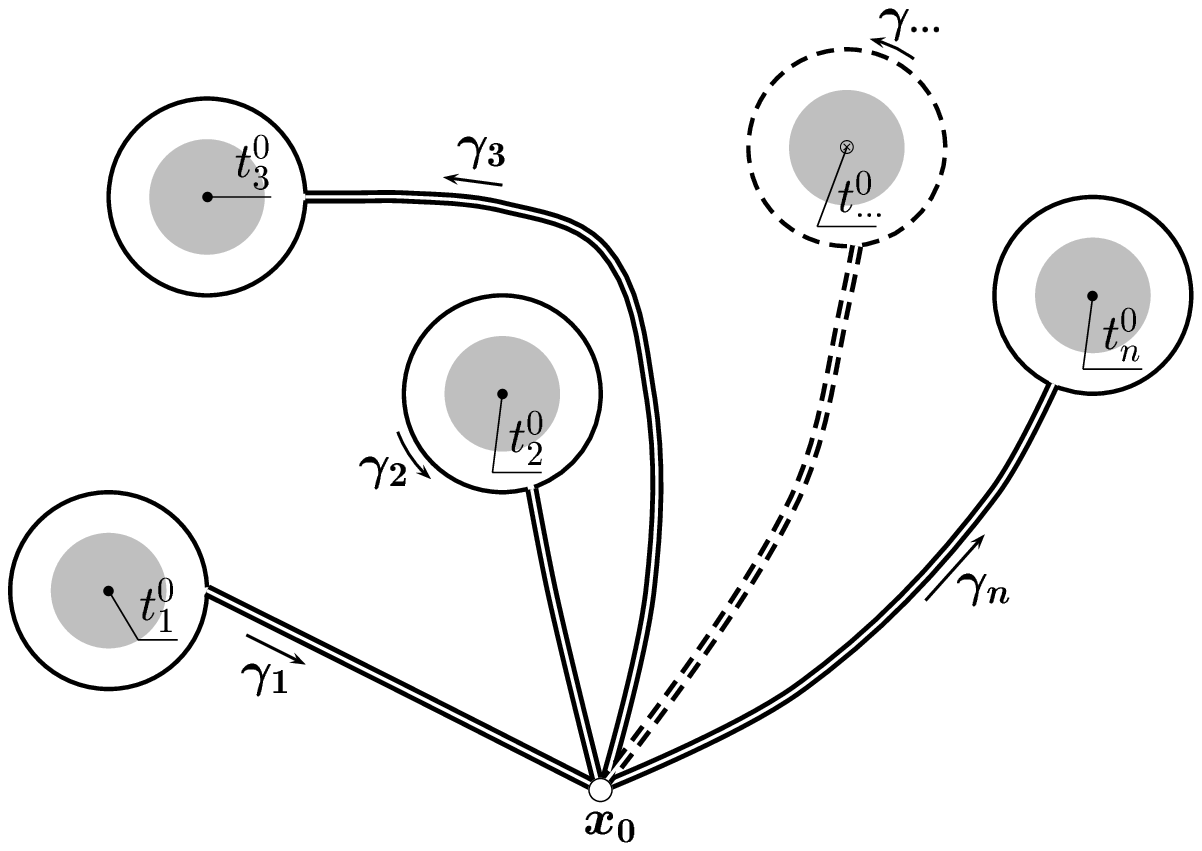,width=\linewidth,clip} \caption{}
\label{BGHV}%
\end{minipage}
\end{wrapfigure}
\(x_0\in\overline{\mathbb{C}}\) is the distinguished point (the
same point which appears in the initial condition(\ref{ICond}),
and which serves as the distinguished point for the fundamental
group \(\pi\big(\overline{\mathbb{C}}\setminus
\{t_1^0,, \dots , t_n^0\}; x_0 \big)\)). Choose and
\textsl{fix} simple loops \(\gamma_j, j=1, \dots ,\gamma_n\),
each of which starts and ends at the point \(x_0\), the loop
\(\gamma_j\)  contains
   the
point \(t_j^0\) ''inside'' and goes around this point
counterclockwise, and all other points \(\{t_p^0\}_{p\not=j}\)
 lie ''outside" the loop \(\gamma_j\), \(j=1, \dots , t_n\).
Let \(\delta>0\) be a small number such that the disk
\(D_{t_j^0, \delta}= \{z\in\mathbb{C}:|z-t_j^0|<\delta\}\) does
not intersect  the loop \(\gamma_j, j=1, \dots , n\) (and
hence, is contained inside the loop \(\gamma_j\)). Then for every
\(\boldsymbol{t}=(t_1, \dots , t_n)\) such that \(t_j\in
D_{t_j^0, \delta}\), the loops
\(\gamma_j, j=1, \dots , \gamma_n,\) generate the fundamental
group \(\pi\big(\overline{\mathbb{C}}\setminus
\{t_1 , \dots , t_n\}; x_0 \big)\)). See Figure \ref{BGHV},
where the disks \(D_{t_j^0, \delta}\) are plotted in the light
gray color.

Assume that the point \(\boldsymbol{t}^0\) belongs to the domain
\(\boldsymbol{\mathcal{D}}\) of the space \(\mathbb{C}^n_{\ast}\),
where the matrix functions \(Q_j(\boldsymbol{t})\) are
holomorphic, and that \(\delta\) is so small that the polydisk
\begin{equation}%
\label{DPD}%
\mathsf{D}(\boldsymbol{t}^0, \delta)=\prod\limits_{p=1}^nD_{t_p^0,\delta},%
\quad\text{where}\quad D_{t_p^0,\delta}=\{t_p: |t_p-t_p^0|<\delta\} ,
\end{equation}%
is contained in \(\boldsymbol{\mathcal{D}}\):
\(\mathsf{D}(\boldsymbol{t}^0, \delta)\subset\boldsymbol{\mathcal{D}}\).
 Then, for \(x\in\gamma_j\), the coefficient matrix
 \(\sum\limits_{1\leq k\leq
n}\frac{Q_j ({\boldsymbol{t}})}{x-t_j}\) of the linear equation (\ref{FDES}) is
holomorphic with respect to \(\boldsymbol{t}\) in
\(\mathsf{D}(\boldsymbol{t}^0, \delta) .\) From the standard results on dependence
of solution of a linear differential equation on parameters  follows

\begin{lemma}%
\label{HolMon} Under above stated assumptions and notation, the
monodromy matrices \(M_{\gamma_j}(\boldsymbol{t}),\ j=1,\dots,n,\)
of the holomorphic family \textup{(\ref{FDES})} of Fuchsian
equations are holomorphic with respect to \(\boldsymbol{t}\) from
the polydisk \(\mathsf{D}(\boldsymbol{t}^0,\delta).\)
\end{lemma}%

\noindent Now we investigate in more detail  the local solution
\(Y_j(x,\boldsymbol{t})\) of the equation (\ref{FDES}) which
corresponds to a neighborhood of the point \(t_j\). We have
already considered such a local solution (see Propositions
\ref{SNRQ} - \ref{SRQ} and Theorems \ref{SNRQI} - \ref{SRQI}).
However, until now we considered \(\boldsymbol{t}\) as
\textsl{fixed} and focused our attention on the multiplicative
representation (\ref{LocS} - (\ref{LocSR}) of the solution
\textsl{considered as a function of} \(x\). Now we emphasize the
dependence of the local solution on \(\boldsymbol{t}\), while
concentrating on the \textsl{non-resonant case}. Given
\(\boldsymbol{t}\in\mathbb{C}^n_{\ast}\) (such that the
coefficients \(Q_p(\boldsymbol{t}), p=1, \dots , n,\) are
defined for this \(\boldsymbol{t}\)), and assuming that the matrix
\(Q_j(\boldsymbol{t})\) is non-resonant for certain \(j\), the
local solution \(Y_j(x,\boldsymbol{t})\) of the equation
(\ref{FDES}) corresponding to the singular point \(t_j\) is sought
in the form
\begin{equation}%
\label{DT}%
Y_j(x,\boldsymbol{t})=H_{j}(x,\boldsymbol{t})(x-t_j)^{Q_j(\boldsymbol{t})} ,
\quad x\in D_{t_j, \rho}\setminus \{t_j\} .
\end{equation}%
The differentiation rule (\ref{DifRu}) for \((x-t_j)^{Q}\) leads to the differential
equation for \(H_j(x,\boldsymbol{t})\):
\begin{equation}%
\label{ForH}%
\frac{d H_j(x,\boldsymbol{t}))}{dx}=
\left(\frac{Q_j(\boldsymbol{t})H_j(x,\boldsymbol{t})-
H_j(x,\boldsymbol{t})Q_j(\boldsymbol{t})}{x-t_j}
+\Phi_j(x,\boldsymbol{t})\right)
 H_j(x,\boldsymbol{t}) ,
\end{equation}%
where
\begin{equation}%
\label{ForPhi}%
\Phi_j(x,\boldsymbol{t})=\sum\limits_{\substack{1\leq p\leq
n,\\
p\not=j}}\frac{Q_p ({\boldsymbol{t}})}{x-t_p} ,
\end{equation}%
the function \(\Phi_j(x,\boldsymbol{t})\) is holomorphic with respect to \(x\)
at \(x=t_j\).
 The equation (\ref{ForH}) does not depend on the
choice of branch of the (multivalued, in general) function
\((x-t_j)^{Q_j(\boldsymbol{t})}\). The right hand side of this equation,
considered as a function of \(x\), has the singularity at the point \(x=t_j\).
Nevertheless, it can be shown that, if the matrix \(Q_j(\boldsymbol{t})\) is
non-resonant, then, under the normalization condition
\begin{equation}%
\label{NLS}%
H_j(x,\boldsymbol{t})_{|x=t_j}=I ,
\end{equation}%
the holomorphic solution \(H_j(x,\boldsymbol{t})\) of the equation
(\ref{ForH}) exists for \(x\) from a disk centered at \(t_j\), and such a
solution is unique. More precise formulation of this result will be done
below.

 In fact, here we repeat what we have already formulated in
the Proposition \ref{SNRQ} and expressed in the formula (\ref{LocS}). The only
difference between the formulae (\ref{DT}) and (\ref{LocS}) is that the
notation in (\ref{DT}) reflects the dependence on \(\boldsymbol{t}\)
explicitly.

\begin{definition}
\label{DNLS}%
 Let \(\boldsymbol{t}=(t_1, \dots , t_n)\in\mathbb{C}^n_{\ast}\)
be given,
and assume that the matrix \(Q_j(\boldsymbol{t})\) is non-resonant.

 The \(k\times k\) matrix function \(Y_j(x,\boldsymbol{t})\), possessing
 the properties
\begin{enumerate}
\item[\textup{i.}] \(Y_j(x,\boldsymbol{t})\) satisfies the differential equation
\textup(\ref{FDES}) with respect to \(x\) in  a punctured neighborhood
\(D_{t_j, \rho}\setminus\{t_j\}\) of the point \(t_j\),\ \
\textup{(\(D_{t_j, \rho}=\{z\in\mathbb{C}: |z-t_j|<\rho\}\)} ,\ %
\(\rho\) is a positive number);
\item[\textup{ii.}]
\(Y_j(x,\boldsymbol{t})\) admits the factorization of the form
\textup{(\ref{DT})}, where the factor \(H_j(x,\boldsymbol{t})\) is a matrix
function holomorphic with respect to \(x\) in the whole (non-punctured) disk
\(D_{t_j, \rho}\) and satisfying the normalizing condition
\textup{(\ref{NLS})},
\end{enumerate}
 is said to be \textsf{the normalized local solution of the differential
 equation \textup{\textrm{(\ref{FDES})}}, corresponding to the singular
 point \(t_j\)}.
\end{definition}

\begin{proposition}%
\label{ELS}%
Let \(Q_p(\boldsymbol{t}), p=1, \dots , n,\) be \(k\times k\)
matrix functions, which are holomorphic with respect to
\(\boldsymbol{t}=(t_1, \dots , t_n)\) for \(\boldsymbol{t}\)
from an open set \(\boldsymbol{\mathcal{D}},
\boldsymbol{\mathcal{D}}\in\mathbb{C}^n_{\ast}.\) Let
\(\boldsymbol{t}^0=(t_1^0, \dots , t_n^0)\) be a point from
\(\boldsymbol{\mathcal{D}}\), and let
\begin{equation}%
\label{DefRho}%
\rho_q(\boldsymbol{t}^0)=\frac{1}{2}\min_{\substack{1\leq p\leq n\\[0.1ex]
p\not=q}}\big|t_{q}^{0}-t_p^0\big| , \quad 1\leq q\leq n .
\end{equation}%
Assume that, for certain \(j\), the matrix
\(Q_j(\boldsymbol{t}^0)\) is non-resonant. Then there exists
\(\delta, \delta>0,\) such that for every
\(\boldsymbol{t}=(t_1, \dots , t_n)\) from the polydisk %
\footnote{The polydisk \(\mathsf{D}(\boldsymbol{t}^0, \delta)\) was defined
in (\ref{DPD}) .} %
\(\mathsf{D}(\boldsymbol{t}^0, \delta)\), the normalized local solution
\(Y_j(x,\boldsymbol{t})\) of the equation (\ref{FDES}), which corresponds to
the singular point \(t_j\),
 exists for \(x\) from the punctured disk %
\footnote{\(D_{t_j, \rho_j(\boldsymbol{t}^0)}
=\{z\in\mathbb{C}:|z-t_j|<\rho_j(\boldsymbol{t}^0)\} .\)}
 \(D_{t_j, \rho_j(\boldsymbol{t}^0)}\setminus\{t_j\}\),
and the factor \(H_j(x,\boldsymbol{t})\) in (\ref{DT}) is a matrix function
holomorphic with respect to \(x,\boldsymbol{t}\) for \(x\in
D_{t_j, \rho_j(\boldsymbol{t}^0)}\),
\(\boldsymbol{t}\in\mathsf{D}(\boldsymbol{t}^0, \delta)\). Moreover, the
matrix \(H_j(x,\boldsymbol{t})\) is invertible for these \(x\) and
\(\boldsymbol{t}\).
\end{proposition}%

\noindent The proof of this proposition will be given in Appendix.

\begin{remark}%
\label{PTDV}%
 Proposition \textup{\ref{ELS}} can be considered as a version of
Proposition \textup{\ref{SNRQ}}, with the emphasis being placed on dependence
on \(\boldsymbol{t}\).
\end{remark}%

\begin{remark}%
\label{NCSP}%
Proposition \ref{ELS} belongs to the class of statements of the following nature: if
coefficients of a differential equations and a normalizing condition depend on some
parameters analytically, then the solution of this equation depends on these
parameters analytically as well. However, in Proposition \ref{ELS} the normalizing
condition (which can be presented in the form \(\lim_{x\to t_j} Y_j(x,
\boldsymbol{t})(x-t_j)^{-Q_j(\boldsymbol{t})}=I\)) is posed at the \textsf{singular}
point \(t_j\).  Therefore, Proposition \ref{ELS} is not a consequence of standard
general results about analytic dependence of a solution of a differential equation
on parameters. The proof  makes explicit use
 of the non-resonance of the matrix
\(Q_j(\boldsymbol{t}^0)\) and without this assumption   the result fails. This shows
that the result is rather delicate.
\end{remark}%

In section \ref{FEBTS} the notion of connection coefficients was introduced. The
connection coefficients relate the solution of the Fuchsian equation, which is
normalized at the distinguished point \(x_0,\) to its local solutions, which are
normalized at the singular points of the equation. However, in section \ref{FEBTS}
we did not  care about the dependence of connection coefficients on
\(\boldsymbol{t}\). Now we focus our attention on this dependence,
concentrating on the non-resonant case.

Given %
\footnote%
{\(\boldsymbol{\mathcal{D}}\) is the open set in
\(\mathbb{C}^n_{\ast}\) where
the coefficients \(Q_p,1\leq p\leq n,\) are defined and holomorphic.} %
 \(\boldsymbol{t}^0,
\boldsymbol{t}^0\in\boldsymbol{\mathcal{D}}\), and given \(j, 1\leq j\leq n\), let
us choose \(\rho_j(\boldsymbol{t}^0)\) according to (\ref{DefRho}). Then choose
\(\delta>0\), satisfying the condition
\begin{equation}%
\label{CFDe}%
2\delta<\rho_j(\boldsymbol{t}^0),
\end{equation}%
and choose and \textsl{fix} the point \(x_j\in\mathbb{C}\),
satisfying the condition
\begin{equation}
\label{GeomConfig}%
\delta<|x_j-t_j^0|<\rho_j(\boldsymbol{t}^0)-\delta .
\end{equation}%

\begin{figure}[hbt]%
\begin{minipage}[H]{0.47\linewidth}%
\centering\epsfig{figure=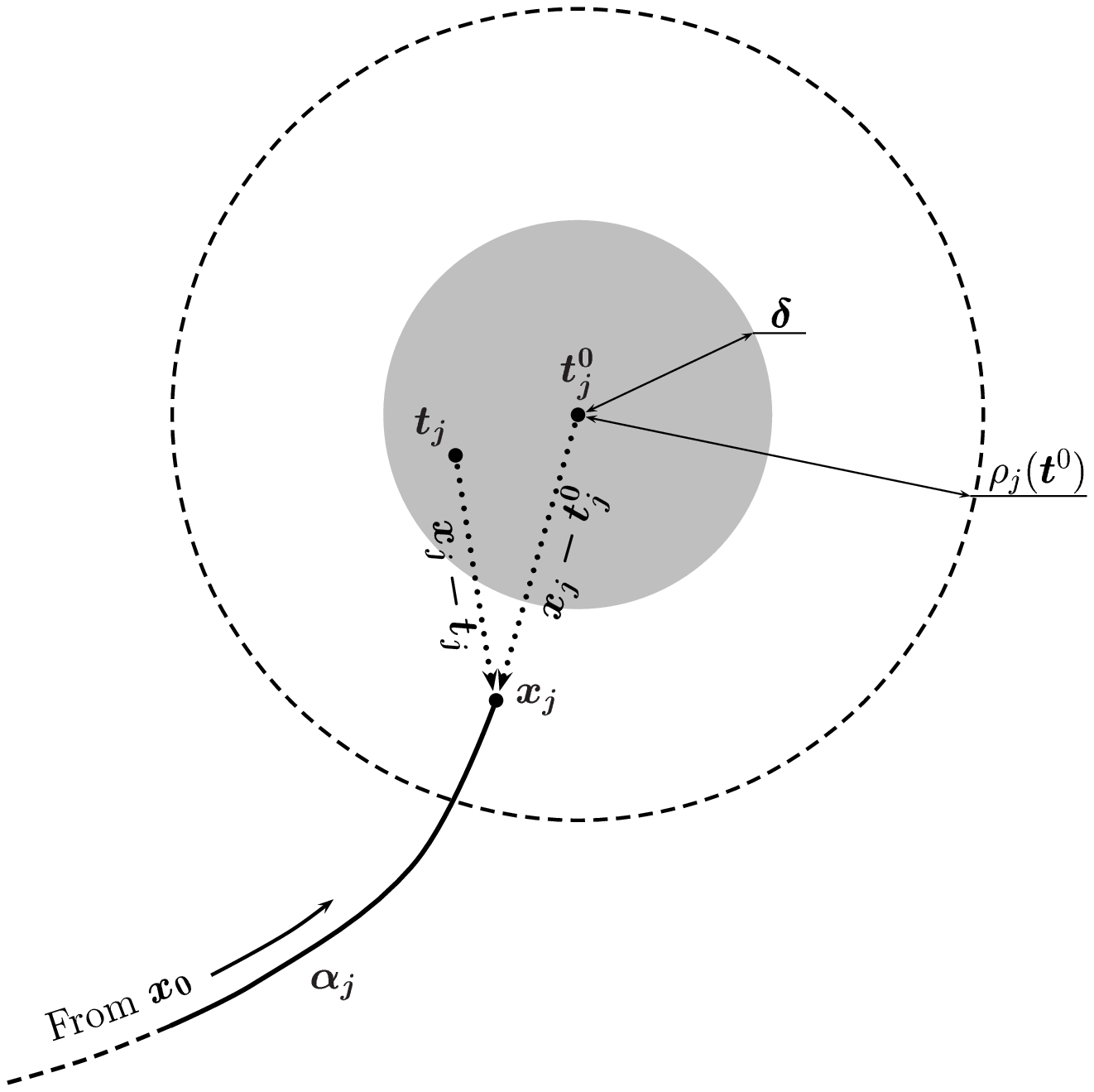,width=1.0\linewidth}%
\caption{a: disk \(D_{ t_j^0, \rho_j(\boldsymbol{t}^0)}\).}
\end{minipage}%
\addtocounter{figure}{-1}%
\begin{minipage}[H]{0.47\linewidth}%
 \centering\epsfig{figure=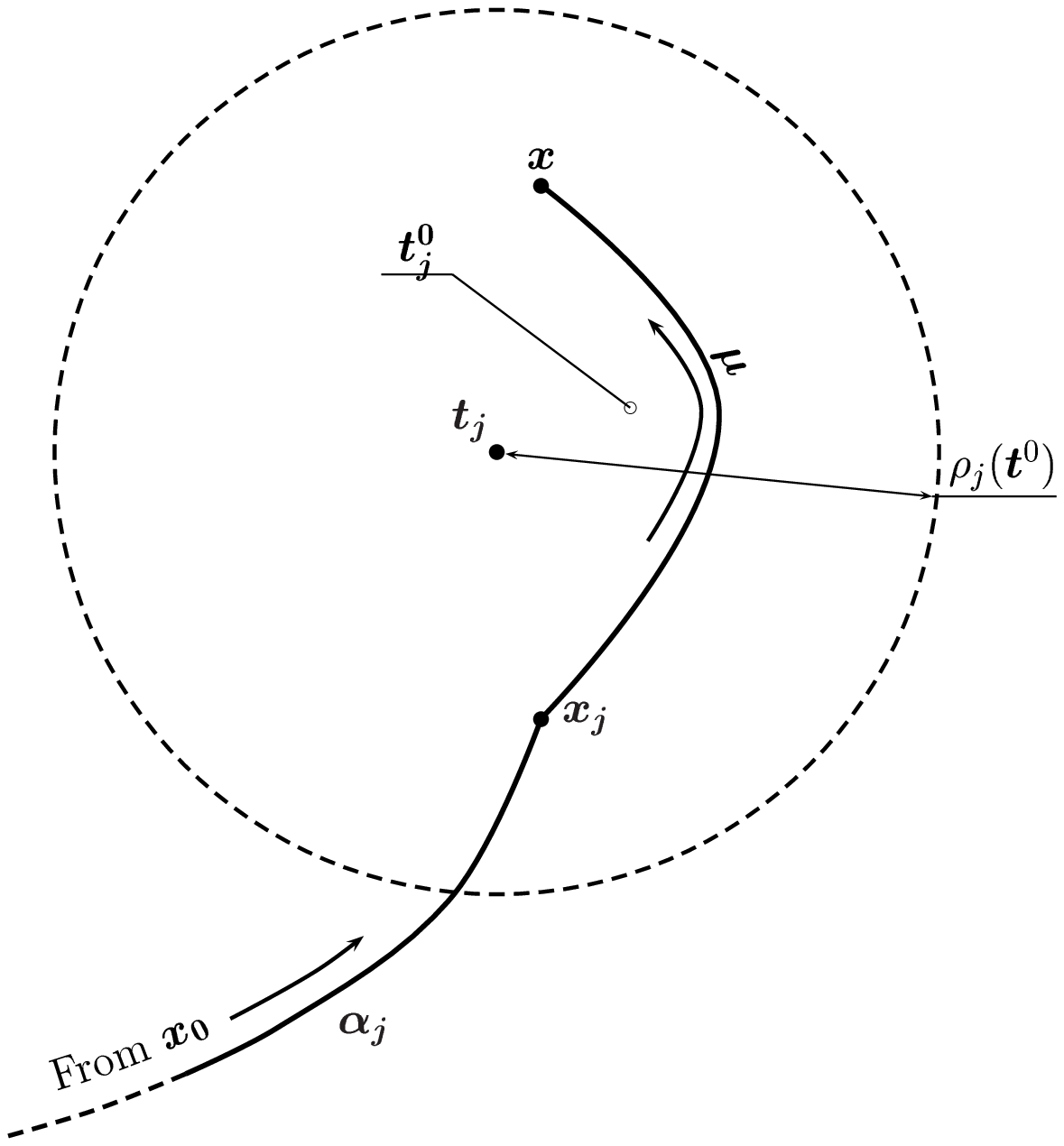,width=0.89\linewidth}%
\caption{b: disk \(D_{ t_j, \rho_j(\boldsymbol{t}^0)}\).}%
\vspace*{1.5ex} \label{CiGeo}
\end{minipage}%
\end{figure}%

With such a choice %
\begin{equation}%
\label{DesGeom}%
x_j\in D_{t_j, \rho_j(\boldsymbol{t}^0)}\quad\text{for every}\quad t_j\in
D_{t_j^0, \delta}.
\end{equation}%
(See Figure \ref{CiGeo}.a, where the disk \(D_{t_j^0, \delta}\)
is plotted in the light gray color, and the dashed circle is the
boundary of the disk \(D_{t_j^0, \rho_j(\boldsymbol{t}^0)}\)).
Choose a simple path \(\alpha_j\), leading from the distinguished
point \(x_0\) to the point \(x_j\) and separated from the points
\(t_p^0, p=1, \dots , t_n\):
\begin{equation}%
\label{Separ}%
|z-t_p^0|>\delta,\quad \forall z\in\alpha_j,\ \ p=1, \dots , n .
\end{equation}

Then \textsl{fix} the choice of the path \(\alpha_j\). (\textsl{In
particular, let the path \(\alpha_j\) be  independent of
\(\boldsymbol{t}\)}.)  From
(\ref{Separ}) it follows that for %
\footnote{ \(\mathsf{D}(\boldsymbol{t}^0, \delta)=
D_{t_1^0, \delta}\times \cdots \times D_{t_n^0, \delta}\) is a polydisk.} %
\(\boldsymbol{t}=(t_1, \dots , t_n)\in\mathsf{D}(\boldsymbol{t}^0, \delta)\),
the condition \(|x-t_p|>0\ \forall
x\in\alpha_j, p=1, \dots , n ,\) holds. Thus, the function
\(\displaystyle\sum\limits_{1\leq p\leq
n}\frac{Q_p ({\boldsymbol{t}})}{x-t_p}\) is holomorphic with
respect to \(x, \boldsymbol{t}\) for \(x\) from a neighborhood of
the path \(\alpha_j\) and  for
\(\boldsymbol{t}\in\mathsf{D}(\boldsymbol{t}^0, \delta)\).
Therefore, the solution \(Y(x,\boldsymbol{t})\) of the
differential equation (\ref{FDES}) that satisfies the initial
condition
\begin{equation}%
\label{ICondT}%
Y(x,\boldsymbol{t})_{x=x_0}=I
\end{equation}%
can be continued analytically along the path \(\alpha_j\) from a neighborhood
\(\mathcal{V}_{x_0}\) of the initial point \(x_0\) to a neighborhood
\(\mathcal{V}_{x_j}\) of the point \(x_j\). (\(\mathcal{V}_{x_j}\) does not
depend on \(\boldsymbol{t\in}\mathsf{D}(\boldsymbol{t}^0, \delta)\)). We
denote by \(Y_{\alpha_j}(x, \boldsymbol{t})\) the function in
\(\mathcal{V}_{x_j}\) which is defined by means of such an analytic
continuation. Since the initial condition in the (regular) point \(x_0\) and
the path \(\alpha_j\) do not depend on \(\boldsymbol{t}\), then for every
\(x\in\alpha_j\), the solution \(Y_{ \alpha_j}(x,\boldsymbol{t})\) is
holomorphic with respect to \(\boldsymbol{t}\) for \(x\in\mathcal{V}_{x_j}\),
\(\boldsymbol{t}\in\mathsf{D}(\boldsymbol{t}^0, \delta)\). In particular, the
matrix \(Y_{\alpha_j}(x_j,\boldsymbol{t})\), as well as the inverse matrix
\(Y^{-1}_{\alpha_j}(x_j,\boldsymbol{t})\),  is holomorphic with respect to
\(\boldsymbol{t}\) for
\(\boldsymbol{t}\in\mathsf{D}(\boldsymbol{t}^0, \delta)\). It is clear that
the solutions \(Y_{\alpha_j^{\prime}}(x, \boldsymbol{t})\) and
\(Y_{\alpha_j^{\prime\prime}}(x, \boldsymbol{t})\) corresponding to the
analytic continuations into \(\mathcal{V}_{x_j}\) along two different paths
\(\alpha_j^{\prime}\) and \(\alpha_j^{\prime\prime}\) (with the same endpoints
\(x_0\) and \(x_j\)), are related by the equality
\begin{equation}%
\label{DifP1}%
Y_{\alpha_j^{\prime\prime}}(x, \boldsymbol{t})=
Y_{\alpha_j^{\prime}}(x, \boldsymbol{t})\cdot M_{\gamma}(\boldsymbol{t}) ,
\end{equation}%
where \(M_{\gamma}(\boldsymbol{t})\) is the monodromy matrix for the equation
(\ref{FDES}) corresponding to the loop \(\gamma\),
\(\gamma=\big(\alpha_j^{\prime}\big)^{-1}\alpha_j^{\prime\prime}\), composed
from the paths \(\alpha_j^{\prime}\), \(\alpha_j^{\prime\prime}\). (The loop
\(\gamma\) does not depend on \(\boldsymbol{t}\), but the equation depends.
So, the monodromy matrix \(M_{\gamma}(\boldsymbol{t})\) can depend on
\(\boldsymbol{t}\)).

 Now we turn our attention to the local solution
\(Y_j(x,\boldsymbol{t})\), corresponding to the singular point
\(t_j\). Assume that the matrix \(Q_j(\boldsymbol{t}^0)\) is
non-resonant. Choosing a smaller \(\delta\), if necessary, we can
apply Proposition \ref{ELS}. According to this proposition, the
local solution \(Y_j(x,\boldsymbol{t})\) exists in the punctured
disk \(D_{t_j, \rho_j(\boldsymbol{t}^0)}\setminus\{t_j\}\). In
particular, it is defined at the point \(x_j\), and even in the
neighborhood \(\mathcal{V}_{x_j}\), if this neighborhood is small
enough. The local solution \(Y_j(x,\boldsymbol{t})\) is the
product of two factors: \(H_j(x,\boldsymbol{t})\) and
\((x-t_j)^{Q_j(\boldsymbol{t})}\). Under the normalizing condition
(\ref{NLS}), the first factor is determined uniquely. The second
factor is not quite unique. In general, the function
\((x-t_j)^{Q_j(\boldsymbol{t})}\) is a multivalued function of
\(x\). This function, and hence the solution
\(Y_j(x,\boldsymbol{t})\), are determined up to right factor
\(\exp\{2\pi i Q_j(\boldsymbol{t})\}\). To avoid this
non-uniqueness, we must come to an agreement on the value of
\(\arg (x-t_j)\) for some fixed \(x\), for example, for
\(x=x_j\). However, now \(\boldsymbol{t}=(t_1, \dots , t_n)\)
is not a constant but a variable with values in the polydisk
\(\mathsf{D}(\boldsymbol{t}^0,\delta)\). In particular, the
\(j\)-th coordinate \(t_j\) takes the values in the disk
\(D_{t_j^0,\rho}\). In order for the function
\((x_j-t_j)^{Q_j(\boldsymbol{t})}\) to be holomorphic with respect
to \(\boldsymbol{t}\), we need to choose \(\arg(x_j-t_j)\)
\textsl{coherently} with respect to \(t_j\).

 Let us choose somehow and \textsl{fix} the value \(\arg(x_j-t_j^0)\), say
\begin{equation}%
\label{ChoiArg0}%
\arg (x_j-t_j^0)=\vartheta_j^0 .
\end{equation}%
For \(t_j\in D_{t_j^0, \delta}\), choose the value of \(\arg(x_j-t_j)\)
according to the rule:
\begin{equation}%
\label{ChoiArg}%
\vartheta_0-\pi/2<\arg(x_j-t_j) < \vartheta_0+\pi/2 .
\end{equation}%
The  choice (\ref{ChoiArg}) is possible and unique (see Figure \ref{CiGeo}.a).
With such a choice, the matrix \((x_j-t_j)^{Q_j(\boldsymbol{t})}\) is
\textsl{determined uniquely} and is holomorphic with respect to
\(\boldsymbol{t}\) for
\(\boldsymbol{t}\in\mathsf{D}(\boldsymbol{t}^0, \delta)\). Therefore, the
matrix
\(Y_j(x_j,\boldsymbol{t})=H_j(x,\boldsymbol{t})(x_j-t_j)^{Q_j(\boldsymbol{t})}\)
is \textsl{determined uniquely} and is holomorphic with respect to
\(\boldsymbol{t}\) for
\(\boldsymbol{t}\in\mathsf{D}(\boldsymbol{t}^0, \delta)\).

Moreover, the function \(Y_j(x,\boldsymbol{t})\) --- the local
solution corresponding to the singular point \(t_j\) --- is
determined uniquely on the universal covering
\\ \(\cov (D_{t_j,\rho_j(\boldsymbol{t}^0)\setminus
\{t_j\}},x_j)\) with the distinguished point \(x_j\). The conditions
(\ref{ChoiArg0})-(\ref{ChoiArg}), together with the condition (\ref{NLS}), can be
considered as normalizing condition for the local solution
\(Y_j(x,\boldsymbol{t})\): the condition (\ref{NLS}) normalized the factor
\(H_j(x,\boldsymbol{t}))\), the condition (\ref{ChoiArg}) normalized the factor
\((x-t_j)^{Q_j(\boldsymbol{t})}\) coherently with respect to \(t_j\).

 Let us define the connection coefficient
\(C_{\alpha_j}(\boldsymbol{t})\) as
\begin{equation}%
\label{ConCoefT}%
C_{\alpha_j}(\boldsymbol{t})=Y^{-1}_j(x_j, \boldsymbol{t})\cdot
Y_{x_j, \alpha_j}(x_j, \boldsymbol{t})
\end{equation}
Since both  factors in the right hand side of (\ref{ConCoefT}) are
holomorphic with respect to
\(\boldsymbol{t}\in\mathsf{D}(\boldsymbol{t}^0, \delta)\) and
invertible, \textsl{the matrix \(C_{\alpha_j}(\boldsymbol{t})\) is
holomorphic with respect to
\(\boldsymbol{t}\in\mathsf{D}(\boldsymbol{t}^0, \delta)\) and
invertible as well}. According to (\ref{DifP1}), the connection
coefficients \(C_{\alpha_j}^{\prime}(\boldsymbol{t})\) and
\(C_{\alpha_j}^{\prime\prime}(\boldsymbol{t})\), corresponding to
two different paths \(\alpha_j^{\prime}\) and
\(\alpha_j^{\prime\prime}\), are related by the equality
\begin{equation}%
\label{DifP2}%
C_{ \alpha_j^{\prime\prime}}(x, \boldsymbol{t})=
C_{\alpha_j^{\prime}}(x, \boldsymbol{t})\cdot M_{\gamma}(\boldsymbol{t}) ,
\end{equation}%
where \(M_{\gamma}(\boldsymbol{t})\) is the monodromy matrix for the equation
(\ref{FDES}) corresponding to the loop \(\gamma\),
\(\gamma=\big(\alpha_j^{\prime}\big)^{-1}\alpha_j^{\prime\prime}\).

 Since the matrix
\(C_{\alpha_j}(\boldsymbol{t})\) does not depend on \(x\), the
product \(Y_j(x,\boldsymbol{t})\cdot C_{\alpha_j}\) is a solution
of the differential equation (\ref{FDES}) for \(x\in D_{t_j,
\rho_j(\boldsymbol{t}^0)}\), in particular, for
\(x\in\mathcal{V}_{x_j}\).  The function \(Y_{x_j,
\alpha_j}(x,\boldsymbol{t})\)  is also a solution of the
differential equation (\ref{FDES}) for \(x\in\mathcal{V}_{x_j}\).
Both these solutions coincide at the point \(x_j\). Hence they
coincide for all \(x\in\mathcal{V}_{x_j}\). Therefore, the product
\(Y_j(x,\boldsymbol{t})\cdot C_{\alpha_j}\) can be considered as
the \textsl{analytic continuation} of the solution \(Y_{x_j,
\alpha_j}(x,\boldsymbol{t})\) from \(\mathcal{V}_{x_j}\) to the
punctured disk \(D_{t_j,
\rho_j(\boldsymbol{t}^0)}\setminus\{t_j\}\) (or, better to say, to
the universal covering \(\cov(D_{t_j,
\rho_j(\boldsymbol{t}^0)}\setminus\{t_j\}; x_j)\) of this
punctured disk with the distinguished point \(x_j\)). In its turn,
the solution \(Y_{x_j, \alpha_j}(x,\boldsymbol{t})\), which is
defined  in \(\mathcal{V}_{x_j}\), was obtained by the
 described above analytic continuation of the solution,``living'' in
\(\mathcal{V}_{x_0}\) and satisfying the initial condition
(\ref{ICondT}). In other words, the solution
\(Y_j(x,\boldsymbol{t})\cdot C_{\alpha_j}\) can be considered as a
solution which is obtained as the result of the ``two step''
analytic continuation from \(\mathcal{V}_{x_0}\) to \(\cov(D_{t_j,
\rho_j(\boldsymbol{t}^0)}\setminus\{t_j\}; x_j)\) of the solution,
satisfying the initial condition (\ref{ICondT}). However, it is
better to talk not about the \textsl{two-step} analytic
continuation but about the analytic continuation \textsl{along the
composed path}. The representative of a point of the universal
covering \(\cov(D_{t_j, \rho_j(\boldsymbol{t}^0)}\setminus\{t_j\};
x_j)\) is a pair \((x, \mu)\), where \(x\) is a point of the
punctured disk \(D_{t_j,
\rho_j(\boldsymbol{t}^0)}\setminus\{t_j\}\), and
\(\mu\subset(D_{t_j, \rho_j(\boldsymbol{t}^0)}\setminus\{t_j\})\)
is a path starting at the point \(x_j\) and ending at \(x\). The
path corresponding to the above mentioned two-step analytic
continuation is the
composition%
\footnote%
{ The path \(\alpha\) starts at the point \(x_0\) and ends at the
point \(x\). Traversing the path \(\alpha\), we traverse first the
path
\(\alpha_j\), and then the path \(\mu\). See Figure \ref{CiGeo}.b}\ %
\(\alpha=\mu\cdot\alpha_j\). Thus, to describe the dependence on
\(\boldsymbol{t}\) of the solution \(Y_{\alpha_j}(x, \boldsymbol{t})\), which
is normalized at the initial point \(x_0\) by the condition
\begin{equation}%
\label{IniCoT}%
Y_{\alpha_j}(x, \boldsymbol{t})_{|x=x_0}=I ,
\end{equation}%
in the neighborhood of the singular point \(t_j\), we continue analytically
this solution from the neighborhood \(\mathcal{V}_{x_0}\) of the initial point
to the punctured neighborhood \(D_{t_j, \rho_j(\boldsymbol{t}^0)}\) along
those paths \(\alpha\) which reach first the point \(x_j\) (\textsl{\(x_j\)
does not depend on \(t_j\)} and is common for all these paths), and then go,
remaining within the punctured neighborhood
\(D_{t_j, \rho_j(\boldsymbol{t}^0)}\setminus\{t_j\}\), to the destination
points .

\begin{definition}%
\label{DefCohCont}%
Such a \textsf{family} \footnote%
{ Indexed by the parameter
\(\boldsymbol{t}\in\mathsf{D}(\boldsymbol{t}^0, \delta)\)}
\(\{Y_{\alpha_j}(x, \boldsymbol{t})\}_{\boldsymbol{t}\in
\mathsf{D}(\boldsymbol{t}^0, \delta)}\)  \textsf{of analytic continuations}
(with respect to \(x\)) of the normalized by \textup{(\ref{IniCoT})} solutions
of the equation (\ref{FDES}) from a neighborhood \(\mathcal{V}_{x_0}\) of the
initial point \(x_0\) to the family %
\(\{D_{t_j, \rho(\boldsymbol{t}^0)}\}_{t_j\in D_{t_j^0, \delta}}\) of the
punctured disk neighborhoods \(\) of the points \(t_j, t_j\in
D_{t_j^0, \delta}\), is said to be \textsf{coherent} with respect to
\(\boldsymbol{t}\). The path \(\alpha_j\) along which such continuations are
carried out is said to be \textsf{the path determining the coherent family of
analytic continuations}.
\end{definition}%

\noindent The solution
\begin{equation}%
\label{BedForm}%
Y_j(x, \boldsymbol{t})\cdot C_{\alpha_j}(\boldsymbol{t})=
H_j(x, \boldsymbol{t})(x-t_j)^{Q_j(\boldsymbol{t})}C_{\alpha_j}(\boldsymbol{t})
\end{equation}%
can be presented in another  form -- more suitable for our goal:
\begin{equation}%
\label{GoodForm}%
Y_{\alpha_j}(x, \boldsymbol{t})=
H_{\alpha_j}(x, \boldsymbol{t})(x-t_j)^{A_{\alpha_j}(\boldsymbol{t})} ,
\end{equation}%
where
\begin{equation}%
\label{ExpEntr}%
\textup{a).}\
A_{\alpha_j}(\boldsymbol{t})=C_{\alpha_j}^{-1}(\boldsymbol{t})Q_j(\boldsymbol{t})
C_{\alpha_j}(\boldsymbol{t}) ;\quad\quad \textup{b).}\
H_{\alpha_j}(x, \boldsymbol{t})=H_j(x, \boldsymbol{t})C_{\alpha_j}(\boldsymbol{t}) .
\end{equation}%
Of course, the relation (\ref{DifP2}) leads immediately to the relations
\begin{equation}%
\label{DifP3}%
\textup{a).}\ A_{\alpha_j^{\prime\prime}}(\boldsymbol{t})=
(M_{\gamma}(\boldsymbol{t}))^{-1}A_{{\alpha_j}^{\prime}}(\boldsymbol{t})
M_{\gamma}(\boldsymbol{t}),\quad\quad \textup{b).}\
H_{\alpha_j^{\prime\prime}}(x, \boldsymbol{t})=
H_{\alpha_j^{\prime}}(x, \boldsymbol{t}) M_{\gamma}(\boldsymbol{t}) ,
\end{equation}%
where \(\alpha_j^{\prime}\) and \(\alpha_j^{\prime\prime}\) are two paths
with the same endpoints \(x_0\) and \(x_j\), and \(M_{\gamma}\) is the
monodromy matrix for the equation (\ref{FDES}) corresponding to the loop
\(\gamma\),
\(\gamma=\big(\alpha_j^{\prime}\big)^{-1}\alpha_j^{\prime\prime}\).

\begin{definition}%
\label{PrReg}%
The matrix functions \(H_{\alpha_j}(x, \boldsymbol{t})\) and
\((x-t_j)^{A_{\alpha_j}(\boldsymbol{t})}\) are said to be the \textsl{regular
factor} and the \textsf{principal factor} of the representation
\textup{(\ref{GoodForm})} of the solution \(Y_{\alpha_j}(x, \boldsymbol{t})\)
near the point \(t_j\).
\end{definition}%

\noindent Let us summarize the  results, obtained above, as

\begin{theorem}%
\label{DepOnT}%
Let the \(k\times k\) matrix functions
\(Q_p(\boldsymbol{t}), p=1, \dots , n,\) which appear in the
Fuchsian system \textup{(\ref{FDES})} be  holomorphic with respect
to \(\boldsymbol{t}=(t_1, \dots , t_n)\) for \(\boldsymbol{t}\)
from an open set
\(\boldsymbol{\mathcal{D}}, \boldsymbol{\mathcal{D}}\in\mathbb{C}^n_{\ast}.\)
Let \(\boldsymbol{t}^0=(t_1^0, \dots , t_n^0)\) be a point from
\(\boldsymbol{\mathcal{D}}\). Assume that, for certain \(j\), the
matrix \(Q_j(\boldsymbol{t}^0)\) is non-resonant, and let
\(\rho_j(\boldsymbol{t}^0)\) be defined according to
(\ref{DefRho})) (with \(q=j\)).

Then there exists a number \(\delta\),
\(0<\delta<\frac{1}{2}\rho_j(\boldsymbol{t}^0)\), such that, for
\(\boldsymbol{t}\in\mathsf{D}(\boldsymbol{t}^0, \delta)\), the solutions
\(\{Y_{\alpha_j}(x, \boldsymbol{t})\}_{\boldsymbol{t}\in
\mathsf{D}(\boldsymbol{t}^0, \delta)}\) of the holomorphic family of equation
(\ref{FDES}) which are distinguished by the properties:
\begin{enumerate}
\item[i.]
 Each of these solutions is normalized by \textup{(\ref{IniCoT})} at the
initial point \(x_0\) ,%
\item[ii]
 Each of these solution is continued analytically from a
neighborhood of \(x_0\) to the family %
\(\{D_{t_j, \rho(\boldsymbol{t}^0)}\}_{t_j\in D_{t_j^0, \delta}}\) of the
punctured disk neighborhoods \(\) of the points \(t_j, t_j\in
D_{t_j^0, \delta}\), and such analytic continuations are performed coherently
with respect to
\(\boldsymbol{t}\in\boldsymbol{t}\in\mathsf{D}(\boldsymbol{t}^0, \delta)\) ,
\end{enumerate}%
are representable in \(D_{t_j,\rho(\boldsymbol{t}^0)}\setminus\{t_j\}\) in the
form \textup{(\ref{GoodForm})}, where:
\begin{enumerate}%
\item
 The matrix function
\(H_{\alpha_j}(x, \boldsymbol{t})\) --- the regular factor of the
representation \textup{(\ref{GoodForm})}  ---    is
 holomorphic together with its inverse
\(H_{ \alpha_j}^{-1}(x, \boldsymbol{t})\)
 with respect to
\(x, \boldsymbol{t}\) for \(x\in D_{t_j, \rho(\boldsymbol{t}^0)}, \),
\(\boldsymbol{t}\in\mathsf{D}(\boldsymbol{t}^0, \delta)\);
\item
 The exponent \(A_{\alpha_j}(\boldsymbol{t})\) of the principal factor
  \((x-t_j)^{A_{\alpha_j}(\boldsymbol{t})}\)
  is a matrix
function holomorphic with respect to \(\boldsymbol{t}\) from the polydisk
\(\boldsymbol{t}\in\mathsf{D}(\boldsymbol{t}^0, \delta)\);
\item
The exponent \(A_{\alpha_j}(\boldsymbol{t})\) is similar to the matrix
\(Q_j(\boldsymbol{t})\), \textup(\ref{ExpEntr}.a), and the similarity matrix
\(C_{\alpha_j}(\boldsymbol{t})\) can be chosen %
\footnote{The similarity matrix is not unique: it can be
multiplied by an arbitrary
matrix commuting with \(Q_j(\boldsymbol{t})\).}\ %
to be holomorphic for
\(\boldsymbol{t}\in\mathsf{D}(\boldsymbol{t}^0, \delta)\) .
\item
By a different choice  of the paths \(\alpha_j^{\prime}\) and
\(\alpha_j^{\prime\prime}\), along which the coherent analytic
continuations are carried out, the matrices
\(A_{\alpha_j}(\boldsymbol{t})\) and
\(H_{\alpha_j}(x, \boldsymbol{t})\) are transformed according to the rule
\textup{(\ref{DifP3})}, where \(M_{\gamma}(\boldsymbol{t})\) is the monodromy
matrix corresponding to a loop \(\gamma\), constructed %
\footnote%
{\label{DiPa} Both paths,\(\alpha_j^{\prime}\) and
\(\alpha_j^{\prime\prime}\),
 have the common starting point \(x_0\).
If their endpoints coincide, then the loop \(\gamma\) is just the composition
\(\gamma=(\alpha_j^{\prime})^{-1}\alpha_j^{\prime\prime}\).
 In the general case, in which these endpoints may be different,
 \(\gamma=(\alpha_j^{\prime})^{-1}\nu \alpha_j^{\prime\prime}\),
 where \(\nu\) is a path, connecting the endpoints.}\ %
 from the paths
\(\alpha_j^{\prime}\) and \(\alpha_j^{\prime\prime}\).
\end{enumerate}%
\end{theorem}%

\noindent Let \(\alpha_j\) be a path, separated (as in (\ref{Separ})) from the
singular points \(t_p, 1\leq p\leq n\), and leading from the point \(x_0\) to a
point \(x_j\in\big(\bigcap D_{t_j, \rho_j(\boldsymbol{t}^0)}\big)\setminus
D_{t_j^0, \delta}\), where the intersection is taken over all \(t_j\in
D_{t_j^0\!, \delta}\). This path can be used for the coherent analytic
continuations of the normalized by (\ref{IniCoT}) solution to the neighborhoods of
the point \(t_j\in D_{t_j^0, \delta}\). The path \(\alpha_j\) also generates a loop
\(\gamma_{\alpha_j}\) which goes around the disk neighborhood
\(D_{t_j^0, \delta}\). This loop does not depend on \(\boldsymbol{t}\) and, for
every \(\boldsymbol{t}=(t_1, \dots , t_n)\),
\(\boldsymbol{t}\in\mathsf{D}(\boldsymbol{t}^0, \delta)\), represents an element of
the fundamental group \(\pi(\overline{\mathbb{C}}\setminus
\{t_1, \dots , t_n\}; x_0 )\). (Compare to the Definition \ref{LoGen}). The loop
\(\gamma_{\alpha_j}\) is the composition of the path \(\alpha_j\), of the circular
arc \(\mu\), starting and ending at the point \(x_j\), and of the path
\(\alpha_j^{-1}\) (the path \(\alpha_j^{-1}\) is the path \(\alpha_j,\) traversed in
the ``opposite direction''). In more detail, \(\mu\) is the circular arc
 of the radius \(|x_j-t_j^0|\), centered at the point \(t_j^0\), which starts
at the point \(x_j\), makes the full  counterclockwise turn around
the point \(t_j^0\), and ends at the same point \(x_j\) (or,
better to say, at the point \(t_j^0+e^{2\pi i}(x_j-t_j^0)\)).
Since \(x_j\in\big(\bigcap
D_{t_j, \rho_j(\boldsymbol{t}^0)}\big)\setminus
D_{t_j^0, \delta}\), \(\mu\subset
D_{t_j, \rho_j(\boldsymbol{t}^0)}\setminus\{t_j\}\) for every
\(t_j\in D_{t_j^0, \delta}\). Since in
\(D_{t_j, \rho_j(\boldsymbol{t}^0)}\setminus\{t_j\}\) the
representation (\ref{GoodForm}) holds,
 Theorem \ref{EMMT} can be formulated in the following
 ``\(\boldsymbol{t}\)-dependent'' version:

\begin{theorem}%
\label{TDV}%
Assume that the matrix \(Q_j(\boldsymbol{t}^0)\) is non-resonant.
Let \(\alpha_j\) be a path leading from the initial point \(x_0\)
to a punctured neighborhood of the singular point \(t_j^0\), and
let \(\gamma_{\alpha_j}\) be the loop generated by the path
\(\alpha_j\) (as described above).

If  \(\boldsymbol{t}\in\mathsf{D}(\boldsymbol{t}^0, \delta)\),
and \(\delta>0\) is small enough, then the monodromy matrix
\(M_{\gamma_{\alpha_j}}(\boldsymbol{t})\) for the equation
\textup{(\ref{FDES})}, corresponding to the loop
\(\gamma_{\alpha_j}\), and the exponent
\(A_{\alpha_j}(\boldsymbol{t})\) of the principal factor of the
representation \textup{(\ref{GoodForm})} of its solution
\(Y_{\alpha_j}(x, \boldsymbol{t})\) are related by the equality
\begin{equation}%
\label{TDE}%
M_{\gamma_{\alpha_j}}(\boldsymbol{t})=\exp\{2\pi
iA_{\alpha_j}(\boldsymbol{t})\} .
\end{equation}%
\end{theorem}%

\noindent The question arises, \textsl{whether the monodromy
matrix \(M_{\gamma_{\alpha_j}}(\boldsymbol{t})\) determines the
exponent \(A_{\alpha_j}(\boldsymbol{t})\) uniquely?} In
particular, \textsl{does the equality
\(A_{\alpha_j}(\boldsymbol{t}^{\prime})=
A_{\alpha_j}(\boldsymbol{t}^{\prime\prime})\) follow from the
equality \(M_{\gamma_{\alpha_j}}(\boldsymbol{t}^{\prime})=
M_{\gamma_{\alpha_j}}(\boldsymbol{t}^{\prime\prime})\), if
\(\boldsymbol{t}^{\prime}, \boldsymbol{t}^{\prime\prime}
\in\mathsf{D}(\boldsymbol{t}^0, \delta)\), and \(\delta\in \) is
small enough?}

The latter question is closely related to the description of the set of
solutions of the matrix equation (with respect to \(A\)):
\begin{equation}%
\label{MatEq}%
\exp\{2\pi i A\}=M,
\end{equation}%
where \(M\) is a given \(k\times k\) matrix. A description of the
set of solutions of (\ref{MatEq}) can be found in \cite{Gant},
Chapter VIII, \S 8. The invertibility of \(M\) is the necessary
and sufficient condition for the solvability of (\ref{MatEq}).
Given \(M, \det M\not=0\), the set  of solutions of (\ref{MatEq})
 is always infinite. Our prime interest is in conditions under
which a certain solution \(A_0\) of the equation (\ref{MatEq}) is
isolated. In general, non-isolated solutions may exist. For
example, let \(M=\textup{diag} [1, 1]\) ,
\(A_0=\textup{diag} [0, 1]\), and let \(T\) be an
\textsl{arbitrary} invertible matrix. Then \(\exp\{2\pi i
(T^{-1}A_0T)\}=M\). If \(T\) is non-diagonal matrix, then
\(T^{-1}A_0T\not= A_0\), but if \(T\) is close to \(I\), then
\(T^{-1}A_0T\) is close to \(A_0\).

\begin{lemma}%
\label{IsSol}%
Given a non-resonant matrix \(A_0\), there exists \(\varepsilon>0\) such that
from the conditions \(\exp\{2\pi i A_0\}=\exp\{2\pi i A\}\) and
\(\|A-A_0\|<\varepsilon\) it follows that \(A=A_0\).
\end{lemma}%

\noindent%
\textsf{PROOF}.%
Since \(A_0\) is non-resonant, there exists an open (not necessarily
connected) set \(\mathcal{O}, \mathcal{O}\subset\mathbb{C}\)
containing the spectrum of \(A_0\) such that the mapping
\(\zeta\mapsto\exp\{2\pi i\zeta\}\) is univalent (schlicht) in
\(\mathcal{O}\). If \(\varepsilon\) is small enough, the set
\(\mathcal{O}\) also contains the spectra of all matrices \(A\)
satisfying the conditions \(\|A-A_0\|<\varepsilon\). The set
\(\mathcal{G}=\exp\{2\pi i \mathcal{O}\}\) is an open set
containing all the spectra of the matrices \(\exp\{2\pi i A\}\)
with \(\|A-A_0\|<\varepsilon\), in particular, the spectrum of the
matrix  \(\exp\{2\pi i A_0\}\). If \(\dfrac{1}{2\pi
i} \textup{Ln} w\), \(\dfrac{1}{2\pi
i} \textup{Ln}: \mathcal{G}\mapsto\mathcal{O}\), is a function,
inverse  to the function \(\exp\{2\pi i
\zeta\}: \mathcal{O}\mapsto\mathcal{G}\), that is
\(\dfrac{1}{2\pi i} \textup{Ln}(\exp\{2\pi i\zeta\})\equiv\zeta\
\ \forall\zeta\in\mathcal{O}\), then \(\dfrac{1}{2\pi
i} \textup{Ln}(\exp\{2\pi iA\})\equiv A\) for all matrices \(A\)
whose spectra are contained in \(\mathcal{O}\), in particular,
\(\forall \, A: \|A-A_0\|<\varepsilon\). \hfill Q.E.D.\\[1ex]

\noindent Since the matrix \(A_{\alpha_j}(\boldsymbol{t}^0)\) is
similar to the matrix \(Q_j(\boldsymbol{t}^0)\), the following
result is the direct consequence of Lemma \ref{IsSol} and Theorem
\ref{TDV}:

\begin{theorem}%
\label{UOE}%
Let the matrix \(Q_j(\boldsymbol{t}^0)\) be non-resonant, and let
\(A_{\alpha_j}(\boldsymbol{t})\) and \(M_{\gamma_{\alpha_j}}(\boldsymbol{t})\)
be the same as in the formulation of Theorem \ref{TDV}. If \(\delta>0\) is
small enough, and if for some
\(\boldsymbol{t}^{\prime}, \boldsymbol{t}^{\prime\prime}\in
\mathsf{D}(\boldsymbol{t}^0, \delta)\), the equalities
\(M_{\gamma_{\alpha_j}}(\boldsymbol{t^{\prime}})=
M_{\gamma_{\alpha_j}}(\boldsymbol{t}^{\prime\prime})\) hold, then the
equalities \(A_{\alpha_j}(\boldsymbol{t}^{\prime})=
A_{\alpha_j}(\boldsymbol{t}^{\prime\prime})\) hold as well. In particular, if
\(M_{\gamma_{\alpha_j}}(\boldsymbol{t})\equiv
M_{\gamma_{\alpha_j}}(\boldsymbol{t}^0)\) in
\(\mathsf{D}(\boldsymbol{t}^0, \delta)\), then
\(A_{\alpha_j}(\boldsymbol{t})\equiv A_{\alpha_j}(\boldsymbol{t}^0)\) in
\(\mathsf{D}(\boldsymbol{t}^0, \delta)\) as well.
\end{theorem}%

\begin{remark}%
\label{GRC}%
If the matrix \(Q_j(\boldsymbol{t}^0)\) is resonant, then, trying to formulate
a ``\(\boldsymbol{t}\)-dependent'' version of Theorem \ref{SRQI}, we should
seek the representation of the solution \(Y_{\alpha_j}(\boldsymbol{t})\) near
\(\boldsymbol{t}^0\) in the form
\begin{equation}%
\label{SRQIT}%
Y_{\alpha_j}(x, \boldsymbol{t})= H_{\alpha_j}(x, \boldsymbol{t})
(x-t_j)^{Z_{\alpha_j}(\boldsymbol{t})}(x-t_j)^{A_{\alpha_j}(\boldsymbol{t})}
, \ x\in \cov\big(D_{t_j,\rho}\setminus\{t_j\},x_j\big) ,
\end{equation}%
where the matrices
\(H_{\alpha_j}(x, \boldsymbol{t}), H_{\alpha_j}^{-1}(x, \boldsymbol{t})\)
are holomorphic for \(x\in D_{t_j,\rho(\boldsymbol{t}^0)}\) and
\(\boldsymbol{t}
\in\mathsf{D}(\boldsymbol{t}^0, \rho_j(\boldsymbol{t}^0))\),  the
matrices
\(Z_{\alpha_j}(\boldsymbol{t}), A_{\alpha_j}(\boldsymbol{t})\)
are holomorphic for \(\boldsymbol{t}
\in\mathsf{D}(\boldsymbol{t}^0, \rho_j(\boldsymbol{t}^0))\), and
moreover the matrices \(Z_{\alpha_j}(\boldsymbol{t})\) are
diagonalizable , the eigenvalues of
\(Z_{\alpha_j}(\boldsymbol{t})\) are integer, the eigenvalues of
\(A_{\alpha_j}(\boldsymbol{t})\) and the eigenvalues of
\(Q_j(\boldsymbol{t})\) are related by equalities analogous to the
equalities \textup{(\ref{EVD})}. However, if the matrix
\(Q_j(\boldsymbol{t}^0)\) is resonant, the representation of the
form \textup{(\ref{SRQIT})}, with \textsf{holomorphic}
\(Z_{\alpha_j}(\boldsymbol{t}), A_{\alpha_j}(\boldsymbol{t})\),
does not exist in general. The reason is that it is possible for
the matrix \(Q_j(\boldsymbol{t}^0)\) to be resonant, while all the
matrices \(Q_j(\boldsymbol{t})\) for
\(\boldsymbol{t}\not=\boldsymbol{t}^0\) are non-resonant. (The
property of  resonance is unstable under small perturbations).
According to Theorem \ref{SNRQI}, the solution
\(Y_{\alpha_j}(\boldsymbol{t})\) must be of the form
\textup{(\ref{GoodForm})} for
\(\boldsymbol{t}\not=\boldsymbol{t}^0\). According to Theorem
\ref{SRQI}, it must be of the form \textup{(\ref{SRQIT})} for
\(\boldsymbol{t}=\boldsymbol{t}^0\), where in general
\(Z_{\alpha_j}(\boldsymbol{t}^0)\not=0\). However, such a
bifurcation is incompatible with the holomorphy of
\(Z_{\alpha_j}(\boldsymbol{t})\).
\end{remark}%

\section{ISOPRINCIPAL AND ISOMONODROMIC DEFORMATIONS.\label{IsoDef}}
Consider a Fuchsian equation (\ref{FDES}), where \(k\times k\)
matrices \(Q_p(\boldsymbol{t}), 1\leq p\leq n,\)  are holomorphic
in \(\boldsymbol{t}\) for \(\boldsymbol{t}\) from a domain
\(\boldsymbol{\mathcal{D}}\), \( \boldsymbol{\mathcal{D}}\subset
\mathbb{C}^n_{\ast}\). The equation (\ref{FDES}) is considered as
a differential equation with respect to the variable \(x\), and
\(\boldsymbol{t}=(t_1, \dots , t_n)\) is considered as a
parameter. Recall that we use the terminology ``\textsl{the
deformation of  Fuchsian equation}''
for a Fuchsian equation depending on a parameter%
\footnote{ We consider only such deformations, where the loci
\(\boldsymbol{t}=(t_1, \dots , t_n)\) serve as a parameter, and the
dependence on the parameter is holomorphic.}.%

\begin{definition}
\label{DIMonD}%
Let (\ref{FDES}) be a deformation of  Fuchsian equation, where
\(k\times k\) matrices \(Q_p(\boldsymbol{t}), 1\leq p\leq n,\) are
holomorphic in a domain \(\boldsymbol{\mathcal{D}},
\boldsymbol{\mathcal{D}}\subset\mathbb{C}^n_{\ast}\), and let
\(\boldsymbol{t}^0=(t_1^0, \dots , t_n^0)\) be a point of
\(\boldsymbol{\mathcal{D}}\). The deformation is said to be
\textsf{isomonodromic with respect to a distinguished point
\(x_0\) at the pole loci \(\boldsymbol{t}^0\)} if for every loop
\(\gamma\) with the distinguished point \(x_0\),
\(\gamma\subset\mathbb{C}\setminus\{t_1^0, \dots , t_n^0\}\),
there exists \(\delta>0\),
\(\delta=\delta(\boldsymbol{t}^0,\gamma)\), such that the
monodromy matrix \(M_{\gamma}(\boldsymbol{t})\) of the deformation
(\ref{FDES}) along the loop \(\gamma\) does not depend on
\(\boldsymbol{t}\) in \(\mathsf{D}(\boldsymbol{t}^0,\delta)\):
\(M_{\gamma}(\boldsymbol{t})\equiv M_{\gamma}(\boldsymbol{t}^0)\,
\forall\,\boldsymbol{t}\in\mathsf{D}(\boldsymbol{t}^0,\delta).\)
(Of course, we assume that
\(\delta<\textup{dist}(\boldsymbol{t}^0, \gamma)\), so that
\(\gamma\subset\mathbb{C}\setminus\{t_1, \dots , t_n\}\) for
\(\boldsymbol{t}\in\mathsf{D}(\boldsymbol{t}^0,\delta)\), and the
monodromy matrix \(M_{\gamma}(\boldsymbol{t})\) is well defined
for such \(\boldsymbol{t}\)). The deformation is said to be
\textsf{isomonodromic with respect to \(x_0\) in
\(\boldsymbol{\mathcal{D}}\)} if it is isomonodromic with respect
to \(x_0\) at every pole loci
\(\boldsymbol{t}\in\boldsymbol{\mathcal{D}}\).
\end{definition}

\begin{remark}%
\label{IMF}%
To check that the deformation \textup{(\ref{FDES})} is
isomonodromic at \(\boldsymbol{t}^0\), there is no need to examine
the monodromy matrix \(M_{\gamma}(\boldsymbol{t})\) for
\textup{all} loops \(\gamma\). It is enough to choose loops
\(\gamma_j, j=1, \dots , n ,
\gamma_j\subset\mathbb{C}\setminus\{t_1^0, \dots , t_n^0\}\), generating %
\footnote{ If \(\boldsymbol{t}\in\mathsf{D}(\boldsymbol{t}^0,\delta)\) and
\(\delta>0\) is small
 enough, then such loops \(\gamma_j\) satisfy the condition
 \(\gamma_j\subset\mathbb{C}\setminus\{t_1, \dots , t_n\}\) and generate
 the fundamental group
 \(\pi(\mathbb{C}\setminus\{t_1, \dots , t_n\},x_0)\). See Figure \ref{BGHV}.}%
 the fundamental group
\(\pi(\mathbb{C}\setminus\{t_1^0, \dots , t_n^0\},x_0)\), and
 to check the the monodromy matrices \(M_{\gamma_j}(\boldsymbol{t})\) do not
 depend on \(\boldsymbol{t}\) for
 \(\boldsymbol{t}\in\mathsf{D}(\boldsymbol{t}^0,\delta)\), \(\delta>0\) is small
 enough.
\end{remark}%

\begin{remark}%
\label{LGlIs}%
If \(\gamma\) is a loop with the distinguished point \(x_0\), and
\(t_j\not\in\gamma, j=1, \dots , n\), then the monodromy
matrix \(M_{\gamma}(\boldsymbol{ . })\) is holomorphic at the
point \(\boldsymbol{t}=(t_1, \dots , t_n)\). From the
uniqueness theorem for holomorphic functions and from elementary
topological considerations it follows that if the deformation
\textup{(\ref{FDES})} is isomonodromic at some
\(\boldsymbol{t}\in\boldsymbol{\mathcal{D}}\), then it is
isomonodromic at every
\(\boldsymbol{t}\in\boldsymbol{\mathcal{D}}\).
\end{remark}%

\begin{definition}%
\label{DefIsosp}%
Deformation \textup{(\ref{FDES})} is said to be
\textit{isospectral}, if for every \(j,\ j=1, \dots , n,\) the
spectrum of the matrix \(Q_j(\boldsymbol{t})\) does not depend on
\(\boldsymbol{t}\) in \(\boldsymbol{\mathcal{D}}\): \(\spec
Q_j(\boldsymbol{t}^{\prime})= \spec
Q_j(\boldsymbol{t}^{\prime\prime})\) (the sets coincides ``with
multiplicities'') for every \(\boldsymbol{t}^{\prime},
 \boldsymbol{t}^{\prime\prime}\in\boldsymbol{\mathcal{D}}\).
\end{definition}%

\begin{theorem}%
\label{IMIS}%
If a deformation \textup{(\ref{FDES})} is isomonodromic with
respect to a distinguished point, then it  is also isospectral.
\end{theorem}%

\noindent \textsf{PROOF}. For given \(j\), the coefficients of the
characteristic polynomials of the matrix \(Q_j(\boldsymbol{t})\)
are holomorphic with respect to \(\boldsymbol{t}\) in
\(\boldsymbol{\mathcal{D}}\). Therefore, it is enough to prove
that for some \(\delta>0\), this polynomial does not depend on
\(\boldsymbol{t}\) for
\(\boldsymbol{t}\in\mathsf{D}(\boldsymbol{t}^0,\delta)\) or, what
is the same, that the spectrum \(\spec Q_j(\boldsymbol{t})\) does
not depend on \(\boldsymbol{t}\) for
\(\boldsymbol{t}\in\mathsf{D}(\boldsymbol{t}^0,\delta)\). Let
\(\gamma_j\) be a loop with the distinguished point \(x_0\) wich
makes one turn counterclockwise around the singular point
\(t_j^0\) and makes no turns around other singular points
\(t^0_p,\ p=1, \dots , n, p\not=j\), as it was described in the
formulation of Theorem \ref{SpMaTe}. The monodromy matrix
\(M_{\gamma_j}(\boldsymbol{t})\) does not depend on
\(\boldsymbol{t}\) for \(\boldsymbol{t}\) which are close to
\(\boldsymbol{t}^0\) (In general,
\(M_{\gamma_j}(\boldsymbol{t})=M_{\gamma_j}(\boldsymbol{t}^0)\) if
it is possible to pass from \(\boldsymbol{t}^0\) to
\(\boldsymbol{t}\) in \(\boldsymbol{\mathcal{D}}\) so that no
\(t_p\) crosses the loop \(\gamma_j\)). According to Theorem
\ref{SpMaTe}, the set \(\exp\{\,\spec Q_j(\boldsymbol{t})\}\) does
not depend on \(\boldsymbol{t}\) in
\(\mathsf{D}(\boldsymbol{t}^0,\delta)\). However, on the one hand,
the spectrum \(\spec Q_j(\boldsymbol{t})\) depends on
\(\boldsymbol{t}\) continuously, and on the other hand, the
mapping \(\zeta\mapsto\exp\{\zeta\}\) is locally univalent
(schlicht). Therefore, the set \(\spec Q_j(\boldsymbol{t})\) does
not depend on \(\boldsymbol{t}\) in
\(\mathsf{D}(\boldsymbol{t}^0,\delta)\) if
\(\delta\) is small enough. \hfill Q.E.D.\\[1ex]

\noindent Let (\ref{FDES}) be an isomodromic with respect to
\(x_0\) deformation of Fuchsian equation, and let
for a certain \(j\), the matrix \(Q_j(\boldsymbol{t}^0)\) be non-resonant %
\footnote{ According to Theorem \ref{IMIS}, the matrix
\(Q_j(\boldsymbol{t})\) in
non-resonant for every \(\boldsymbol{t}\in\boldsymbol{\mathcal{D}}\).}. %
Let \(\alpha_j\),
\(\alpha_j\subset\mathbb{C}\setminus\{t_1^0, \dots , t_n^0\}\),
be a path leading from the point \(x_0\) to a small neighborhood
of the point \(t_j^0\). According to Theorem \ref{UOE}, for
\(\boldsymbol{t}\in\mathsf{D}(\boldsymbol{t}^0,\delta)\), the
solution \(Y_{\alpha_j}(x,\boldsymbol{t})\), normalized at the
distinguished point \(x_0\) by the condition (\ref{ICondT}) and
continued analytically to the neighborhood of \(t_j^0\), is
representable near \(t_j\) in the form
\begin{equation}%
\label{LIMD}%
Y_{\alpha_j}(x,\boldsymbol{t})=
H_{\alpha_j}(x,\boldsymbol{t})(x-t_j)^{A_{\alpha_j}},\quad x\in
D_{t_j, \rho},  \boldsymbol{t}\in\mathsf{D}(\boldsymbol{t}^0,\delta) ,
\end{equation}%
where \(\delta\) and \(\rho\) are some positive numbers, the
\textsl{regular factor} \(H(x,\boldsymbol{t})\) is holomorphic in
\(x,\boldsymbol{t}\) and invertible for \(x\in D_{t_j, \rho}\),
\(\boldsymbol{t}\in\mathsf{D}(\boldsymbol{t}^0,\delta)\), and the
exponent \(A_{\alpha_j}\) of the \textsl{principal factor}
\((x-t_j)^{A_{\alpha_j}}\) \textsf{\textsl{does not depend on}}
\(\boldsymbol{t}\). Thus it is natural to give the following
definition:

\begin{definition}%
\label{IDNS}%
The deformation \textup{(\ref{FDES})} is said to be
\textsf{isoprincipal with respect to a distinguished point \(x_0\)
in the narrow sense} at the pole loci
\(\boldsymbol{t}^0=(t_1^0, \dots , t_n^0)\),
\((\boldsymbol{t}^0\in\boldsymbol{\mathcal{D}}\)), if for every
\(j=1, \dots , n\) and for some paths
\(\alpha_1, \dots , \alpha_n\) leading from \(x_0\) to
neighborhoods of the points \(t_1^0, \dots , t_n^0\)
respectively, the solution \(Y_{\alpha_j}(x,\boldsymbol{t})\)
normalized at \(x_0\) by \textup{(\ref{ICondT})} is representable
in the form \textup{(\ref{LIMD})} near \(t_j\), where the matrices
\(A_{\alpha_j}\) do not depend on \(\boldsymbol{t}\) for
 \(\boldsymbol{t}\in\mathsf{D}(\boldsymbol{t}^0,\delta)\), \(\delta>0\) is small
 enough.
\end{definition}%

\noindent We emphasize that we do not require in this definition
that the matrices \(Q_j(\boldsymbol{t}^0)\) are non-resonant. Of
course, in order for this definition to be correct, we must  check
that it does not depend on the choice of the paths \(\alpha_j\).
This can be done easily. Below (Lemma \ref{IDWD}) we prove this
independence in a more general situation. Let us summarize the
results in the form of two lemmas:

\begin{lemma}%
\label{ImIpr}%
In the non-resonant case (i.e., if all the matrices \(Q_j(\boldsymbol{t}^0)\)
are non-resonant), every isomonodromic deformation is isoprincipal in the
narrow  sense as well.
\end{lemma}%

\noindent This lemma is the immediate consequence of Theorem
\ref{UOE}. The converse is true always, without any assumption of
non-resonance.

\begin{lemma}%
\label{IpIso}%
If the deformation \textup{(\ref{FDES})} is isoprincipal (in the
``narrow sense'') at some  \(\boldsymbol{t}^0\) (no assumptions on
the matrices \(Q_j(\boldsymbol{t}^0), j=1, \dots , n\) are
made), then this deformation is isomonodromic as well.
\end{lemma}%

\noindent%
\textsf{PROOF}. %
 If the solution \(Y_{\alpha_j}(x,\boldsymbol{t})\) of
(\ref{FDES}) is representable in the form (\ref{LIMD}) with a
constant \(A_{\alpha_j}\), then
\(M_{\gamma_{\alpha_j}}(\boldsymbol{t})= \exp{\{2\pi i
A_{\alpha_j}\}}\), where \(M_{\gamma_{\alpha_j}}\) is the
monodromy matrix along the loop \(\gamma_{\alpha_j}\) generated by
the path \(\alpha_j\). Therefore, \(M_{\gamma_{\alpha_j}}\) does
not depend on \(\boldsymbol{t}\). \hfill Q.E.D.

\begin{remark}%
\label{RCWC}%
If a deformation \textup{(\ref{FDES})} is isomonodromic, but
 some of the matrices \(Q_j(\boldsymbol{t})\) are not non-resonant, then this
 deformation may be not isoprincipal. As we shall see later, there is a rich
 class of deformations of Fuchsian equations all fundamental solutions
  \(Y(x,\boldsymbol{t})\) of which are rational functions of \(x\).
  Thus, such deformations are isomonodromic. (Their monodromies are trivial).
  However, most of these deformations are not isoprincipal.
\end{remark}%

\noindent From Lemmas \ref{ImIpr} and \ref{IpIso} it follows that
in the non-resonant case the classes of isomonodromic deformations
and of isoprincipal (in the ``narrow sense'') deformations
coincide. This is a useful observation. On the one hand, it is the
property of "isoprincipalness" that implies the Schlesinger system
for the matrix functions \(Q_p(\boldsymbol{t})\) appearing in the
deformation. (We prove this implication below). On the other hand,
there are methods, which are based on solving the so called
Riemann-Hilbert problem, that allow to construct deformations of
Fuchsian equations with the prescribed monodromy, but not with the
prescribed principal factors of the solutions \(Y_{\alpha_j}\). In
particular, such methods allow (under non-resonance condition and
certain other restrictions) to construct directly isomonodromic
(but not isoprincipal) deformations.   The proof of the fact that
the resulting \(Q_j(\boldsymbol{t})\) satisfy the Schlesinger
system uses essentially the equivalency of these two classes of
deformations in the non-resonant case.

However, the above given ``narrow'' definition of the
isoprincipalness, which is quite natural in the non-resonant case,
seems to be a little bit artificial and too restrictive in the
general case. Indeed, in the general case the solutions
\(Y_{\alpha_j}(\boldsymbol{t})\) of the equation (\ref{FDES}) are
of the form (\ref{SRQIT}), where exponents
\(Z_{\alpha_j}(\boldsymbol{t})\) and
\(A_{\alpha_j}(\boldsymbol{t})\) depend on \(\boldsymbol{t}\) in
general and posses some additional properties. In particular, the
matrix \(A_{\alpha_j}(\boldsymbol{t})\) is non-resonant, the
matrix \(Z_{\alpha_j}(\boldsymbol{t})\) is diagonalizable, with
integer eigenvalues. (These exponents may be non-holomorphic, and
even discontinuous functions of \(\boldsymbol{t}\). See Remark
\ref{GRC}). If these exponents are constant then we have good
reason to call such a deformation isoprincipal. However, we have
no need to go into detail. All that is important for us is that
the principal factor \(\displaystyle
(x-t_j)^{Z_{\alpha_j}}\cdot(x-t_j)^{A_{\alpha_j}}\) depends on the
difference \(x-t_j\) only if the exponents \(Z_{\alpha_j}\) and
\(A_{\alpha_j}\) do not depend on \(\boldsymbol{t}\).

\begin{definition}%
\label{GIPD}%
The deformation \textup{(\ref{FDES})} is said to be
\textsf{isoprincipal with respect to a distinguished point
\(x_0\)} (in the wide sense) at a pole loci \(\boldsymbol{t}^0=
(t_1^0, \dots , t_n^0)\),
\((\boldsymbol{t}^0\in\boldsymbol{\mathcal{D}}\)), if for every
\(j=1, \dots , n\) and for some paths
\(\alpha_1, \dots , \alpha_n\) leading from \(x_0\) to
neighborhoods of the points \(t_1^0, \dots , t_n^0\)
respectively, the solution \(Y_{\alpha_j}(x,\boldsymbol{t})\)
normalized at \(x_0\) by \textup{(\ref{ICondT})} is representable
in the form
\begin{equation}%
\label{EIpD}%
Y_{\alpha_j}(x,\boldsymbol{t})=H_{\alpha_j}(x,\boldsymbol{t})E_{\alpha_j}(x-t_j),
\quad x\in D_{t_j, \rho}\setminus\{t_j\},
\boldsymbol{t}\in\mathsf{D}(\boldsymbol{t}^0,\delta) ,
\end{equation}%
where the matrix functions \(H_{\alpha_j}(x,\boldsymbol{t})\),
\(\big(H_{\alpha_j}(x,\boldsymbol{t})\big)^{-1}\)
 are holomorphic in \(x,\boldsymbol{t}\) for \(x\in D_{t_j, \rho},
\boldsymbol{t}\in\mathsf{D}(\boldsymbol{t}^0,\delta)\), (\(\rho\) and
\(\delta\) are some positive numbers), and
\(E_{\alpha_j}(\zeta), j=1, \dots , n,\) are matrix functions possessing
the properties:
\begin{enumerate}
\item
\(E_{\alpha_j}(\zeta)\) are holomorphic on the Riemann surface of
logarithm\footnote{ The Riemann surface of logarithm is the
universal covering
\(\text{cov}\big(\mathbb{C}\setminus\{0\}\big)\). Considered as
functions defined in \(\mathbb{C}\setminus\{0\}\),
\(E_{\alpha_j}(\zeta)\) are multivalued in general.} and take
invertible values there;
\item
\(E_{\alpha_j}(\zeta)\) do not depend on \(\boldsymbol{t}\) and
\(x\).
\end{enumerate}
\end{definition}%

\begin{remark}%
\label{CGID}%
Of course, to formulate Definition \ref{GIPD} more accurately, we
have to describe the geometric configuration more accurately. For
every \(j\), we have to choose the point \(x_j\) according to
(\ref{GeomConfig}), assuming that \(\delta\) is small enough, in
particular, (\ref{CFDe}) is satisfied. (The value
\(\rho_j(\boldsymbol{t}^0)\) is defined according to
(\ref{DefRho})). By such a choice of \(x_j\), the condition
(\ref{DesGeom}) is satisfied. Then we have to choose the arguments
\(\arg (x_j-t_j^0)\) somehow, (\ref{ChoiArg0}), and then \(\arg
(x_j-t_j)\) for \(t_j\in D_{t_j^0,\delta}\) according to
(\ref{ChoiArg}). A point of the universal covering
\(\cov(D_{t_j,\rho_j(\boldsymbol{t}^0)}\setminus \{t_j\}, x_j)\)
is a pair \((x,\mu)\), where \(x\) is a point from the punctured
disk \(D_{t_j,\rho_j(\boldsymbol{t}^0)})\setminus \{t_j\}\), and
\(\mu\) is a path which starts at the point \(x_j\) and ends at
the point \(x\), \(\mu\in
D_{t_j,\rho_j(\boldsymbol{t}^0)})\setminus \{t_j\}\). The path
\(\mu\), together with the choice of \(\arg(x_j-t_j)\), determines
the value of \(\arg (x-t_j)\), and hence, the value of
\(E_{\alpha_j}(x-t_j)\). On the other hand, the pair
\((x,\mu)\in\cov(D_{t_j,\rho_j(\boldsymbol{t}^0)}\setminus
\{t_j\}, x_j)\) determines the value of the analytic continuation
\(Y_{\alpha_j}\) along the path \(\alpha_j\) at this point
\((x,\mu)\). By definition, the value
\(Y_{\alpha_j}((x,\mu),\boldsymbol{t})\) is the result of the
analytic continuation of the solution normalized at \(x_0\) by
(\ref{IniCoT}) along the composed path \(\mu\cdot\alpha_j\) from a
neighborhood of \(x_0\) to the point \(x\). Since both values,
\(Y_{\alpha_j}(x,\boldsymbol{t})\) and \(E_{\alpha_j}(x-t_j)\) are
determined for every choice of the path \(\mu\) (these values may
depend on \(\mu\)), the value
\(H_{\alpha_j}(x,\boldsymbol{t})\stackrel{\textup{\tiny def}}{=}
Y_{\alpha_j}(x,\boldsymbol{t}) \big(E_{\alpha_j}(x-t_j)\big)^{-1}
\) is determined for \(x\in
D_{t_j,\rho_j(\boldsymbol{t}^0)}\setminus \{t_j\}\), and it may
depend on \(\mu\). The definition \ref{GIPD} means firstly, that
the value \(H_{\alpha_j}(x,\boldsymbol{t})\) is determined only by
the point \(x\in D_{t_j,\rho_j(\boldsymbol{t}^0)}\setminus
\{t_j\}\) itself, but does not depend on the path \(\mu\) leading
to this point, thus the function
\(H_{\alpha_j}(x,\boldsymbol{t})\) as a function of \(x\) is
single-valued in the punctured disk
\(D_{t_j,\rho_j(\boldsymbol{t}^0)}\setminus \{t_j\}\). Secondly,
this definition means that the singularity \(x=t_j\) of the
function \(H_{\alpha_j}(x,\boldsymbol{t})\) is removable, that is
there exists the limit \(\lim_{x\to t_j, x\not=
t_j}H_{\alpha_j}(x,\boldsymbol{t})\stackrel{\textup{\tiny
def}}{=}H_{\alpha_j}(t_j,\boldsymbol{t})\). Thirdly, the
definition \ref{GIPD} means that the value
\(H_{\alpha_j}(t_j,\boldsymbol{t})\) is an invertible matrix. (For
\(x\not=t_j\), the values \(H_{\alpha_j}(t_j,\boldsymbol{t})\) are
invertible, because both the values
\(Y_{\alpha_j}(x,\boldsymbol{t})\) and \(E_{\alpha_j}(x-t_j)\) are
invertible). So, the accurate definition of the isoprincipal
deformation is rather bulky.
\end{remark}%

\noindent Below (Lemma \ref{IDWD}) we prove that the notion of
isoprincipal deformation is well defined, that is the property of
the deformation to be isoprincipal at some point does not depend
on the choice of the paths \(\alpha_j\).

\begin{theorem}%
\label{GenMon}%
If \textup{(\ref{FDES})} is a deformation of Fuchsian equation
which is isoprincipal with respect to a distinguished point in the
wide sense (i.e, in the sense of Definition \ref{GIPD}) at a pole
loci \(\boldsymbol{t}^0= (t_1^0, \dots , t_n^0)\),
\((\boldsymbol{t}^0\in\boldsymbol{\mathcal{D}}\)), then the
matrices \(E_{\alpha_j}(\zeta), j=1, \dots , n\), representing
the principle factors, must satisfy the monodromic relations of
the form
\begin{equation}%
\label{MonRel}%
E_{\alpha_j}(\zeta e^{2\pi i})=E_{\alpha_j}(\zeta)M_{\alpha_j},\quad
\forall \zeta\not=0 ,
\end{equation}%
where \(M_{\alpha_j}\) are constant (not depending on \(\zeta\)) invertible
matrices. For the monodromy matrix \(M_{\gamma_{\alpha_j}}(\boldsymbol{t})\)
of the equation \textup{(\ref{FDES})} along the loop \(\gamma_{\alpha_j}\)
generated by the paths \(\alpha_j\), (see Figure \ref{BGHj}), the equality
\begin{equation}%
\label{EMMo}%
M_{\gamma_{\alpha_j}}(\boldsymbol{t})=M_{\alpha_j}, \quad
\boldsymbol{t}\in\mathsf{D}(\boldsymbol{t}^0,\delta),\ \ j=1, \dots , n .
\end{equation}%
holds. In particular, this deformation is isomonodromic with
respect to \(x_0\).
\end{theorem}%

\noindent\textsf{PROOF.} %
For fixed
\(\boldsymbol{t}\in\mathsf{D}(\boldsymbol{t}^0,\delta)\), and for
\(x\in D_{t,\rho}\setminus\{t_j\}\), consider two matrix functions
of \(x\): \(Y_{\alpha_j}(x,\boldsymbol{t})\) and
\(Y_{\alpha_j}(t_j+(x-t_j)e^{2\pi i},\boldsymbol{t})\). Here \(x\)
and \(t_j+(x-t_j)e^{2\pi i}\) are considered as points running
over the universal covering \(\cov
(D_{t_j\!,\rho}\setminus\{t_j\}, x_j)\) of the punctured disk
\(D_{t_j\!,\rho}\setminus\{t_j\}\) with the distinguished point
\(x_j\). These two points of the universal covering have the same
projections in the complex plane. Therefore, the matrices
\(Y_{\alpha_j}(x,\boldsymbol{t})\) and
\(Y_{\alpha_j}(t_j+(x-t_j)e^{2\pi i},\boldsymbol{t})\) can be
considered as the values of two different solutions of the same
linear differential equation at the same point. From \(\det
Y(x,\boldsymbol{t})_{|x=x_j}\not=0\) (see (\ref{ICondT})) it
follows that \(\det Y_{\alpha_j}(x,\boldsymbol{t})\not=0\), \(\det
Y_{\alpha_j}(t_j+(x-t_j)e^{2\pi i},\boldsymbol{t})\not=0\). Hence,
these solutions are ``proportional'' as functions of \(x\):
\begin{equation}%
\label{Prop}%
Y_{\alpha_j}(t_j+(x-t_j)e^{2\pi
i},\boldsymbol{t})=Y_{\alpha_j}(x,\boldsymbol{t})M_{\alpha_j}(\boldsymbol{t}) ,
\quad x\in
D_{t_j,\rho}\setminus\{t_j\} ,\boldsymbol{t}\in\mathsf{D}(\boldsymbol{t}^0,\rho) ,
\end{equation}%
where \(M_{\alpha_j}(\boldsymbol{t})\) are constant (with respect to \(x\))
non-degenerate  matrices, which may
depend on \(\boldsymbol{t}\). Since the factor \(H_{\alpha_j}\) is univalued
(\(H_{\alpha_j}(t_j+(x-t_j)e^{2\pi
i},\boldsymbol{t})=H_{\alpha_j}(x,\boldsymbol{t})\)), from (\ref{EIpD}) and
(\ref{Prop}) it follows that %
\begin{equation}%
\label{QMR}%
 E_{\alpha_j}((x-t_j)e^{2\pi
i})=E_{\alpha_j}(x-t_j)M_{\alpha_j}(\boldsymbol{t})
\end{equation}%
Since \(x\) is an arbitrary point of \(D_{t_j,\rho}\setminus \{t_j\}\),
\(x-t_j\) is an arbitrary point of \(\mathbb{C}\setminus\{0\}\). Since the
function \(E_{\alpha_j}\) does not depend on \(\boldsymbol{t}\), the matrix
\(M_{\alpha_j}\) does not depend on \(\boldsymbol{t}\) as well. Moreover, the
monodromic relation (\ref{MonRel}) holds at lest for
\(\zeta: 0<|\zeta|<\rho\). However, the function \(E_{\alpha_j}\) is
holomorphic in the whole Riemann surface of logarithm. Therefore, the relation
(\ref{MonRel}) holds for all \(\zeta: |\zeta|>0\).

The expression for the monodromy matrix
\(M_{\gamma_{\alpha_j}}(\boldsymbol{t})\) along the loop \(\gamma_{\alpha_j}\)
was already obtained (see (\ref{EFM})):
\begin{equation}%
\label{ExpMon}%
 M_{\gamma_{\alpha_j}}(\boldsymbol{t})=Y_{x_j,\alpha_j}(x_j,\boldsymbol{t})^{-1}\cdot
Y_{x_j,\alpha_j}(t_j+(x_j-t_j)e^{2\pi i},\boldsymbol{t}) .
\end{equation}%
The equality (\ref{EMMo}) is the consequence of (\ref{Prop}) and
(\ref{ExpMon}). \hfill Q.E.D.\\[1ex]

\noindent Let for a certain \(p\), \(\alpha_p\) be a path leading
from the point \(x_0\) to a neighborhood of the point \(t_p\).
 Let \(\alpha_p^{\prime}\) be another path   leading from
\(x_0\) to the same neighborhood. Let \(Y_{\alpha_p}\) and
\(Y_{\alpha_p^{\prime}}\) be the analytic continuations of the normalized at
\(x_0\) solution of (\ref{FDES}) along the paths \(\alpha_p\) and
\(\alpha_p^{\prime}\) respectively. Then the solutions \(Y_{\alpha_p}\) and
\(Y_{\alpha_p^{\prime}}\) are related by the relation %
\(Y_{\alpha_p^{\prime}}(x,\boldsymbol{t})= Y_{\alpha_p}(x,\boldsymbol{t})\cdot
M_{\gamma}(\boldsymbol{t}), \ x\in D_{t_p,\rho}\setminus\{t_j\}\), %
 where \(M_{\gamma}(\boldsymbol{t})\) is the monodromy matrix along the loop
\(\gamma\), constructed from the paths \(\alpha_p\) and
\(\alpha_p^{\prime}\). (See footnote \ref{DiPa}). If the deformation
\textup{(\ref{FDES})} is isoprincipal in the sense of the definition
\ref{GIPD} by a certain choice of the paths \(\alpha_j, j=1, \dots , n\),
then, according to Theorem \ref{GenMon}, this deformation is isomonodromic. In
particular, the monodromy matrix \(M_{\gamma}(\boldsymbol{t})\) in (\ref{RDS})
does not depend on \(\boldsymbol{t}\):
\(M_{\gamma}(\boldsymbol{t})=M_\gamma\), thus
\begin{equation}%
\label{RDS}%
Y_{\alpha_p^{\prime}}(x,\boldsymbol{t})= Y_{\alpha_p}(x,\boldsymbol{t})\cdot
M_{\gamma}, \quad x\in D_{t_p,\rho}\setminus\{t_j\},
\end{equation}%
 Therefore, the  solution
\(Y_{\alpha_p^{\prime}}\) is representable in the form
\begin{equation}%
\label{ReC}%
Y_{\alpha_p^{\prime}}(x,\boldsymbol{t})=H_{\alpha_p^{\prime}}(x,\boldsymbol{t})
\cdot E_{\alpha_p^{\prime}}(x-t_j),\quad x\in D_{t_j,\rho}\setminus \{t_j\},
\boldsymbol{t}\in\mathsf{D}(t_j,\rho) ,
\end{equation}%
with
\(H_{\alpha_p^{\prime}}(x,\boldsymbol{t})=H_{\alpha_p}(x,\boldsymbol{t})M_{\gamma}\),
\(E_{\alpha_p^{\prime}}(\zeta)=M_{\gamma}^{-1}E_{\alpha_p}(\zeta)M_{\gamma}\).
Since the matrices \(E_{\alpha_p}(\zeta)\) and \(M_{\gamma}\) do not depend on
\(\boldsymbol{t}\), the matrix \(E_{\alpha_p^{\prime}}(\zeta)\) does not
depend on \(\boldsymbol{t}\) as well. So, we prove the following

\begin{lemma}%
\label{IDWD}%
The notion of the isoprincipal deformation at a point is well
defined. If the deformation \textup{(\ref{FDES})} is isoprincipal
at the point \(\boldsymbol{t}^0\) for some choice of the paths
\(\alpha_j,\ j=1, \dots , n\) leading from \(x_0\) to
neighborhoods of the singularities \(t_j^0\), then this
deformation is isoprincipal at this point for any other choice of
such paths.
\end{lemma}

\begin{definition}%
\label{GGIPD}%
The deformation \textup{(\ref{FDES})}, where the coefficients
\(Q_p(\boldsymbol{t}), p=1, \dots , n ,\) are holomorphic for
\(\boldsymbol{t}\in\boldsymbol{\mathcal{D}}\),
\(\boldsymbol{\mathcal{D}}\) is an open set in
\(\mathbb{C}^n_{\ast}\), is said to be \textsf{isoprincipal} (in
the wide sense) \textsf{in} \(\boldsymbol{\mathcal{D}}\) if it is
isoprincipal at every point
\(\boldsymbol{t}\in\boldsymbol{\mathcal{D}}\).
\end{definition}

\begin{theorem}%
\label{LIPGIP} If the deformation \textup{(\ref{FDES})} is
isoprincipal with respect  to \(x_0\) at some
\(\boldsymbol{t}^0\),
\(\boldsymbol{t}^0\in\boldsymbol{\mathcal{D}}\), and the open set
\(\boldsymbol{\mathcal{D}}\) is connected, then this deformation
is isoprincipal with respect  to \(x_0\) in
\(\boldsymbol{\mathcal{D}}\) (i.e., at every
\(\boldsymbol{t}\in\boldsymbol{\mathcal{D}}\)).
\end{theorem}%

\noindent We omit proof of this theorem.

\section{ ISOPRINCIPALNESS OF A DEFORMATION IMPLIES  THE SCHLESINGER SYSTEM FOR
ITS COEFFICIENTS} \noindent
\begin{lemma}%
\label{IICC}%
Let the deformation\,%
\footnote{%
The condition \( \sum\limits_{1\leq p\leq n}Q_p
({\boldsymbol{t}})\equiv 0, \ \boldsymbol{t}\in\mathcal{D}, \) is
assumed.
}%
 \textup{(\ref{FDES})} be isoprincipal with
respect to a distinguished point \(x_0\) (in the wide sense) at
some pole loci \(\boldsymbol{t}^0\),
\(\boldsymbol{t}^0\in\boldsymbol{\mathcal{D}}\), and let \(Y(x,
\boldsymbol{t}),\) \(
x\in\text{cov} (\overline{\mathbb{C}}\setminus\{t_1, \dots , t_n\}), %
\boldsymbol{t}\in\boldsymbol{\mathcal{D}},\) be the solution of
the differential equation (\ref{FDES}) which is normalized by the
condition \textup{(\ref{ICondT})} at the distinguished point
\(x_0\).

Then the equalities
\begin{equation}%
\label{DESS}%
\frac{\partial Y(x, \boldsymbol{t})}{\partial
t_k}=Q_k(\boldsymbol{t})\left(\frac{1}{x_0-t_k}-\frac{1}{x-t_k}\right)
Y(x, \boldsymbol{t}) ,\quad
k=1, \dots , n
\end{equation}%
 are satisfied. In particular, if \(x_0=\infty\), then the
 equations \textup{(\ref{ODS2})} are satisfied.
\end{lemma}%

\noindent%
\textsf{PROOF.} The key to the proof are the (trivial, but
nevertheless very useful) equalities
\begin{equation}%
\label{RaVa}%
\frac{d }{dx} E_{\alpha_j}(x-t_j)=- \frac{\partial  }{\partial
t_j} E_{\alpha_j}(x-t_j) ,\quad j=1, \dots , n ,
\end{equation}%
which express the fact that the functions \(E_{\alpha_j}(x-t_j)\)
depend only on the difference of the arguments \(x\) and \(t_j\).
 Expressing the logarithmic derivative
\(\dfrac{dY}{dx} Y^{-1}\) from the differential equation
(\ref{FDES}): \(\dfrac{dY}{dx} Y^{-1}=\displaystyle
\sum\limits_p \frac{Q_p(\boldsymbol{t})}{x-t_p}\) , we derive
its behavior in a neighborhood of the singular point \(t_j\):
\begin{multline}%
\label{AB1}%
\frac{dY}{dx} (x, \boldsymbol{t}) Y^{-1}(x, \boldsymbol{t})=\\
=\frac{Q_j(\boldsymbol{t})}{x-t_j}
+ (\text{a function of \(x\), holomorphic with respect  to \(x\)}),\\
\ \ x\in D_{t_j, \rho} , \ \ j=1, \dots n .
\end{multline}%
(\(D_{t_j, \rho}=\{z\in\mathbb{C}: |z-t_j|<\rho \},\ \rho\) is
a positive number).
 We restrict our consideration to
\(\boldsymbol{t}\in\mathsf{D}(\boldsymbol{t}^0,\delta)\), where
\(\delta>0\) is small enough, so that the representations
(\ref{EIpD}) of the solution \(Y(x,\boldsymbol{t})\) hold in the
appropriate neighborhoods \(D_{t_j, \rho}\setminus\{t_j\}\) of
the singular points \(t_j\).
 Using the representation
(\ref{EIpD}), we see that
\begin{multline}%
\label{AB2}%
\frac{dY_{\alpha_j}}{dx} (x, \boldsymbol{t}) Y_{\alpha_j}^{-1}(x, \boldsymbol{t})=
H_{\alpha_j}(x, \boldsymbol{t}) \left(\frac{d }{dx}
E_{\alpha_j}(x-t_j)E_{\alpha_j}^{-1}(x-t_j)\right)
H_{\alpha_j}^{-1}(x, \boldsymbol{t})\\
+ \left(\frac{d }{dx} H_{\alpha_j}(x, \boldsymbol{t})) \right)
 H_{\alpha_j}^{-1}(x, \boldsymbol{t}) .
\end{multline}%
Since the function \(H_{\alpha_j}(x, \boldsymbol{t})\) is
holomorphic and invertible in an entire (non-punctured)
neighborhood of the singular point \(t_j\), comparing (\ref{AB1})
and (\ref{AB2}), we obtain
\begin{multline}
\label{AB3}%
H_{\alpha_j}(x, \boldsymbol{t}) \left(\frac{d }{dx}
E_{\alpha_j}(x-t_j)\big(E_{\alpha_j}(x-t_j)\big)^{-1}\right)
H_{\alpha_j}(x, \boldsymbol{t})^{-1}=\\
=\frac{Q_j}{x-t_j}  +
(\text{a function of \(x\), holomorphic with respect  to
\(x\)}),\\ x\in D_{t_j, \rho} , \ \ j=1, \dots n .
\end{multline}
Let us investigate the logarithmic derivatives \(\dfrac{\partial
Y}{\partial t_k} (x, \boldsymbol{t}) Y^{-1}(x, \boldsymbol{t}) ,\)
as  functions of \(x\). First of all we remark that for fixed
\(\boldsymbol{t}\), these function are single valued and
holomorphic with respect to \(x\) in the domain
\(\overline{\mathbb{C}}\setminus\{t_1, \dots , t_n\}\). Indeed,
let \(\mathcal{O}\) be a simply connected open set in
\(\overline{\mathbb{C}}\) separated from the points \(t_1^0, \dots
, t_n^0\): \(\text{dist} (\mathcal{O}, \{t_1^0, \dots ,
t_n^0\})>0.\) Let \(Y_1(x, \boldsymbol{t}), Y_2(x,
\boldsymbol{t})\) be two branches of the solution \(Y(x,
\boldsymbol{t})\) in \(\mathcal{O}\), that is the functions
\(Y_1(x, \boldsymbol{t}), Y_2(x, \boldsymbol{t})\) are obtained by
the analytic continuation of the solution of (\ref{FDES}) along
two different paths \(\gamma_1, \gamma_2\) leading from \(x_0\) to
\(\mathcal{O}\). These branches \(Y_1(x, \boldsymbol{t}), Y_2(x,
\boldsymbol{t})\) must satisfy the relation \(Y_2(x,
\boldsymbol{t})=Y_1(x,
\boldsymbol{t})M_{\gamma}(\boldsymbol{t})\), where
\(M_{\gamma}(\boldsymbol{t})\) is the monodromy matrix
corresponding to the loop \(\gamma\) constructing from the paths
\(\gamma_1\) and \(\gamma_2\). In general, this monodromy matrix
\(M_{\gamma}(\boldsymbol{t})\) may depend on \(\boldsymbol{t}\).
However, according to Theorem \ref{GenMon}, the considered
deformation is isomonodromic. In particular, the matrix
\(M_{\gamma}(\boldsymbol{t})\) does not depend on
\(\boldsymbol{t}\), that is \(Y_2(x, \boldsymbol{t})=Y_1(x,
\boldsymbol{t})M_{\gamma}\), where \(M_\gamma\) does not depend on
\(\boldsymbol{t}\) and is invertible. (\(M_\gamma\) may depend on
the choice of the branches \(Y_1\) and \(Y_2\).) Differentiating
the last equality with respect to \(t_k\), we obtain
\(\dfrac{\partial Y_2}{\partial t_k}=\dfrac{\partial Y_1}{\partial
t_k} M_{\gamma.}\) From the last two equalities it follows that
\(\dfrac{\partial Y_1}{\partial t_k} Y_1^{-1}=\dfrac{\partial
Y_2}{\partial t_k} Y_2^{-1}\), that is the logarithmic derivative
\(\dfrac{\partial Y}{\partial t_k} Y^{-1}\) is a single-valued
holomorphic function of \(x\) in
\(\overline{\mathbb{C}}\setminus\{t_1, \dots , t_n\}\). Using the
representation (\ref{EIpD}) in a neighborhood of the singular
point \(t_j\), we see that for \(j\not=k\)
\begin{equation*}%
\label{CuSP}%
\frac{\partial Y}{\partial
t_k}((x, \boldsymbol{t})) Y^{-1}((x, \boldsymbol{t}))=\frac{\partial
H_{\alpha_j}(x, \boldsymbol{t})}{\partial
t_k}  H_{\alpha_j}^{-1}(x, \boldsymbol{t}) , \quad x\in
D_{t_j,\rho}.
\end{equation*}%
(For \(j\not=k\), the second factor \(E_{\alpha_j}(x-t_j)\) in
(\ref{EIpD}) does not depend on \(t_k\).)
 Since the function \(H_{\alpha_j}(x, \boldsymbol{t})\) is
holomorphic with respect  to \((x, \boldsymbol{t})\) for \(x\in
D_{t_j, \rho} ,
\boldsymbol{t}\in\mathsf{D}(\boldsymbol{t}^0,\delta)\),
so is the function \( \frac{\partial Y}{\partial
t_k}((x, \boldsymbol{t})) Y^{-1}((x, \boldsymbol{t}))\)
for \( j\not=k\).

For \(j=k\), we have to take into account the factor
\(E_{\alpha_j}(x-t_j)\):
\begin{multline*}%
\frac{\partial Y}{\partial
t_k}((x, \boldsymbol{t})) Y^{-1}((x, \boldsymbol{t}))=\\
=H_{\alpha_k}(x, \boldsymbol{t})\left(\frac{\partial }{\partial
t_k
} E_{\alpha_k}(x-t_k)  E_{ \alpha_k}^{-1}(x-t_k)\right)H_{\alpha_k}^{-1
}(x, \boldsymbol{t})\\ +\frac{\partial
H_{\alpha_k}(x, \boldsymbol{t})}{\partial
t_k}  H_{\alpha_k}^{-1}(x, \boldsymbol{t}) , \quad x\in
D_{t_k,\rho} .
\end{multline*}
The function \(\dfrac{\partial
H_{\alpha_k}(x, \boldsymbol{t})}{\partial
t_k}  H_{\alpha_k}^{-1}(x, \boldsymbol{t})\) is holomorphic for
\(x\in D_{t_j, \rho}\). Taking into account (\ref{RaVa}) and
(\ref{AB3}), we come to the conclusion:
\begin{multline}%
\label{CuSP3}%
\frac{\partial Y}{\partial
t_k}((x, \boldsymbol{t})) Y^{-1}((x, \boldsymbol{t}))=
-\frac{Q_k}{x-t_k} +\\
+(\text{a function holomorphic with respect  to \(x\)}),\ \ x\in
D_{t_k, \rho}.
\end{multline}%
Thus, the function \(\dfrac{\partial Y}{\partial
t_k}((x, \boldsymbol{t})) Y^{-1}((x, \boldsymbol{t}))\),
considered as a function of the variable \(x\), is holomorphic in
the entire extended complex plane \(\overline{\mathbb{C}}\),
except the point \(t_k\). Near \(t_k\) this function behaves as it
described in (\ref{CuSP3}). Therefore,
\begin{equation*}%
\label{CuSP4} \dfrac{\partial Y}{\partial
t_k}((x, \boldsymbol{t})) Y^{-1}((x, \boldsymbol{t}))=-\frac{Q_k}{x-t_k}+C_k(\boldsymbol{t}) ,
\end{equation*}%
where \(C_k(\boldsymbol{t})\) does not depend on \(x\). To
determine \(C_k(\boldsymbol{t})\), we use the normalizing
condition (\ref{ICondT}). Since \(Y(x_0, \boldsymbol{t})\equiv
I\), \(\dfrac{\partial Y(x_0, \boldsymbol{t})}{\partial t_k}=0\),
i.e. \(-\dfrac{Q_k(\boldsymbol{t})}{x-t_k}+C_k(\boldsymbol{t})=0\)
at \(x=x_0\). Thus, (\ref{DESS}) holds. %
\hfill Q.E.D.

\begin{theorem}%
\label{IPISchS}%
Let a \(k\times k\) matrix functions \(Q_j(\boldsymbol{t})\),
\(\boldsymbol{t}=(t_1, \dots , t_n) , j=1, \dots , n\),
are holomorphic for \(\boldsymbol{t}\in\boldsymbol{\mathcal{D}}\),
where \(\boldsymbol{\mathcal{D}}\subset\mathbb{C}^n_{\ast}\) is an
open connected set, and satisfying the condition
\begin{equation*}
\label{VaCo}%
\sum\limits_{1\leq j\leq n}Q_j(\boldsymbol{t})=0 ,\quad
\boldsymbol{t}\in\boldsymbol{\mathcal{D}} .
\end{equation*}
Assume that the deformation
\begin{equation*}
\label{IPrD}%
\frac{dY}{dx}=\bigg(\sum\limits_{1\leq j\leq
n}\frac{Q_j(\boldsymbol{t})}{x-t_j}\bigg) Y
\end{equation*}
is isoprincipal (in the wide sense) with respect to the
distinguished point \(x_0=\infty\) at least for one pole loci
\(\boldsymbol{t}^0\in\boldsymbol{\mathcal{D}}\).

Then the matrix functions
\(Q_j(\boldsymbol{t}), j=1, \dots , n\), satisfy the
Schlesinger system (\ref{Sch}) for all
\(\boldsymbol{t}\in\boldsymbol{\mathcal{D}}\).
\end{theorem}%

\noindent%
\textsf{PROOF}. Let \(Y(x, \boldsymbol{t})\) be the solution of
the differential equation (\ref{FDES}) satisfying the normalizing
condition \(Y(x, \boldsymbol{t})_{x=\infty}=I .\) According to
Lemma \ref{IICC}, the equations (\ref{ODS2}) are satisfied for
\(x\in\overline{\mathbb{C}}\setminus\{t_1, \dots , t_n\} ,
\boldsymbol{t}\in\mathsf{D}(\boldsymbol{t}^0,\delta)\). The
equation (\ref{ODS1}) is the same as the equation (\ref{FDES}).
According to Theorem \ref{SEC}, the matrix functions
\(Q_j(\boldsymbol{t})\),
\(\boldsymbol{t}=(t_1, \dots , t_n) , j=1, \dots , n\)
satisfy the Schlesinger system (\ref{Sch}) for
\(\boldsymbol{t}\in\mathsf{D}(\boldsymbol{t}^0,\delta)\). Since
matrix functions in both left and right hand sides of (\ref{Sch})
are holomorphic in \(\boldsymbol{\mathcal{D}}\) and coincide in
\(\mathsf{D}(\boldsymbol{t}^0,\delta)\), they coincide everywhere
in \(\boldsymbol{\mathcal{D}}\). \hfill Q.E.D.%

\section{CONSTRUCTION OF THE ISOPRINCIPAL DEFORMATIONS
 AS A TOOL FOR SOLVING THE SCHLESINGER SYSTEM.}

Theorem \ref{IPISchS} opens a way for constructing solutions of
the Schlesinger system. The \textsf{Cauchy problem} for this
system can be formulated as follows.

\noindent\textsl{Given \(\boldsymbol{t}^0\in\mathbb{C}^n_{\ast}\)
and \(k\times k\) matrices \(Q_1^0, \dots , Q_n^0\), one needs
to find matrix functions
\(Q_1(\boldsymbol{t}), \dots , Q_n(\boldsymbol{t})\) of the
variable \(\boldsymbol{t}\) satisfying the Schlesinger system
(\ref{Sch}) and the initial condition}
\begin{equation}%
\label{ICSchS}%
Q_j(\boldsymbol{t^0})=Q_j^0,\quad j=1, \dots ,n .
\end{equation}%

\noindent If we succeed in constructing a family of Fuchsian
equations
\begin{equation}%
\label{FDES8}%
\frac{dY}{dx}=\bigg(\sum\limits_{1\leq j\leq
n}\frac{Q_j(\boldsymbol{t})}{x-t_j}\bigg) Y ,
\end{equation}%
 enumerated by the pole loci
\(\boldsymbol{t}=(t_1,\dots,t_n)\), with holomorphic coefficients
\(Q_1(\boldsymbol{t}),\) \(\dots, Q_n(\boldsymbol{t})\) (a
holomorphic deformation) such that \textsl{this holomorphic
deformation is isoprincipal with respect to \(x_0=\infty\)}, then,
according to Theorem \ref{IPISchS}, these coefficients will
satisfy the Schlesinger system. In order for the initial
conditions \((\ref{ICSchS})\) to be satisfied, the family must
contain the equation
\begin{equation}%
\label{FDES81}%
\frac{dY}{dx}=\bigg(\sum\limits_{1\leq j\leq
n}\frac{Q_j^0}{x-t_j^0}\bigg) Y ,
\end{equation}%
which corresponds to the parameter value
\(\boldsymbol{t}=\boldsymbol{t}^0\).

If all the matrices \(Q_1^0, \dots , Q_n^0\) are non-resonant
(and hence, according to Theorems \ref{GenMon} and \ref{IMIS}, all
the matrices
\(Q_1(\boldsymbol{t}), \dots , Q_n(\boldsymbol{t})\) are
non-resonant as well),  this deformation is isoprincipal if, and
only if, it is isomonodromic (Theorem \ref{GenMon} and Lemma
\ref{ImIpr}). Therefore, if all the matrices
\(Q_1^0, \dots , Q_n^0\) are non-resonant, the Cauchy problem
with initial conditions (\ref{ICSchS}) is reduced to the problem
of constructing the \textsl{isomonodromic} deformation
(\ref{FDES8}) which in particular contains the equation
(\ref{FDES81}). An approach to constructing such  isomonodromic
deformation lies in using \textsl{the Riemann-Hilbert problem.}

The Riemann-Hilbert problem can be formulated as follows.

\noindent\textsl{Given
\(\boldsymbol{t}=(t_1,\dots , t_n)\in\mathbb{C}^n_{\ast}\),
loops
\(\gamma_1(\boldsymbol{t}),\dots , \gamma_n(\boldsymbol{t})\)
with the distinguished point \(x_0\) that form the geometric
configuration as it is plotted in Figure \ref{BGH}, and matrices
\(M_1, \dots, M_n\) satisfying the
condition\footnote{\(\{\pi_1,\dots , \pi_n\}\) is a permutation
of the indices \(\{1, \dots , n\}\) which is determined by the
geometric configuration of the paths
\(\alpha_1, \dots , \alpha_n.\) For the configuration plotted
in Figure \ref{BGH} the permutation is trivial: \(\pi_k=k .\)
}\(M_{\pi(1)}\cdot \cdots \cdot M_{\pi(n)}=I\). One needs to
construct the Fuchsian differential equation (\ref{FDES8}) with
poles \(t_1, \dots , t_n\) whose monodromy matrices take the
prescribed values}
\begin{equation}%
\label{PMoMa}%
M_{\gamma_j(\boldsymbol{t})}=M_j,\quad j=1, \dots , n .
\end{equation}%

\noindent An approach for solving the Riemann-Hilbert problem was
developed by J. Plemelj in 1908, and, a little later but
independently, by G.D. Birkhoff (see \cite{Plem1}, \cite{Plem2},
\cite{Birk}). The first step of this approach requires to solve a
factorization problem for a matrix function on a closed contour.
The contour passes through points \(t_1,\ldots,t_n\), and the
matrix function is piece-wise constant, constructed from the
monodromy matrices \(M_1,\ldots, M_n.\) This factorization problem
can be reduced to a singular linear matrix integral equation on
the contour. On the second step the solution of the factorization
problem undergoes certain additional transformations in order to
obtain the solution of the original Riemann-Hilbert problem. The
second step is not completely painless: there are reefs here.  In
fact, not long ago A. A. Bolibruch has discovered examples of data
for which the Riemann-Hilbert problem has no solution (see
\cite{AnBo}, \cite{Bol1}, \cite{Bol2}, \cite{Bol3}).

 In order to apply the
Riemann-Hilbert problem for the construction of the isomonodromic deformation,
first we have  to determine for the equation \eqref{FDES81} (whose coefficients
-- matrices \( Q_1^0,\ldots,Q_n^0\) -- are the initial values in the Cauchy
problem for the Schlesinger system)
 the
monodromy matrices \(M_1^0,\ldots, M_n^0,\) corresponding to the loops
\(\gamma_1(\boldsymbol{t}^0),\dots , \gamma_n(\boldsymbol{t}^0).\) Then we
have to solve the Riemann-Hilbert problem for the following data:
 an arbitrary point
\(\boldsymbol{t}=(t_1,\ldots,t_n)\in\mathbb{C}_*^n\),  the
appropriate loops
\(\gamma_1(\boldsymbol{t}),\dots , \gamma_n(\boldsymbol{t})\),
and the matrices \(M_1^0,\ldots, M_n^0.\) As a result we shall
obtain the matrices
\(Q_1(\boldsymbol{t}),\ldots,Q_n(\boldsymbol{t}).\) Repeating this
procedure for every \(\boldsymbol{t},\) (and, therefore, solving
the infinite family of the Riemann-Hilbert problems) we construct
matrix {\em functions}
\(Q_1(\boldsymbol{t}),\ldots,Q_n(\boldsymbol{t}).\) If we succeed
in performing this construction in such a way that these matrix
functions are holomorphic then we obtain the desired isomonodromic
deformation. Finally, if the initial data are non-resonant, we
obtain the solution of the Cauchy problem for the Schlesinger
system. Thus solving the Cauchy problem for a {\em non-linear}
system is reduced to solving a family of {\em linear} problems.

However, solving the Schlesinger system in the non-resonant case is not the
ultimate goal of the present manuscript. We would like to consider the case
which is, in a certain sense, just the opposite. More precisely, we consider
the initial values \( Q_1^0,\ldots,Q_n^0\) such that  fundamental solutions of
the Fuchsian system \eqref{FDES81} are rational matrix functions of \(x\) in
general position. Of course, a rational function is uni-valued in the complex
plane, hence the monodromy in this case is trivial: \(M_1=\ldots=M_n=I.\) In
this case the spectrum of each matrix \( Q_j^0\) is either \(\{0,1\}\) or
\(\{0,-1\}\) and \(\mathrm{rank}(Q_j^0)=1.\) In particular, all the matrices
 \( Q_1^0,\ldots,Q_n^0\) are resonant. Therefore, the considerations above,
 concerning the isomonodromic deformations, are not applicable in this case.
Moreover, any family of Fuchsian systems with rational solutions is
isomonodromic (the monodromies are trivial). Thus there are "too many"
isomonodromic deformations here, and only one of them leads to the solution of
the Schlesinger system. This is the isoprincipal deformation.

Let us recall Theorem~\ref{SRQI} and, in particular, representation
\eqref{RSNRQ}. If system \eqref{FDE} has a generic rational solution \(Y(x)\)
then in \eqref{RSNRQ} we have \(A_j=0\),
 and hence \(Y(x)\)
is locally representable in a neighborhood of \(x=t_j\) as
\begin{equation}
\label{faczj}
 Y(x)=H_j(x)(x-t_j)^{Z_j},
\end{equation}
with \(H_j(x)\)  holomorphic and invertible there, and
\begin{equation}
\label{projz}
Z_j^2=\pm Z_j.
\end{equation}
  Furthermore, under condition
\eqref{projz} the factor \((x-t_j)^{Z_j}\) has the following simple form:
\begin{equation}
\label{RaPF}
(x-t_j)^{Z_j}=\left\{\begin{array}{l@{  \text{ if }  }l} I-Z_j+(x-t_j)Z_j,& Z_j^2=Z_j,\\
I+Z_j-\dfrac{1}{x-t_j}Z_j,& Z_j^2=-Z_j.
\end{array}\right.
\end{equation}
Thus the factorization \eqref{faczj} takes the form
\begin{eqnarray}
\label{FARS1}
Y(x)=H_j(x)\Big( I-Z_j+(x-t_j)Z_j\Big),&\text{ if }& Z_j^2=Z_j,\\
\label{FARS2}
 Y(x)=H_j(x)\Big(I+Z_j-\dfrac{1}{x-t_j}Z_j\Big),&\text{ if }&
Z_j^2=-Z_j.
\end{eqnarray}
Factorizations of the form \eqref{FARS1}, \eqref{FARS2} are traditional in the
theory of rational matrix functions. The factors \eqref{RaPF} are just the
principal factors of the rational solution \(Y(x)\) at the singular point
\(t_j\). If \(Y(x,\boldsymbol{t})\) is a family of rational solutions for a
family of Fuchsian systems which gives an isoprincipal deformation then for
every \(\boldsymbol{t}=(t_1,\ldots,t_n)\) at each singular point \(x=t_j\)
\(Y(x,\boldsymbol{t})\) must admit the factorization
\begin{eqnarray}
\label{FARS1t}
Y(x,\boldsymbol{t})=H_j(x,\boldsymbol{t})\Big( I-Z_j+(x-t_j)Z_j\Big),&\text{ if }& Z_j^2=Z_j,\\
\label{FARS2t}
 Y(x,\boldsymbol{t})=H_j(x,\boldsymbol{t})\Big(I+Z_j-\dfrac{1}{x-t_j}Z_j\Big),&\text{ if }&
Z_j^2=-Z_j,
\end{eqnarray}
where \(Z_j\) do not depend on \(\boldsymbol{t},\) and
\(H_j(x,\boldsymbol{t})\) are holomorphic and invertible at \(x=t_j.\) The
realization theory for rational matrix functions allows to construct families
of Fuchsian systems whose solutions \(Y(x,\boldsymbol{t})\) admit the
factorizations \eqref{FARS1t}, \eqref{FARS2t}. Development of the realization
theory appropriate for this goal originated with L.A. Sakhnovich \cite{Sakhn}.
This topic was further developed in \cite{GKRM}. In \cite{Kats1} and
\cite{Kats2} this theory was adapted for applications concerning the
Schlesinger system. The forthcoming second part of the present manuscript will
be dedicated to detailed construction of rational solutions of the Schlesinger
system, based on the realization theory for rational matrix functions.

\newpage
\vspace{2.6ex}%
\renewcommand{\thesection}{\Alph{section}}
\setcounter{section}{0}%
\centerline{\( \boldsymbol{\mathcal{APPENDIX}}\)}

\vspace{1.8ex}%
\noindent
\section{THE REPRESENTATION OF A SOLUTION
 Of A FUCHSIAN
EQUATION IN A NEIGHBORHOOD OF A SINGULAR POINT IN THE NON-RESONANT
CASE.}
\renewcommand{\theequation}{\Alph{section}.\arabic{equation}}

\noindent%
\textsf{PROOF OF THE PROPOSITION \ref{ELS}}. Choose a positive \(\delta\) so
small that the closure of the polydisk
\(\mathsf{D}(\boldsymbol{t}^0, \delta)\) is contained in the set
\(\boldsymbol{\mathcal{D}}\), where the matrix functions
\(Q_p(\boldsymbol{t}), p=1, \dots , n,\) are defined and holomorphic. We
also impose the condition \(\delta<\frac{1}{2}\rho_j(\boldsymbol{t}^0)\). So,
if \(|t_j-t_j^0|<\delta,\) and \(|t_p-t_p^0|<\delta, p\not=j\), then
\(|t_j-t_p|>\rho_j(\boldsymbol{t}^0)\). Since the matrix
\(Q_j(\boldsymbol{t}^0)\) is non-resonant, choosing, if necessary, a smaller
\(\delta\), we can also assume that the matrix \(Q_j(\boldsymbol{t})\) is
non-resonant for all
\(\boldsymbol{t}\in\mathsf{D}(\boldsymbol{t}^0, \delta)\). Moreover, if
\(\delta>0\) is small enough, then the differences of eigenvalues
\(\lambda_p(\boldsymbol{t})-\lambda_q(\boldsymbol{t})\) of the matrix
\(Q_j(\boldsymbol{t})\) are separated from non-zero integers:
\begin{equation}%
\label{Sep}%
|\lambda_p(\boldsymbol{t})-\lambda_q(\boldsymbol{t})-m|\geq\varepsilon>0 \quad
\forall \ p,q=1, \dots , k,\quad\forall\ m\in\mathbb{Z}\setminus\{0\},
\ \forall \ \boldsymbol{t}\in\mathsf{D}(\boldsymbol{t}^0, \delta) ,
\end{equation}%
where \(\varepsilon>0\) does not depend on
\(\boldsymbol{t}\in\mathsf{D}(\boldsymbol{t}^0,\delta)\),
\(m\in\mathbb{Z}\setminus\{0\}\).

Changing variable
\begin{equation}%
\label{ChVar}%
x\to x+t_j,
\end{equation}%
we reduce the differential equation (\ref{FDES}) to the form
\begin{equation}%
\label{RFDES}%
\frac{dZ(x,\boldsymbol{t})}{d x}=
\bigg(\frac{Q_j(\boldsymbol{t})}{x}+\Phi(x,\boldsymbol{t})\bigg)Z(x,\boldsymbol{t}) ,
\end{equation}%
where
\begin{equation}%
\label{NF}%
Z(x,\boldsymbol{t})=Y(x+t_j,\boldsymbol{t}),\quad
\Phi(x,\boldsymbol{t})=\sum_{\substack{1\leq p\leq n\\
p\not=j}}\frac{Q_p(\boldsymbol{t})}{x+t_j-t_p} .
\end{equation}
In view of the relation \(|t_j-t_p|>\rho_j(\boldsymbol{t}^0)\), the function
\(\Phi(x,\boldsymbol{t})\) admits expansion of the form
\begin{equation}%
\label{EPhi}%
\Phi(x,\boldsymbol{t})=\sum\limits_{0\leq r<\infty}\Phi_r(\boldsymbol{t})x^r,
\quad x\in\mathbb{C},\ |x|<\rho_j(\boldsymbol{t}^0) ,
\end{equation}%
where \(\Phi_r(\boldsymbol{t})\) are \(k\times k\) matrix functions, which are
holomorphic with respect to \(\boldsymbol{t}\) for
\(\boldsymbol{t}\in\mathsf{D}(\boldsymbol{t}^0, \delta)\) and admit the
estimate
\begin{equation}%
\label{EFC}%
\big\|\Phi_r(\boldsymbol{t})\big\|\leq
C_1\cdot(\rho_j(\boldsymbol{t}^0))^{-r}, \quad
\boldsymbol{t}\in\mathsf{D}(\boldsymbol{t}^0, \delta),\ \ 0\leq r<\infty,
\end{equation}%
where \(C_1<\infty\) does not depend on \(\boldsymbol{t}\) and \(r\).

We seek a normalized solution of (\ref{RFDES}), corresponding to the singular
point \(x=0\), in the form
\begin{equation}%
\label{Psi}%
Z(x)=\Psi(x, \boldsymbol{t})x^{Q_j(\boldsymbol{t})},
\end{equation}%
where \(\Psi(x, \boldsymbol{t})\) is  power series in \(x\):
\begin{equation}%
\label{PS}%
\Psi(x, \boldsymbol{t})=\sum\limits_{0\leq
r<\infty}\Psi_r(\boldsymbol{t})x^r,\quad
\text{with}\quad\Psi_0(\boldsymbol{t})=I .
\end{equation}
Substituting the expressions (\ref{Psi}) - (\ref{PS}) and (\ref{EPhi}) into
(\ref{RFDES}) and using the differentiation rule
\begin{equation}%
\label{DiRu}%
\frac{d}{dx} x^{Q_j(\boldsymbol{t})}=
\frac{Q_j(\boldsymbol{t})}{x} x^{Q_j(\boldsymbol{t})} ,
\end{equation}%
we obtain for \(r=0, 1, 2, \dots\)  the equalities
\begin{equation}%
\label{RecS}%
(r+1)\Psi_{r+1}(\boldsymbol{t})+\Psi_{r+1}(\boldsymbol{t})Q_j(\boldsymbol{t})
-Q_j(\boldsymbol{t})\Psi_{r+1}(\boldsymbol{t})=\sum_{\substack{l+m=r\\l\geq
0, m\geq 0}}\Phi_l(\boldsymbol{t})\Psi_m(\boldsymbol{t}),
\end{equation}%
relating the given sequence \(\{\Phi_r(\boldsymbol{t})\}_{0\leq r<\infty}\)
and the sequence \(\{\Psi_r(\boldsymbol{t})\}_{0\leq r<\infty}\), which has to
be found. The relations (\ref{RecS}) express the fact that
(\ref{Psi}) - (\ref{PS}) is a \textsl{formal solution} of the differential
equation (\ref{FDES}). If the series (\ref{PS}) converge in a disk
\(\{x\in\mathbb{C}: |x|<\rho\}\), with \(\rho>0 ,\) then (\ref{Psi}), where
\(\Psi(x, \boldsymbol{t})\) is the sum of these series, is an actual solution
of (\ref{FDES}) in the punctured disk
\(\{x\in\mathbb{C}: |x|<\rho\}\setminus\{0\}.\)

Since the constant term \(\Psi_0(\boldsymbol{t})\) is given:
\(\Psi_0(\boldsymbol{t})=I\), the equations (\ref{RecS}) can be considered as
a recursive system for successive determining of \(\Psi_1(\boldsymbol{t}),\)
\(\Psi_2(\boldsymbol{t})\), etc. The \(r\)-th equation of the system
(\ref{RecS}) is of the form
\begin{equation}%
\label{CommE}%
\lambda X - [Q,X]=Y ,
\end{equation}%
where \([Q,X]=QX-XQ\)  ---  the commutator of the matruces \(Q\) and \(X\).
In (\ref{RecS}),
\begin{equation}%
\label{Decod}%
 \lambda=r+1,\quad  Q=Q_j(\boldsymbol{t}),\quad X=\Psi_{r+1}(\boldsymbol{t}),\quad
Y=\sum_{\substack{l+m=r\\l\geq 0, m\geq
0}}\Phi_l(\boldsymbol{t})\Psi_m(\boldsymbol{t}) .
\end{equation}%

Now we interrupt the proof of \textsf{Proposition \ref{ELS}} to discuss the
matrix equation (\ref{CommE}).

\begin{lemma}%
\label{SpComm}%
Let \(\mathfrak{M}_k\) be the set of all \(k\times k\) matrices with complex
entries. Given a matrix \(Q\in\mathfrak{M}_k\), let us associate with this
matrix the operator \(\textup{ad}_{ Q}:\)
\begin{equation}
\label{DefAd}%
\textup{ad}_{ Q}X=QX-XQ,\quad
\textup{ad}_{ Q}: \mathfrak{M}_k\to\mathfrak{M}_k .
\end{equation}
Let \(\lambda_1(Q), \dots ,\lambda_k(Q)\) be the set of all eigenvalues of
the matrix \(Q\) (enumerated with multiplicities). Then the set of all
eigenvalues (enumerated with multiplicities) of the operator
\(\textup{ad}_{ Q}\) is the set of all the differences
\(\{\lambda_p(Q)-\lambda_q(Q)\}_{1\leq p, q\leq k} .\)
\end{lemma}%

\noindent%
\textsf{PROOF}. Assume first that all eigenvalues of the matrix \(Q\) are
pairwise different. Let \(u_p\) and \(v_q\) be, respectively, eigenvector
columns and eigenvector rows of the matrix \(Q\):
\begin{equation*}%
\label{EVCR}%
 Qu_p=\lambda_pu_p ,\quad v_qQ=\lambda_qv_q , \quad 1\leq p,q\leq k .
\end{equation*}%
Then the matrices \(X_{p,q}\stackrel{\scriptscriptstyle
\textup{def}}{=}u_p\otimes v_q\) are eigenvectors of the operator
\(\textup{ad}_{ Q}\):
\begin{equation*}%
\label{EVAd}%
\textup{ad}_{ Q}X_{p, q}=(\lambda_p-\lambda_q)X_{p, q},\quad 1\leq
p, q\leq k .
\end{equation*}%
Thus, the numbers \(\lambda_p(Q)-\lambda_q(Q)\) are eigenvalues of
the operator \(\textup{ad}_{ Q}\). Since the eigenvalues of the
matrix \(Q\) are pairwise different, the vectors \(\{u_p\}_{1\leq
p\leq k}\) and \(\{v_q\}_{1\leq q\leq k}\) form bases in the
spaces of all \(k\)-vector columns and \(k\)-vector rows
respectively. Since dimension of the space \(\mathfrak{M}_k\)  is
equal to \(k^2\), the set of matrices \(\{X_{p, q}\}_{1\leq
p, q\leq k}\) is a basis of the space \(\mathfrak{M}_k\).
Therefore, the set \(\{\lambda_p(Q)-\lambda_q(Q)\}_{1\leq
p, q\leq k}\) is the set of \textsl{all} eigenvalues of the
operator \(\textup{ad}_{ Q}\) (enumerated with multiplicities).
This result can be extended to general matrices \(Q\), i.e. to
matrices  whose eigenvalues are not necessarily pairwise
different. Such an extension can be done using the approximation
reasoning. \hfill Q.E.D.\\[1ex]

\noindent From the above result it follows that if \(\lambda\)
satisfies the condition \(\lambda\not=
\lambda_p(Q)-\lambda_q(Q),\\  1\leq p ,q\leq k\), then the
equation (\(Q, X, Y\) are \(k\times k\) matrices)
\begin{equation}%
\label{AdE}%
\lambda X-\textup{ad}_QX=Y
\end{equation}%
is solvable with respect to \(X\) for every \(Y\). However, for what follows
one needs not only the solvability of this equation but also an estimate for
its solution.

\begin{lemma}%
\label{SWE}%
Let \(\mathcal{L}\) be a \(K\)-dimensional vector space equipped with a norm
\(\| . \|_{\mathcal{L}}\), and let \(\mathcal{A}\) be a linear operator
acting in \(\mathcal{L}\). Assume that the eigenvalues
\(\lambda_m(\mathcal{A})\) of the operator \(\mathcal{A}\) are separated from
zero:
\begin{equation}%
\label{BAFZ}%
|\lambda_m(\mathcal{A})|\geq \varepsilon ,\quad 1\leq m\leq K .
\end{equation}%
Then the inverse operator \(\mathcal{A}^{-1}\) admits the estimate from above:
\begin{equation}%
\label{EFA}%
\|\mathcal{A}^{-1}\|\leq 2^{K}\|\mathcal{A}\|^{K-1}\varepsilon^{-K} ,
\end{equation}%
where \(\|\mathcal{A}\|=\sup
\|\mathcal{A}x\|_{\mathcal{L}}, \|\mathcal{A}^{-1}\|=
\sup\|\mathcal{A}^{-1}x\|_{\mathcal{L}}\), and \(\sup\)'s are taken over all
\(x\in\mathcal{L}\) such that \(\|x\|_{\mathcal{L}}\leq 1\).
\end{lemma}%

\noindent%
\textsf{PROOF.} %
Let \(\chi(\zeta)=\det(\zeta \mathcal{I}-\mathcal{A})\) be the characteristic
polynomial of the operator \(\mathcal{A}\),
\begin{equation*}%
\label{CharPol}%
\chi(\zeta)=\zeta^K-\sigma_1 \zeta^{K-1}+ \dots +
(-1)^{K-1}\sigma_{K-1} \zeta+(-1)^{K}\sigma_K ,
\end{equation*}%
where \(\sigma_l ,1\leq l\leq K,\) denotes  the  \(l\)-th elementary
symmetric polynomial of in \(\lambda_1(\mathcal{A}) , \dots ,\\
 \lambda_K(\mathcal{A})\). According to Hamilton-Cayley theorem,
\(\chi(\mathcal{A})=0\), hence
\begin{equation}%
\label{CarPol}%
\mathcal{A}^{-1}=(-1)^{K-1}\frac{1}{\sigma_K} \mathcal{A}^{K-1}+
(-1)^{K-2}\frac{\sigma_1}{\sigma_K} \mathcal{A}^{K-2}+ \cdots +
\frac{\sigma_{K-1}}{\sigma_K} \mathcal{I} .
\end{equation}%
>From (\ref{CarPol}) it follows that
\begin{equation*}%
\label{EstInv}%
\|\mathcal{A}^{-1}\|\leq
\left|\frac{\sigma_0}{\sigma_K}\right| \|\mathcal{A}\|^{K-1}+
\left|\frac{\sigma_1}{\sigma_K}\right| \|\mathcal{A}\|^{K-2}+ \cdots +
\left|\frac{\sigma_{K-1}}{\sigma_K}\right| .\quad
(\sigma_0\stackrel{\scriptscriptstyle \textup{def}}{=}1)
\end{equation*}%
It is clear that the fractions \(\displaystyle \frac{\sigma_{l}}{\sigma_K}\)
are the elementary symmetric polynomials in the inverse eigenvalues
\(\displaystyle \frac{1}{\lambda_1(\mathcal{A})},
 \dots ,\frac{1}{\lambda_K(\mathcal{A})} .\) From (\ref{BAFZ}) we conclude
that
\begin{equation*}%
\label{EKCP}%
\left|\frac{\sigma_{l}}{\sigma_K}\right|\leq C_{K}^{ l} \varepsilon^{l-K},\
\ \textup{where }C_{K}^{ l}\text{ are binomial coefficients: }
C_{K}^{ l}=\frac{K!}{l!(K-l)!} .
\end{equation*}%
Therefore,
\begin{equation*}%
\label{EstInvA}%
\|\mathcal{A}^{-1}\|\leq C_{K}^{ 0} \varepsilon^{-K}\|\mathcal{A}\|^{K-1}+
C_{K}^{ 1}\varepsilon^{1-K}\|\mathcal{A}\|^{K-2}+ \cdots +
 C_{K}^{K-1}\varepsilon^{-1}\|\mathcal{A}\|^{0} .
\end{equation*}%
Using the binomial formula, we rewrite the previous inequality as
\begin{equation}%
\label{EFAU}%
\|\mathcal{A}^{-1}\|\leq
\|\mathcal{A}\|^{-1}\left(\left(1+\varepsilon^{-1}\|\mathcal{A}\|\right)^K
-1\right) ,
\end{equation}%
Since \(|\lambda_m(\mathcal{A})|\leq \|\mathcal{A}\|\) for any eigenvalue
\(\lambda_m(\mathcal{A})\), inequality (\ref{BAFZ}) implies that \(1\leq
\varepsilon^{-1}\|\mathcal{A}\|\), and hence \(1+
\varepsilon^{-1}\|\mathcal{A}\|\leq 2\varepsilon^{-1}\|\mathcal{A}\|\).
 Weakening the inequality (\ref{EFAU}), we come to the more rough, but more
 simple inequality (\ref{EFA}).  \hfill Q.E.D.

\begin{lemma}%
\label{SWEM}%
Let \(Q\) be a \(k\times k\) matrix, with eigenvalues
\(\lambda_1, \dots ,\lambda_k\), enumerated with multiplicities. Let
\(\lambda\) be a complex number, and let for some \(\varepsilon>0\) the
inequalities
\begin{equation}%
\label{SepSp}%
|\lambda-(\lambda_p-\lambda_q)|>\varepsilon, \quad 1\leq p, q\leq k ,
\end{equation}%
hold true. Then for every \(Y\), the equation (\ref{AdE}) (or, what is the
same, the equation (\ref{CommE})) is solvable with respect to \(X\), and the
estimate
\begin{equation}%
\label{EFSol} \|X\|\leq 2^{k^2}\varepsilon^{-k^2}
(|\lambda|+2\|Q\|)^{k^2-1} \|Y\| ,
\end{equation}%
holds, where \(\| . \|\) is an arbitrary norm on \(\mathfrak{M}_k\)
with the property: \(\|A B\|\leq \|A\| \|B\|\) for every
\(A, B\in\mathfrak{M}_k\).
\end{lemma}%

\noindent \textsf{PROOF}. %
We use Lemma \ref{SWE} for the operator \(\mathcal{A}=\lambda
\mathcal{I}-\textup{ad}_{Q},\) acting on the space \(\mathfrak{M}_k.\) It is
clear that \(\|\mathcal{A}\|\leq |\lambda|+2\|Q\|\). According to Lemma
\ref{SpComm}, the set \(\{\lambda-(\lambda_p-\lambda_q)\}_{1\leq p, q\leq
k}\) is the set of all eigenvalues of \(\mathcal{A}\) (enumerated with
multiplicities). The inequality (\ref{SepSp}) provides an estimate of these
eigenvalues from below.
 Applying
Lemma \ref{SpComm} to this operator \(\mathcal{A}\) (and taking
into account that \(\textup{dim} \mathfrak{M}_k=k^2\)), we obtain
the estimate (\ref{EFSol}). \hfill Q.E.D.\\[1ex]

\noindent The following lemma gives an estimate, which is more
precise than (\ref{EFSol}) for large \(\lambda\).

\begin{lemma}%
\label{SWER}%
Let \(Q\) be a \(k\times k\) matrix, with eigenvalues
\(\lambda_1, \dots ,\lambda_k\), enumerated with multiplicities. Let
\(\lambda\) be a complex number which is different from the numbers
\(\lambda_p-\lambda_q, 1\leq p, q\leq k\),
 and let the estimates for \(\|Q\|\) from above and for
 \(|\lambda-(\lambda_p-\lambda_q)|\) from below be given:
 \begin{equation}%
 \label{EFAB}%
 |\lambda-(\lambda_p-\lambda_q)|\geq \varepsilon,\quad 1\leq p, q\leq k, \qquad
 \|Q\|\leq \mu , \quad (\varepsilon>0,\  \mu<\infty) .
 \end{equation}%
 Then for every \(Y\in\mathfrak{M}_k\), the equation (\ref{AdE}) (or, what is the same, the
equation (\ref{CommE})) is solvable with respect to \(X\in\mathfrak{M}_k\),
and the estimate
\begin{equation}%
\label{EFSolR} \|X\|\leq
\frac{C(\varepsilon, \mu, k)}{1+|\lambda|} \|Y\| ,
\end{equation}%
holds, where \(\| . \|\) is an arbitrary norm on \(\mathfrak{M}_k\)
with the property: \(\|A B\|\leq \|A\| \|B\|\) for every
\(A, B\in\mathfrak{M}_k\), and \(C(\varepsilon, \mu, k)<\infty\) is a
constant depending only on \(\varepsilon, \mu\) and \(k\).
\end{lemma}%

\noindent%
\textsf{PROOF}. %
 If \(\lambda X-\textup{ad}_QX=Y\), then
 \(|\lambda|||X\|\leq \|\textup{ad}_QX\|+\|Y\|\leq2\|Q\|\|X\|+
 \|Y\|\leq 2\mu\|X\|+\|Y\|\).
 Thus,
\begin{equation}%
\label{IFBL}%
\|X\|\leq \frac{1}{|\lambda|-2\mu}  \|Y\|,\quad \text{for}\ \
|\lambda|>2\mu .
\end{equation}%
On the other hand, if we assume that \(|\lambda|\leq 3\mu\) and replace in
 the inequality (\ref{EFSol}) \(|\lambda|\) and \(\|Q\|\) by the
larger values \(3\mu\) and \(\mu\) respectively, we come to the inequality
\begin{equation}%
\label{WeE}%
\|X\|\leq 10^{k^2}\varepsilon^{-k^2}\mu^{k^2-1} \|Y\| , \quad \text{for
}|\lambda|\leq 3\mu .
\end{equation}%
Unifying the inequalities  (\ref{IFBL})and (\ref{WeE}), we come to the
inequality (\ref{EFSolR}) with a suitable constant
\(C(\varepsilon, \mu, k)\) (which, of course, can be found explicitly).%
\hfill Q.E.D.\\[1ex]

\noindent Now, after we have investigated the matrix equation
(\ref{CommE}), we resume the proof of \textsf{Proposition
\ref{ELS}}. We apply Lemma \ref{SWEM} to the equality (\ref{RecS})
considered as an equation of the form (\ref{CommE}), with
\(X, Y, Q, \lambda\) defined in (\ref{Decod}). The number
\(\max\limits_{\boldsymbol{t}\in
\mathsf{D}(\boldsymbol{t}^0)}\|Q_j(\boldsymbol{t})\|\) serves as
\(\mu\), the number \(\varepsilon\) is taken from (\ref{Sep}).
 Substituting these
\(X, \lambda\) and the estimate
\begin{equation}%
\label{ERHS}%
\|Y\|\leq \sum_{\substack{0\leq m\leq
r}}\|\Phi_{r-m}(\boldsymbol{t})\| \|\Psi_m(\boldsymbol{t})\|
\end{equation}%
 for \(\|Y\|\) into (\ref{EFSolR}), we obtain, for \(\boldsymbol{t}\in\mathsf{D}
 (\boldsymbol{t}^0,\delta),\) \(
r=0, 1, 2, \dots\) the estimate
\begin{equation}%
\label{RecEst}%
\|\Psi_{r+1}(\boldsymbol{t})\|\leq
\frac{C(\varepsilon, \mu, k)}{r+2} \sum_{\substack{0\leq m\leq
r}}\|\Phi_{r-m}(\boldsymbol{t})\| \|\Psi_m(\boldsymbol{t})\| .
\end{equation}%
Substituting the upper estimate (\ref{EFC}) for
\(\|\Phi_{r-m}(\boldsymbol{t})\|\) into (\ref{RecEst}), we obtain:
\begin{equation}%
\label{RESt}%
\|\Psi_{r+1}(\boldsymbol{t})\| \rho^{r+1}\leq \frac{C_2}{r+2}
\sum\limits_{0\leq m\leq r} \|\Psi_{m}(\boldsymbol{t})\| \rho^{m} , \quad
\forall  \boldsymbol{t}\in\mathsf{D}(\boldsymbol{t}^0,\delta),\ \
r=0, 1, 2, \dots ,
\end{equation}%
where
\begin{equation}%
\label{CExp}%
C_2=\rho\cdot  C_1\cdot C(\varepsilon, \mu, k),
\end{equation}%
\(\rho=\rho_j(\boldsymbol{t}^0)\) is defined in (\ref{DefRho}), \(C_1\) and
\(C(\varepsilon, \mu, k)\) are the values from (\ref{EFC}) and
(\ref{EFSolR}) respectively. It is important to note that \textsf{the value
\(C_2<\infty\) does not depend on \(r\) and on
\(\boldsymbol{t}\in\mathsf{D}(\boldsymbol{t}^0)\).}

\begin{lemma}%
\label{PEst}%
Let \(\{a_r\}_{0\leq r<\infty}\) be a sequence of non-negative numbers,
satisfying the conditions
\begin{equation}%
\label{SuIn}%
a_{r+1}\leq\frac{d}{r+2}\sum\limits_{0\leq m\leq r}a_m, \quad
r=0, 1, 2, \dots ;\quad a_0=1,
\end{equation}%
where \(d, 0<d<\infty,\) does not depend on \(r\).

Then the sequence \(\{a_r\}_{0\leq r<\infty}\) admits the upper estimate
\begin{equation}%
\label{FUE}%
a_r\leq(r+1)^{d-1},\quad r=0, 1, 2, \dots .
\end{equation}%
\end{lemma}%

\noindent%
 \textsf{PROOF.}\ %
For \(r=0\), the assertion (\ref{FUE}) is true: \(a_0=1\). Assume that the
inequalities (\ref{FUE}) hold for \(r=0, 1, \dots , m\). Then, according
to (\ref{FUE}), the inequality
\begin{equation}%
\label{Indu}%
a_{m+1}\leq\frac{d}{m+2}\sum\limits_{0\leq r\leq m}(r+1)^{d-1}
\end{equation}%
holds. If \(d\geq 1\), then \[(r+1)^{d-1}\leq
\int\limits_{r+1}^{r+2}t^{d-1}dt,\ \  %
\sum\limits_{0\leq r\leq
m}(r+1)^{d-1}\leq\int\limits_{1}^{m+2}t^{d-1}dt\leq\frac{(m+2)^d}{d} ,\] %
and the inequality (\ref{FUE}) holds for \(r=m+1\). If \(0<d<1\), then
\[(r+1)^{d-1}\leq
\int\limits_{r}^{r+1}t^{d-1}dt,\ \  %
\sum\limits_{0\leq r\leq
m}(r+1)^{d-1}\leq\int\limits_{0}^{m+1}t^{d-1}dt\leq\frac{(m+2)^d}{d} ,\] %
and the inequality (\ref{FUE}) holds for \(r=m+1\) as well. By
induction, the inequality (\ref{FUE}) holds for all \(r\). \hfill
Q.E.D.\\[1ex]

\noindent Applying Lemma \ref{PEst} to (\ref{RESt}), with
\(a_r=\|\Psi_{r}(\boldsymbol{t})\| \rho^{r}\), we obtain the
estimate
\begin{equation}%
\label{FESTI}%
\|\Psi_{r}(\boldsymbol{t})\|\leq (r+1)^d \rho^{ -r}, \quad
r=0, 1, 2, \dots ,\quad \forall
 \boldsymbol{t}\in\mathsf{D}(\boldsymbol{t}^0, \delta),
\end{equation}%
where \(\rho=\rho_j(\boldsymbol{t}^0)\) is defined in (\ref{DefRho}),
\(d=C_2\) from (\ref{CExp}), \(d<\infty\) does not depend on \(r\) and on
\(\boldsymbol{t}\in\mathsf{D}(\boldsymbol{t}^0, \delta)\).

>From the estimate (\ref{FESTI}) it follows, that the power series (\ref{PS})
converge locally uniformly in the open circle \(\{x: |x|<\rho\}\), and the
convergence is uniform with respect to the parameter
\(\boldsymbol{t}\in\mathsf{D}(\boldsymbol{t}^0, \delta)\). Therefore, the sum
\(\Psi(x, \boldsymbol{t})\) of these power series is a function, holomorphic
with respect to \(x,\boldsymbol{t}\) for
\(\{x: |x|<\rho_j(\boldsymbol{t}^0)\}\),
\(\boldsymbol{t}\in\mathsf{D}(\boldsymbol{t}^0, \delta)\). Let us prove that
the function \(\Psi^{-1}(x, \boldsymbol{t})\) is holomorphic for these %
\footnote{Of course, in view of the condition \(\Psi_0(\boldsymbol{t})=I\),
the function \(\Psi^{-1}(x, \boldsymbol{t})\) is holomorphic for \(x\) which
are small enough. However, it is very easy to prove the holomorphy of the
function \(\Psi^{-1}\) in the whole disc \(\{x: |x|<\rho\}\). Such a proof
can be done "for free", using the previously established estimates.} %
\(x,\boldsymbol{t}\)  as well. Since \(\Psi_0(\boldsymbol{t})=I\), we can
consider the \textsl{formal} power series
\begin{equation}%
\label{PSIn}%
\Psi^{\text{inv}}(x, \boldsymbol{t})=\sum\limits_{0\leq
r<\infty}\Psi^{\text{inv}}_{ r}(\boldsymbol{t})x^r,\quad
\end{equation}%
which are \textsl{inverse} to the series (\ref{PS}), i.e.
\begin{equation}%
\label{FormInvS}%
\begin{array}{rcl}
\left(\sum\limits_{0\leq
r<\infty}\Psi_{ r}(\boldsymbol{t})x^r\right)&\cdot&\left(\sum\limits_{0\leq
r<\infty}\Psi^{\text{inv}}_{ r}(\boldsymbol{t})x^r\right)\\&&=\left(\sum\limits_{0\leq
r<\infty}\Psi^{\text{inv}}_{ r}(\boldsymbol{t})x^r\right)\cdot\left(\sum\limits_{0\leq
r<\infty}\Psi_{ r}(\boldsymbol{t})x^r\right)=I .\end{array}
\end{equation}%
It is clear that
\begin{equation}%
\label{CoTe}%
\Psi^{\text{inv}}_{ 0}(\boldsymbol{t})=I.
\end{equation}%
The expressions   (\ref{RFDES}), (\ref{Psi}), (\ref{PS}) and (\ref{FormInvS})
imply that the product
\begin{equation}%
\label{PsiInv}%
 Z^{\text{inv}}(x,\boldsymbol{t}) \stackrel{\textup{def}}{=}
x^{-Q_j(\boldsymbol{t})}\Psi^{\text{inv}}(x, \boldsymbol{t})
\end{equation}%
is a formal solution of the differential equation
\begin{equation}%
\label{RFDESInv}%
\frac{dZ^{\textup{inv}}}{d x}=
-Z^{\text{inv}}\bigg(\frac{Q_j(\boldsymbol{t})}{x}+\Phi(x,\boldsymbol{t})\bigg) ,
\end{equation}%
As before (see (\ref{RecS})), substituting the expression (\ref{PsiInv}) into
(\ref{RFDESInv}), we obtain the infinite recursive system
\begin{equation}%
\label{RecSInv}%
(r+1)\Psi^{\textup{inv}}_{r+1}(\boldsymbol{t})+
\Psi^{\textup{inv}}_{r+1}(\boldsymbol{t})Q_j(\boldsymbol{t})
-Q_j(\boldsymbol{t})\Psi^{\textup{inv}}_{r+1}(\boldsymbol{t})=
-\sum_{\substack{l+m=r\\l\geq 0, m\geq
0}}\Psi^{\textup{inv}}_{ m}(\boldsymbol{t})\Phi_l(\boldsymbol{t}) .
\end{equation}%
Like the \(r\)-th equation of the system (\ref{RecS}), the \(r\)-th equation
of the system (\ref{RecSInv}) is of the form (\ref{CommE}), with \textsl{the
same} \(\lambda\) and \(Q\) as in (\ref{RecS}):
\(\lambda=r+1, Q=Q_j(\boldsymbol{t})\). The only difference is that the
right-hand side of (\ref{CommE}) is \(Y=\sum\Phi_l(\boldsymbol{t})
\Psi_m(\boldsymbol{t})\) for (\ref{RecS}),
and  \(Y=-\sum\Psi^{\textup{inv}}_{ m}(\boldsymbol{t})
\Phi_l(\boldsymbol{t})\) for
(\ref{RecSInv}). Because \(\lambda\) and \(Q\) for both matrix equations,
corresponding to the systems (\ref{RecS}) and (\ref{RecSInv}), are the same,
the estimates (\ref{EFSolR}) for both equations hold with \textsf{the same}
constant \(C(\varepsilon,\mu,k)\) (Of course, \(\varepsilon,\mu,k\) are the
same). Instead of the estimate (\ref{ERHS}), the estimate
\begin{equation}%
\label{ERHSInv}%
 \|Y\|\leq \sum_{\substack{0\leq m\leq
r}}\|\Phi_{r-m}(\boldsymbol{t})\| \|\Psi^{\text{inv}}_{ m}(\boldsymbol{t})\|
\end{equation}%
for the right-hand side of (\ref{CommE}) should be used now. Finally, we come
to the recursive system of inequalities
\begin{equation}%
\label{RecEstInv}%
\|\Psi^{\text{inv}}_{r+1}(\boldsymbol{t})\|\leq
\frac{C(\varepsilon, \mu, k)}{r+2} \sum_{\substack{0\leq m\leq
r}}\|\Phi_{r-m}(\boldsymbol{t})\| \|\Psi^{\text{inv}}_m(\boldsymbol{t})\| ,
\quad \forall
 \boldsymbol{t}\in\mathsf{D}(\boldsymbol{t}^0,\delta).
\end{equation}%
with respect to  \(\|\Psi^{\text{inv}}_{m}(\boldsymbol{t})\|\). However, the
latter system of inequalities coincides with the system (\ref{RecEst}) of
inequalities with respect to \(\|\Psi_{m}(\boldsymbol{t})\|\). Since
\(\|\Psi^{\text{inv}}_{ 0}(\boldsymbol{t})\|=
\|\Psi_{0}(\boldsymbol{t})\|=\|I\|\), the estimates for
\(\|\Psi_{m}(\boldsymbol{t})\|\) and
\(\|\Psi^{\text{inv}}_{m}(\boldsymbol{t})\|\), which follow from
(\ref{RecEst}) and (\ref{RecEstInv}), coincide. Thus, the estimates
\begin{equation}%
\label{FESTIInv}%
\|\Psi^{\text{inv}}_{r}(\boldsymbol{t})\|\leq (r+1)^d \rho^{ -r}, \quad
r=0, 1, 2, \dots ,\quad \forall
 \boldsymbol{t}\in\mathsf{D}(\boldsymbol{t}^0, \delta),
\end{equation}%
hold with the same \(\rho\) and \(d\) that in (\ref{FESTI}). Therefore, the
series (\ref{PSIn}) converge locally uniformly in the disk
\(\{x: |x|<\rho_j(\boldsymbol{t})\}\) and represent the matrix
\textsl{function} \(\Psi^{\text{inv}}(x, \boldsymbol{t})\) which is
holomorphic with respect to \(x, \boldsymbol{t}\) for
\(x: |x|<\rho_j(\boldsymbol{t}^0)\),
\(\boldsymbol{t}\in\mathsf{D}(\boldsymbol{t}^0, \delta)\). The identities
(\ref{FormInvS}) for formal power series imply the identities
\begin{equation}%
\label{FunkInv}%
\Psi(x, \boldsymbol{t})\Psi^{\text{inv}}(x, \boldsymbol{t})=
\Psi^{\text{inv}}(x, \boldsymbol{t})\Psi(x, \boldsymbol{t})=I \quad \forall\
\  x:|x|<\rho_j(\boldsymbol{t}^0),
\boldsymbol{t}\in\mathsf{D}(\boldsymbol{t}^0, \delta)
\end{equation}%
for the appropriate holomorphic matrix functions. Thus, the values of the
matrix function \(\Psi(x, \boldsymbol{t})\)
 for \(x: |x|<\rho_j(\boldsymbol{t}^0)\),
\(\boldsymbol{t}\in\mathsf{D}(\boldsymbol{t}^0, \delta)\) are invertible, and
\(\Psi^{-1}(x, \boldsymbol{t})=\Psi^{\text{inv}}(x, \boldsymbol{t})\).

Performing the change of variables \(x\to x-t_j\), inverse to the change
(\ref{ChVar}), we obtain the representation (\ref{DT}) for the local solution
\(Y_j(x,\boldsymbol{t})\), with
\(H_j(x,\boldsymbol{t}))=\Psi(x-t_j,\boldsymbol{t})\). The properties of the
function \(H_j(x,\boldsymbol{t})\), stated in the formulation of Proposition
\ref{ELS}, express the properties of \(\Psi(x, \boldsymbol{t})\), as
established above. This completes the proof of Proposition \ref{ELS}.
 \hfill Q.E.D.

%

\newpage

\vspace{1.0cm}
\end{document}